\baselineskip=12pt
\font\blowup=cmbx10 scaled\magstep1
\def\adots{\mathinner{\mskip1mu\raise1pt\hbox{.}\mskip2mu\raise4pt\hbox{.}\mskip2mu\raise7pt\vbox{\kern7pt\hbox{.}}\mskip1mu}}

\def\End{\mathop{\rm End}\nolimits}
\def\Hom{\mathop{\rm Hom}\nolimits}

\def\tr{\mathop{\rm tr}\nolimits}

\mathchardef\bfplus="062B
\mathchardef\bfminus="067B
\font\title=cmbx10 scaled\magstep5
\font\chapter=cmbx10 scaled\magstep4
\font\section=cmbx10 scaled\magstep2

\def\~#1{{\accent"7E #1}}
\def\bull{\vrule height .9ex width .8ex depth -.1ex}
\def\sqr#1#2{{\vcenter{\hrule height.#2pt \hbox{\vrule width.#2pt height#1pt \kern#1pt \vrule width.#2pt}\hrule height.#2pt}}}

\def\hk#1#2{{\vcenter{\hrule height0.0pt \hbox{\vrule width0.0pt \kern#1pt \vrule width.#2pt height#1pt}\hrule height.#2pt}}}

\def\operp{\mathrel{\bigcirc \hskip -11.5pt \perp}}
\def\scoperp{\mathrel{\bigcirc \hskip -7.25pt \perp}}

\def\mapstoleft#1{\smash{\mathop{\longmapsto}\limits^{#1}}}

\centerline{\chapter Algebraic Classification of the Ricci}
\centerline{\chapter Curvature Tensor and Spinor}
\centerline{\chapter for Neutral Signature in Four Dimensions}
\vskip 24pt
\noindent {\section Abstract}\hfil\break
I apply the algebraic classification of self-adjoint endomorphisms of ${\bf R}^{2,2}$ provided by their Jordan canonical form to the Ricci curvature tensor of four-dimensional neutral manifolds and relate this classification to an algebraic classification of the Ricci curvature spinor. These results parallel similar results well known in four-dimensional Lorentzian geometry. The classification is summarized in Table 2 at the end of the paper.
\vskip 24pt
\noindent Peter R Law. 4 Mack Place, Monroe, NY 10950, USA. prldb@member.ams.org\hfil\break
\vskip 24pt
\vfill\eject
\noindent {\section 1. Introduction}
\vskip 12pt
Algebraic classifications of Weyl and Ricci curvature in the context of general relativity are well known. The algebraic classification of the Weyl curvature tensor was provided by Petrov, see e.g., Petrov (1969), and of the Weyl curvature spinor by Penrose, see e.g., Penrose \& Rindler (1986). Analogues of these classifications for neutral signature in four dimensions were presented in Law (1991) and (2006), respectively. Derdzinski (2000), \S 39, has also given a refinement of Petrov's classification, which he calls Petrov-Segre classes, based also on the eigenvalues and eigenvectors of the Weyl curvature endomorphisms. In the Lorentzian case, it coincides with the Penrose classification; in the neutral case it is of course very similar to, but does not exactly coincide with, the classification given in Law (2006) (Derdzinski's classification recognizes the Law 2006 types: $\{31\}$III; $\{1\overline{1}11\}$Ib; $\{22\}$Ia; $\{2\overline{2}\}$Ia; lumps $\{211\}$II and $\{1\overline{1}2\}$II into a single class; lumps $\{1111\}$Ia and $\{1\overline{1}1\overline{1}\}$Ia into a single class; and divides $\{4\}$II into two classes). The algebraic classification of the Ricci tensor in general relativity goes back at least to Churchill (1932); see Kramer et al. (1980), Ch. 5, for a modern treatment. Penrose provided an algebraic classification of the Ricci curvature spinor and related it to the classification of the Ricci curvature tensor, see, e.g., Penrose \& Rindler (1986). As for Weyl curvature, the classification of Ricci curvature for neutral signature is sufficiently different from the Lorentz-signature case that an explicit account is useful in the study of four-dimensional neutral geometry. The primary purpose of this paper is to provide that account.

By ${\bf R}^{p,q}$, I denote ${\bf R}^{p+q}$ equipped with the scalar product $s_{p,q}$ that makes the standard basis $\{e_1,\ldots,e_n\}$ pseudo-orthonormal ($\Psi$-ON) of signature $(p,q)$. A vector $v$ is called: {\sl time/space-like} if $s_{p,q}(v,v)$ is positive/negative; a {\sl unit vector} if $s_{p,q}(v,v) = \pm1$; and {\sl null} if $s_{p,q}(v,v) = 0$. By ${\bf R}^{2n}_{\rm hb}$, I denote ${\bf R}^{2n}$ equipped with the scalar product $s_{\rm hb}$ such that $s_{\rm hb}(e_i,e_{n+i}) = 1$, $i=1,\ldots,n$, are the only nontrivial scalar products amongst the elements of the standard basis. This space is called {\sl hyperbolic space} in Porteous (1981) but {\sl ultra-hyperbolic space} by some other authors; I shall follow Porteous. ${\bf R}^{2n}_{\rm hb}$ is isomorphic to ${\bf R}^{n,n}$. Any basis $\{E_1,\ldots,E_{2n}\}$ of ${\bf R}^{n,n}$ for which $s_{n,n}(E_i,E_{n+i}) = 1$, $i=1,\ldots,n$, are the only nontrivial scalar products is called a {\sl Witt basis} and models ${\bf R}^{n,n}$ as ${\bf R}^{2n}_{\rm hb}$. The group of orthogonal automorphisms of ${\bf R}^{p,q}$ is denoted by {\bf O(p,q)} and of ${\bf R}^{2n}_{\rm hb}$ by {\bf O(2n;hb)}. Of course, ${\bf O(2n;hb)} \cong {\bf O(n,n)}$; the matrix representation of {\bf O(n,n)} with respect to a Witt basis is of the same form as the matrix representation of {\bf O(2n;hb)} with respect to the standard basis. In particular, if an endomorphism $S$ of ${\bf R}^{n,n}$ has block matrix representation ${\underline S} = \left({A \atop B}{C \atop D}\right)$, $A$, $B$, $C$, and $D \in {\bf R}(n)$, with respect to a basis, then the adjoint $^*\! S$ has matrix representation with respect to the same basis:
$${^*\! {\underline S}} = \pmatrix {{^\tau\! A}&-{^\tau\! B}\cr -{^\tau\! C}&{^\tau\! D}\cr},\eqno(1.1)$$
when the basis is $\Psi$-ON; and
$${^*\! {\underline S}} = \pmatrix{{^\tau\! D}&{^\tau\! C}\cr {^\tau\! B}&{^\tau\! A}\cr},\eqno(1.2)$$
when the basis is Witt. 

The identity-connected components of {\bf O(p,q)} and {\bf O(2n;hb)} are denoted, respectively, by ${\bf SO^\bfplus(p,q)}$ and ${\bf SO^\bfplus(2n;hb)}$. Let $K_j$ be the orthogonal automorphism that leaves the standard basis unchanged except that $e_j \mapsto -e_j$ (i.e., hyperplane reflexion in $\langle e_j \rangle_{\bf R}^\perp$). I will also denote ${\bf SO^\bfplus(p,q)}$ by ${\bf O^\bfplus_\bfplus(p,q)}$ and define ${\bf O^\bfminus_\bfminus(p,q)} := {\bf SO(p,q)} \setminus {\bf SO^\bfplus(p,q)}$, ${\bf O^\bfplus_\bfminus(p,q)} := K_n{\bf SO^\bfplus(p,q)}$ and ${\bf O^\bfminus_\bfplus(p,q)} := K_1{\bf SO^\bfplus(p,q)}$, so that the four connected components of {\bf O(p,q)} are ${\bf O^\bfplus_\bfplus(p,q)}$, ${\bf O^\bfminus_\bfminus(p,q)}$, ${\bf O^\bfplus_\bfminus(p,q)}$ and ${\bf O^\bfminus_\bfplus(p,q)}$. There are various well known characterizations of the connected components of {\bf O(p,q)}. I simply record here the following useful facts. By the Cartan-Dieudonn\'e Theorem, each $L \in {\bf O(p,q)}$ can be expressed as a product of hyperplane reflexions: $L = K_{u_1} \circ \cdots \circ K_{u_m}$, for unit vectors $u_1,\ldots,u_m$, where $K_{u_j}$ is the hyperplane reflexion in $\langle u_j \rangle_{\bf R}^\perp$. Define $\tau(L) := \pm1$ according as there are an even/odd number of time-like unit vectors in the decomposition; define $\sigma(L) := \pm1$ according as there are an even/odd number of space-like unit vectors in the decomposition. It can be shown that $\bigl(\tau(L),\sigma(L)\bigr)$ is independent of the particular decomposition of $L$ into a product of hyperplane reflexions. Now let $\left({A \atop B}{C \atop D}\right)$ be the matrix representation of $L$ with respect to some given $\Psi$-ON basis, $A \in {\bf R}(p)$, $C \in {\bf R}(p,q)$, $B \in {\bf R}(q,p)$ and $D \in {\bf R}(q)$. Let 
$$t(L) := {\det(A) \over \vert \det(A)\vert} \hskip 1.25in s(L) := {\det(D) \over \vert \det(D) \vert}.\eqno(1.3)$$
Then it can also be shown that $\bigl(t(L),s(L)\bigr)$ is independent of the choice of $\Psi$-ON basis, that
$$\bigl(\tau(L),\sigma(L)\bigr) = \bigl(t(L),s(L)\bigr)\eqno(1.4)$$
and
$$\displaylines{{\bf O^\bfplus_\bfplus(p,q)} = \{\,L \in {\bf O(p,q)}:\bigl(t(L),s(L)\bigr) = (1,1)\,\}\cr
{\bf O^\bfplus_\bfminus(p,q)} = \{\,L \in {\bf O(p,q)}:\bigl(t(L),s(L)\bigr) = (1,-1)\,\}\cr
\noalign{\vskip -6pt}
\hfill\llap(1.5)\cr
\noalign{\vskip -6pt}
{\bf O^\bfminus_\bfplus(p,q)} = \{\,L \in {\bf O(p,q)}:\bigl((t(L),s(L)\bigr) = (-1,1)\,\}\cr
{\bf O^\bfminus_\bfminus(p,q)} = \{\,L \in {\bf O(p,q)}:\bigl(t(L),s(L)\bigr) = (-1,-1)\,\}.\cr}$$

I employ the abstract index notation of Penrose \& Rindler (1984); in particular, italic indices denote abstract indices, while upright bold indices are concrete, i.e., take numerical values; see \S 2. Henceforward, I shall refer to Penrose \& Rindler (1984) as PRI and Penrose \& Rindler (1986) as PRII. Common notations and terms from PRI \& II will also be employed, often without further explanation. My conventions, detailed in Law \& Matsushita (2008), Appendix One, follow those of PRI with one exception: I take the Ricci tensor to be the negative of PRI's definition and modify PRI (4.6.20--23) by replacing their Ricci tensor by its negative, thus leaving the definitions of $\Phi_{ABA'B'}$ and $\Lambda$ unchanged.

${\bf R}^{2,2}$ is the space of primary interest in this paper. In \S 2, I establish further notation for certain elementary notions, thereby illustrating the abstract index notation. In \S 3 I review spinor theory for neutral signature. \S 4 contains the classification of self-adjoint endomorphisms of ${\bf R}^{2,2}$ by means of the Jordan canonical form; \S 5 contains the algebraic classification of the spinor translates of traceless self-adjoint endomorphisms, and \S 6 relates these classifications.
\vskip 24pt
\noindent {\section 2. Abstract Indices; Linear Algebra}
\vskip 12pt
Let $V$ be a right {\bf K}-linear space, ${\bf K} = {\bf R}$ or {\bf C}, of finite dimension; though {\bf K} is commutative I employ notation appropriate for modules over non-commutative (division) rings, e.g., the quaternions {\bf H}. For endomorphisms $T$ of $V$, write $T(v\lambda) = T(v)\lambda$. With respect to a basis $\{v_1,\ldots,v_n\}$ of $V$, the matrix representation of $T$ is defined by
$$T(v_{\bf j}) =: v_{\bf i}T^{\bf i}{}_{\bf j} \hskip 1.25in v =: \sum_{\bf j}\,v_{\bf j}\lambda^{\bf j} \mapsto v_{\bf i}T^{\bf i}{}_{\bf j}\lambda^{\bf j}.\eqno(2.1)$$

The dual space $V_\bullet$ has a natural left {\bf K}-linear structure: $\bigl(\lambda.\phi\bigr)(v) = \lambda\phi(v)$. For $L \in \End(V_\bullet)$, write $L(\mu.\phi) = \mu L(\phi)$. For a basis $\{\phi^1,\ldots,\phi^n\}$ of $V_\bullet$, write the matrix representation of $L \in \End(V_\bullet)$ as
$$L(\phi^{\bf j}) =: L^{\bf j}{}_{\bf i}\phi^{\bf i} \hskip 1.25in \phi =: \sum_{\bf j}\,\mu_{\bf j}\phi^{\bf j} \mapsto \mu_{\bf j}L^{\bf j}{}_{\bf i}\phi^{\bf i}.\eqno(2.2)$$

There is a unique mapping
$$\Phi: V \otimes_{\bf K}V_\bullet \to \End(V)\qquad\hbox{satisfying}\qquad \Phi(v \otimes \phi)(w) = v\phi(w);\eqno(2.3{\rm a})$$
it is injective, and surjective in finite dimensions. If $\{v_1,\ldots,v_n\}$ and $\{\phi^1,\ldots,\phi^n\}$ are dual bases, then under $\Phi$:
$$\sum_{{\bf i},{\bf j}}\,T^{\bf i}{}_{\bf j}v_{\bf i} \otimes \phi^{\bf j}\ \leftrightarrow\ T.\eqno(2.3{\rm b})$$
The mapping
$$(v,\phi) \mapsto \phi(v),\eqno(2.4)$$
induces a unique {\bf K}-bilinear mapping
$${\cal C}:V \otimes V_\bullet \to {\bf K} \qquad\hbox{satisfying}\qquad v \otimes \phi \mapsto \phi(v),\eqno(2.5)$$
under which
$$T\ \leftrightarrow\ \sum_{{\bf i},{\bf j}}\,T^{\bf i}{}_{\bf j}v_{\bf i} \otimes \phi^{\bf j}\ \mapstoleft{\cal C}\ \sum_{\bf i}\, T^{\bf i}{}_{\bf i} = \tr(T).\eqno(2.6)$$
As the identification $\Phi$ permits one to view an endomorphism $T$ of $V$ as an element of $V \otimes V_\bullet$, using abstract indices write $T^a{}_b$ so that $T(v) = T^a{}_bv^b$. For an element $L \in \End(V_\bullet)$, write $L^b{}_a$ so that $L(\phi) = \phi_bL^b{}_a$. The abstract index notation assumes a notion of {\sl total reflexivity}, which entails in particular that the canonical inclusion of $V$ into its double dual is an isomorphism, see PRI; total reflexivity is automatic in finite dimensions.

The linear dual of $T \in \End(V)$ is defined by
$$T_\bullet(\phi) := \phi \circ T \hskip 1.25in \phi_b(T_\bullet)^b{}_a := \phi_bT^b{}_a.\eqno(2.7)$$
Thus, abstractly, $T_\bullet$ and $T$ may be viewed as the same object with two actions, consistent with the identifications $\End(V_\bullet) \cong V_\bullet \otimes V_{\bullet\bullet} \cong V_\bullet \otimes V \cong V \otimes V_\bullet \cong \End(V)$. As regards the matrix representation,
$$T^{\bf j}{}_{\bf k} = \phi^{\bf j}(v_{\bf i}T^{\bf i}{}_{\bf k}) = (\phi^{\bf j} \circ T)(v_{\bf k}) = \bigl(T_\bullet(\phi^{\bf j})\bigr)(v_{\bf k}) = \bigl((T_\bullet)^{\bf j}{}_{\bf i}\phi^{\bf i}\bigr)(v_{\bf k}) = (T_\bullet)^{\bf j}{}_{\bf k},\eqno(2.8)$$
but note that the indices play different roles in (2.1) and (2.2) so that in effect the matrix for the dual of $L$ is the transpose of the matrix of $L$ as usual.

Consider a {\bf K}-linear, nondegenerate, correlation $\xi:V \to V_\bullet$, i.e., an injective element of $\Hom_{\bf K}(V,V_\bullet)$. Define
$$g: V \times V \to {\bf K} \hskip 1.25in g(u,v) := \bigl(\xi(u)\bigr)(v).\eqno(2.9)$$
This definition is made so that $g$ is automatically {\bf K}-linear in its second argument when the assumption that $\xi$ is {\bf K}-linear is relaxed (e.g., to $\overline{\bf C}$-linearity), this choice being preferable for right linear spaces; $g$ is bilinear and represents $\xi$ as an element of $V_\bullet \otimes V_\bullet$.

Now $\xi$ is right-to-left linear: $\xi(v\lambda) = \lambda\xi(v)$. One therefore writes the matrix representation as
$$\xi(v_{\bf j}) =: \xi_{\bf ji}\phi^{\bf i} \hskip 1.25in v =: \sum_{\bf j}\,v_{\bf j}\lambda^{\bf j} \mapsto \lambda^{\bf j}\xi_{\bf ji}\phi^{\bf i}.\eqno(2.10)$$
Hence, abstractly, write $u \mapsto g_{ba}u^b = u^b\xi_{ba}$. The matrix representation of $g$ is
$$g_{\bf ij} := g(v_{\bf i},v_{\bf j}) = \bigl(\xi(v_{\bf i})\bigr)(v_{\bf j}) = (\xi_{\bf ik}\phi^{\bf k})(v_{\bf j}) = \xi_{\bf ij}.\eqno(2.11)$$
Abstractly, $\xi$ and $g$ are just different aspects of the same object; note that no symmetry property of $g$ has been assumed.
The inverse correlation $\xi^{-1}:V_\bullet \to V$ defines the bilinear form $G:V_\bullet \times V_\bullet \to {\bf K}$
$$G(\phi,\psi) := (\phi)\bigl(\xi^{-1}(\psi)\bigr).\eqno(2.12)$$
and thus corresponds to an element of $V \otimes V$. This definition is made so that $G$ is automatically {\bf K}-linear in its first argument in accord with the fact that $V_\bullet$ is left linear. Abstractly, $\xi^{-1}$ is represented by $\psi_b \mapsto G^{ab}\psi_b = (\xi^{-1})^{ab}\psi_b$, so the fact that $\xi^{-1}$ is the inverse of $\xi$ can be expressed as
$$G^{ac}g_{bc} = \delta^a_b = g_{cb}G^{ca};\hbox{  or  }\xi^{-1} \circ \xi = 1\ \Leftrightarrow\ (\xi^{-1})^{\bf ik}\xi_{\bf jk} = \delta^{\bf i}_{\bf j} \hbox{  and  } \xi \circ \xi^{-1} = 1\ \Leftrightarrow\ \xi_{\bf ki}(\xi^{-1})^{\bf kj} = \delta^{\bf j}_{\bf i},\eqno(2.13)$$
where $\xi^{-1}$ is left-to-right {\bf K}-linear. I write the matrix representation as
$$\xi^{-1}(\phi^{\bf j}) = v_{\bf i}(\xi^{-1})^{\bf ij} \hskip 1.25in \sum_{\bf j}\,\mu_{\bf j}\phi^{\bf j} \mapsto v_{\bf i}(\xi^{-1})^{\bf ij}\mu_{\bf j}.\eqno(2.14)$$
The matrix representation of $G$ is
$$G^{\bf ij} := G(\phi^{\bf i},\phi^{\bf j}) = \phi^{\bf i}\bigl(\xi^{-1}(\phi^{\bf j})\bigr) = \phi^{\bf i}(v_{\bf k}(\xi^{-1})^{\bf kj}) = (\xi^{-1})^{\bf ij}.\eqno(2.15)$$
Note that with these conventions,
$$G(\phi,\psi) = (\phi)\bigl(\xi^{-1}(\psi)\bigr) = \bigl((\xi \circ \xi^{-1})(\phi)\bigr)\bigl(\xi^{-1}(\psi)\bigr) = g\bigl(\xi^{-1}(\phi),\xi^{-1}(\psi)\bigr).\eqno(2.16
)$$
Hence, using abstract indices, one writes $G^{ab}$ as $g^{ab}$.

Now let $W$ be another (finite-dimensional) {\bf K}-linear space; then the linear dual $T_\bullet \in \Hom_{\bf K}(W_\bullet,V_\bullet)$ of $T \in \Hom_{\bf K}(V,W)$ is defined as for endomorphisms. When $\eta:W \to W_\bullet$ is a {\bf K}-linear, nondegenerate correlation, the linear dual $T_\bullet$ can be recast as the adjoint $^*\!T \in \Hom_{\bf K}(W,V)$ of $T$:
$${^*\! T} := \xi^{-1} \circ T_\bullet \circ \eta,\eqno(2.17)$$
and $^*\! T$ is the unique element of $\Hom(W,V)$ satisfying
$$g\bigl({^*\! T}(w),v\bigr) = h\bigl(w,T(v)\bigr),\eqno(2.18)$$
for all $v \in V$ and $w \in W$, where $h$ is the bilinear form associated to $\eta$. If $\{w_1,\ldots,w_m\}$ is a basis for $W$ with dual basis $\{\psi^1,\ldots,\psi^m\}$, the matrix representation ${^*\! {\underline T}}$ of $^*\! T$ is given by:
$$(^*\! T)(w_{\bf q}) = (\xi^{-1} \circ T_\bullet \circ \eta) (w_{\bf q}) = (\xi^{-1} \circ T_\bullet)(\eta_{\bf qp}\psi^{\bf p}) = \xi^{-1}(\eta_{\bf qp}T^{\bf p}{}_{\bf k}\phi^{\bf k}) = v_{\bf i}(\xi^{-1})^{\bf ik}T^{\bf p}{}_{\bf k}\eta_{\bf qp}.\eqno(2.19)$$
If $S := (\xi_{\bf ij})$ and $R := (\eta_{\bf qp})$, then by (2.13),
$${^* \!{\underline T}} = {^\tau\! S}^{-1}.{^\tau \!{\underline T}}.{^\tau \! R},\eqno(2.20)$$
from which (1.1--2) can be obtained. Abstractly,
$$^*\! T\ = \xi^{-1} \circ T_\bullet \circ \eta\ \leftrightarrow\ g^{ad}T^c{}_dh_{bc} = h_{bc}T^c{}_dg^{ad} =: (^*\! T)^a{}_b.\eqno(2.21{\rm a})$$
When $h$ is symmetric/skew, (2.21a) gives
$$(^*\! T)^a{}_b = \cases{T_b{}^a := h_{cb}T^c{}_dg^{ad},&\cr -T_b{}^a := -h_{cb}T^c{}_dg^{ad},&\cr}\eqno(2.21{\rm b})$$
i.e., the adjoint is $\pm T$ with its indices raised and lowered. The two cases that provide the background for this paper are $(V,\xi) = {\bf R}^{2,2}$ and $(V,\xi) = {\bf R}^2_{\rm sp}$, the symplectic plane, the latter being the two-component spinor space for ${\bf R}^{2,2}$.

The action of $\End(V)$ on $V_\bullet$ that preserves the natural pairing between $V$ and $V_\bullet$ is expressed abstractly by
$$T^a{}_b \in \End(V) \mapsto {\cal T}(T) := T_\bullet^{-1},\hbox{ acting as } \phi_b \mapsto \phi_b(T^{-1})^b{}_a.\eqno(2.22)$$
When restricting the action to the group of correlated automorphisms of a {\bf K}-bilinear symmetric or skew scalar product, the correlation intertwines between the action on $V$ and the induced action on $V_\bullet$; in other words, the two actions of orthogonal automorphisms commute with raising and lowering of indices, e.g., for a correlated automorphism $T$, of a symmetric/skew scalar product, acting on $P^a_b \in V \otimes V_\bullet$:
$$g_{fb}[T^a{}_cT^f{}_dP^{cd}] = \pm T^a{}_cP^c_dT_b{}^d = T^a{}_cP^c_d{^*\! T}^d{}_b = T^a{}_cP^c_d(T^{-1})^d{}_b,\eqno(2.23)$$
by (2.21b).

For any automorphism $T$ of $V$ and frame $\{v_1,\ldots,v_n\}$, $w_j := T(v_j) = v_iT^i{}_j$, i.e., the image of the frame $\{v_1,\ldots,v_n\}$ under $T$ is the image under the natural right action of the matrix representation of $T$ with respect to $\{v_1,\ldots,v_n\}$. Let $\underline{\bf O(p,q)}$ denote the group of standard matrix representations of {\bf O(p,q)}, i.e., the matrix representations with respect to any $\Psi$-ON basis. Let {\bf G} be a subgroup of {\bf O(p,q)} and let $G$ be its subgroup in $\underline{\bf O(2,2)}$ of standard matrix representations. $G$ defines an equivalence relation on the collection ${\cal F}$ of $\Psi$-ON frames: two $\Psi$-ON frames are $G$-{\sl related} iff one is the image of the other under the natural right action of $\underline{\bf O(p,q)}$ restricted to $G$. The equivalence classes of this $G$-relation on $\cal F$ correspond to the left cosets of $G$ in $\underline{\bf O(p,q)}$. Each equivalence class of the $G$-relation is the image of the standard basis by a left coset of {\bf G} in {\bf O(p,q)}. If {\bf G} is a normal subgroup of ${\bf O(p,q)}$, and its cosets are disconnected from each other in the natural topology, the $G$-relation is called {\bf G}-{\sl orientation} and the equivalence classes {\bf G}-{\sl orientation classes}. Normality of {\bf G} ensures that: {\bf G} has matrix representations in $G$ with respect to any $\Psi$-ON basis; {\bf G}-orientation classes are the images of {\sl any} fixed $\Psi$-ON basis by left cosets of {\bf G}, not just of the standard basis; two $\Psi$-ON bases have the same {\bf G}-orientation iff one is the image of the other by an element of {\bf G}. The topological condition ensures that, by regarding $\cal F$ as homeomorphic with {\bf O(p,q)}, the {\bf G}-orientation classes are disconnected from each other, whence {\bf G}-orientation cannot change along a continuous path in $\cal F$. One has the following notions of orientation: {\bf SO(p,q)} defines the usual notion of orientation; ${\bf O^\bfplus(p,q)} := {\bf O^\bfplus_\bfplus(p,q)} \amalg {\bf O^\bfplus_\bfminus(p,q)}$ provides a notion of {\sl time}-orientation; ${\bf O_\bfplus(p,q)} := {\bf O^\bfplus_\bfplus(p,q)} \amalg {\bf O_\bfplus^\bfminus(p,q)}$ provides a notion of {\sl space}-orientation; and ${\bf SO^\bfplus(p,q)}$ provides a notion, with four orientation classes, that I shall refer to as ${\bf SO^\bfplus}$-orientation 
\vskip 24pt
\noindent {\section 3. Spinor Theory}
\vskip 12pt
I turn now to an account of spinors for ${\bf R}^{2,2}$. Beginning with complex Clifford algebra ${\bf C}_{n}$ for the complexification of ${\bf R}^{n,0}$, ${\bf C}_{2k} \cong {\bf C}(2^k)$ as algebras, with the even part of the algebra ${\bf C}^0_{2k} \cong {\bf C}_{2k-1} \cong \left({{\bf C}(2^{k-1}) \atop 0}{0 \atop {\bf C}(2^{k-1})}\right)$, which entails that {\bf Spin(2k;C)} acts reducibly on ${\bf C}^{2k-1} \oplus {\bf C}^{2k-1}$. Writing these summands as $\cal S$ and ${\cal S}'$, then ${\bf C}_{2k} \cong \End_{\bf C}({\cal S} \oplus {\cal S}')$. Moreover ${\bf C}^{2k}$ itself has a copy lying in the odd part of ${\bf C}_{2k}$: $\left({0 \atop \Hom({\cal S},{\cal S}')}{\Hom({\cal S}',{\cal S}) \atop 0}\right)$. Now, the equation $\dim_{\bf C}\bigl(\Hom({\cal S}',{\cal S})\bigr) = 2^{2k-2} = 2k = \dim_{\bf C}({\bf C}^{2k})$ has the unique integral solution $k=2$, i.e., only for $k=2$ does one obtain ${\bf C}^{2k} \cong \Hom({\cal S}',{\cal S})$ thereby permitting an identification of $V = {\bf C}^4$ with ${\cal S} \otimes {\cal S}'_\bullet$. This observation is the basis of the particular utility of spinors in dimension four.

For the following discussion of real Clifford algebras, see Porteous (1981, 1995), but note that the scalar product space I denote by ${\bf R}^{p,q}$, Porteous denotes ${\bf R}^{q,p}$. Let ${\bf R}_{p,q}$ denote the Clifford algebra of ${\bf R}^{p,q}$; the copy of ${\bf R}^{p,q}$ within ${\bf R}_{p,q}$ will be denoted by $X$. For $a \in {\bf R}_{p,q}$: $a \mapsto \hat a$ denotes the {\sl involution} of ${\bf R}_{p,q}$, called the {\sl main involution}, uniquely defined by the fact that its restriction to $X$ is $-1_X$; $a \mapsto a^-$ denotes the {\sl anti-involution} of ${\bf R}_{p,q}$, called {\sl (Clifford) conjugation}, uniquely defined by the fact that its restriction to $X$ is also $-1_X$. An element $ a \in {\bf R}_{p,q}$ for which $\hat a = \pm a$ is called {\sl even/odd}; the even elements form a subalgebra, denoted ${\bf R}^0_{p,q}$, and the odd elements a linear subspace. The {\sl Clifford group} $\Gamma_{p,q}$ consists of those $g \in {\bf R}_{p,q}$ which are invertible and such that $gx\hat g^{-1} \in X$ for all $x \in X$; in which case $x \mapsto gx\hat g^{-1}$ is an orthogonal automorphism of $X$. Let ${\bf R}^+$ be the multiplicative group of positive real numbers; it has a canonical copy in ${\bf R}_{p,q}$, and one defines
$${\bf Pin(p,q)} := {\Gamma_{p,q} \over {\bf R}^+} \hskip 1.25in {\bf Spin(p,q)} := {\Gamma^0_{p,q} \over {\bf R}^+},\eqno(3.1)$$
where $\Gamma^0_{p,q}$ is the (normal) subgroup of even elements of $\Gamma_{p,q}$.
One can prove, however,
$${\bf Pin(p,q)} \cong \{\,g \in \Gamma_{p,q}:gg^- = \pm 1\,\} \hskip 1.25in {\bf Spin(p,q)} \cong \{\,g \in \Gamma^0_{p,q}:gg^- = \pm1\,\},\eqno(3.2)$$
whence ${\bf Pin(p,q)}$ and ${\bf Spin(p,q)}$ are typically identified with these subgroups of $\Gamma_{p,q}$ and this identification will be understood hereafter. The {\sl vector representation} is the surjective homomorphism 
$$\Upsilon:{\bf Pin(p,q)} \to {\bf O(p,q)} \hskip 1.25in \Upsilon(g) = gx\hat g^{-1},\eqno(3.3)$$
with kernel ${\bf Z}_2$. For any $g \in {\bf Pin(p,q)}$, as $\Upsilon(g) \in {\bf O(p,q)}$, by the Cartan-Dieudonn\'e Theorem there exist unit vectors $u_1,\ldots,u_k \in {\bf R}^{p,q}$ such that $\Upsilon(g)$ is the product of the hyperplane reflexions $K_{u_i}$ defined by the $u_i$, $i=1,\ldots,k$. Let $u_i$ also denote the copy of $u_i$ in $X$; then $\Upsilon(u_i) = K_{u_i}$, whence $\Upsilon(g) = \Upsilon(u_1\ldots u_k)$, i.e., $g = (\pm u_1)\ldots u_k$. It follows that every element of {\bf Spin(p,q)} can be expressed as a product of an even number of unit vectors in $X$ and that every element of ${\bf Pin(p,q)} \setminus {\bf Spin(p,q)} =: {\bf ASpin(p,q)}$ can be expressed as a product of an odd number of unit vectors in $X$ and is consequently odd.

Of particular value in low dimensions is the result
$${\rm Spin(p,q)} = \{\,g \in {\bf R}^0_{p,q}:gg^- = \pm 1\,\},\eqno(3.4)$$
valid when $\dim(X) \leq 5$ (Porteous 1981, (13.58); Porteous 1995, (16.15)). This result is easily extended to ${\bf Pin}(X)$ as follows.
\vskip 24pt
\noindent {\bf 3.5 Lemma}\hfil\break
When $\dim(X) \leq 5$,
$${\bf Pin(p,q)} = \{\,g \in {\bf R}_{p,q}:gg^- = \pm 1\,\}.$$

Proof. By (3.2), ${\bf Pin(p,q)} \leq \{\,g \in {\bf R}_{p,q}:gg^- = \pm 1\,\}$. Let $g \in {\bf R}_{p,q}$ satisfy $gg^- = \pm 1$. If $g$ is even, then it belongs to ${\bf Spin(p,q)}$ by (3.4) and thus to ${\bf Pin(p,q)}$. So, assume $g$ is odd. Choose a unit element $e \in X$, $e^2 = \mp 1$ (signs paired with $gg^- = \pm 1$), so $e^{-1} = \mp e$. Then $h = eg$ is even and $hh^- = 1$, whence $h \in {\bf Spin(p,q)}$. Therefore, $hx\hat h^{-1} \in X$, for all $x \in X$. Now $y := hx\hat h^{-1} = e(gx\hat g^{-1})\hat e^{-1}$, so $gx\hat g^{-1} = e^{-1}y\hat e = ey\hat e^{-1}$. But $y \mapsto ey\hat e^{-1}$ is reflexion in the hyperplane orthogonal to $e$, whence $ey\hat e^{-1} \in X$, for all $y \in X$. Hence, $gx\hat g^{-1} \in X$, for all $x$, and thus $g \in {\bf Pin(p,q)}$.\bull
\vskip 24pt
Another useful result is that ${\bf R}^0_{p+1,q} \cong {\bf R}_{p,q}$. Each ${\bf R}_{p,q}$ can be realized as the real algebra of endomorphisms of a right $A$-linear space $V$, $A$ the real algebra {\bf R}, {\bf C}, {\bf H}, $^2{\bf R}$ or $^2{\bf H}$, and Clifford conjugation is realized as the adjoint with respect to a certain sesquilinear scalar product on $V$. In particular, ${\bf R}_{2,2} \cong {\bf R}(4)$ and the relevant sesquilinear form is a real, skew form.

I now present a model of ${\bf R}_{2,2}$. Putting
$$K := \pmatrix{0&1\cr -1&0\cr} \hskip .75in M := \pmatrix{0&1\cr 1&0\cr} \hskip .75in N := \pmatrix{-1&0\cr 0&1\cr},\eqno(3.6)$$
$K$, $M$ and $N$ anti-commute with each other and $KMN = -1_2$. Hence, with
$$E_1 := \pmatrix{0_2&K\cr K&0_2\cr} \hskip .5in E_2 := \pmatrix{0_2&1_2\cr -1_2&0_2\cr} \hskip .5in E_3 := \pmatrix{0_2&M\cr M&0_2\cr} \hskip .5in E_4 := \pmatrix{0_2&N \cr N&0_2\cr},\eqno(3.7)$$
the $E_i$ anti-commute with each other, and satisfy $-1_4 = (E_1)^2 = (E_2)^2 = -(E_3)^2 = -(E_4)^2$ and
$$\Lambda := E_1E_2E_3E_4 = \pmatrix{1_2&0_2\cr 0_2&-1_2\cr}.\eqno(3.8)$$
Thus, $\{E_1,E_2,E_3,E_4\}$ generate ${\bf R}(4)$ as an algebra and serve as a $\Psi$-ON basis for a copy $X$ of ${\bf R}^{2,2}$ within ${\bf R}(4) \cong {\bf R}_{2,2}$. Explicitly,
$$uE_1 + vE_2 + xE_3 + yE_4 = \pmatrix{0_2&Z\cr-{^*\! Z}&0_2\cr};\hskip .75in Z := \pmatrix{v-y&x+u\cr x-u&v+y\cr},\eqno(3.9)$$
where for any matrix in ${\bf R}(2)$
$$^{\raise 6pt\hbox{$*$}\!}\pmatrix{a&c\cr b&d\cr} = \pmatrix{d&-c\cr -b&a\cr},\eqno(3.10)$$
which is the adjoint of an endomorphism of the symplectic plane ${\bf R}^2_{\rm sp}$ for matrix representations with respect to symplectic bases. By inspection, one confirms that, with $\alpha$, $\beta$, $\gamma$ and $\delta \in {\bf R}(2)$, the main involution and Clifford conjugation for this ${\bf R}(4)$ model of ${\bf R}_{2,2}$ are, respectively,
$$\pmatrix{\alpha&\gamma\cr \beta&\delta\cr}{\raise 6pt\hbox{$\widehat{}$}} = \pmatrix{\alpha&-\gamma\cr -\beta&\delta\cr} \hskip 1.25in \pmatrix{\alpha&\gamma\cr \beta&\delta\cr}^{\hbox{--}} = \pmatrix{{^*\! \alpha}&{^*\! \beta}\cr {^*\! \gamma}&{^*\! \delta}\cr},\eqno(3.11)$$
from which the even and odd elements are easily recognized. The element $\Lambda$ commutes precisely with the even elements and anti-commutes with precisely the odd elements. From (3.4) and (3.11), one easily deduces that
$$\eqalignno{{\bf Spin^\bfplus(2,2)} := \{\;g \in {\bf R}^0_{2,2}:gg^- = 1\,\} &\cong \left\{\,\pmatrix{\alpha&0_2\cr 0_2&\delta\cr}:\alpha,\ \delta \in {\bf Sp(2;R)} = {\bf SL(2;R)}\,\right\}\cr
\noalign{\vskip 6pt}
&\cong {\bf SL(2;R)} \times {\bf SL(2;R)}.&(3.12)\cr}$$

Let ${\bf P} \cong {\bf R}^4$ and regard ${\bf R}_{2,2} \cong {\bf R}(4)$ as $\End_{\bf R}({\bf P})$; {\bf P} is called the {\sl pinor space}. The element $\Lambda$ decomposes {\bf P} into its two eigenspaces, called {\sl spinor spaces}, each of which is isomorphic to ${\bf R}^2$. Write ${\bf P} = S \oplus S'$ or, with abstract indices, ${\bf P}^\alpha = S^A \oplus S^{A'}$, in which case write an element of the Clifford algebra as $\gamma^\alpha{}_\beta \in \End({\bf P})$. In light of (3.12), one regards $S^A$ and $S^{A'}$ as copies of ${\bf R}^2_{\rm sp}$, not just ${\bf R}^2$. The symplectic forms on $S^A$ and $S^{A'}$ are denoted $\epsilon_{AB}$ and $\epsilon_{A'B'}$ respectively.

More generally, let $G_{\bf i}$, ${\bf i}=1,\ldots,4$, be elements of $\Hom(S^{A'},S^A)$, so that with
$$\gamma_{\bf i} := \pmatrix{0_2&G_{\bf i}\cr -{^*\! G}_{\bf i}&0_2\cr}, \qquad\hbox{then}\qquad \gamma_{\bf i}\gamma_{\bf j} = \pmatrix{-G_{\bf i}{^*\! G}_{\bf j}&0_2\cr 0_2&-{^*\! G}_{\bf i}G_{\bf j}\cr},\eqno(3.13)$$
for ${\bf i} \not= {\bf j}$ and where $*$ is the adjoint for endomorphisms of ${\bf R}^2_{\rm sp}$. Assume
$$G_{\bf i}{^*\! G}_{\bf j} \in \End_-(S^A) \hskip 1.25in {^*\! G}_{\bf i}G_{\bf j} \in \End_-(S^{A'}),\eqno(3.14)$$
i.e., that each is skew adjoint, so that $\gamma_{\bf i}$ anti-commutes with $\gamma_{\bf j}$ (${\bf i} \not= {\bf j}$). Further, suppose
$$G_{\bf i}{^*\! G}_{\bf i} = \cases{1 \in \End(S^A),&${\bf i}=1$, 2;\cr -1 \in \End(S^A),&${\bf i}=3$, 4;\cr} \hskip 1.25in {^*\! G}_{\bf i}G_{\bf i} = \cases{1 \in \End(S^{A'}),&${\bf i}=1$, 2;\cr -1 \in \End(S^{A'}),&${\bf i}=3$, 4\cr}\eqno(3.15)$$
and
$$G_1{^*\! G}_2G_3{^*\! G}_4 = 1 \in \End(S^A) \hskip 1.25in {^*\! G}_1G_2{^*\! G}_3G_4 = -1 \in \End(S^{A'}).\eqno(3.16)$$
Observe that, by (2.21b), $-{^*\! G}^{A'}{}_B = G_B{}^{A'}$, so one can write
$$\gamma_{\bf i} = \pmatrix{0_2&(G_{\bf i })^A{}_{B'}\cr (G_{\bf i})_B{}^{A'}&0_2\cr},\eqno(3.17)$$
and the conditions (3.14--15) can be succinctly summarized as
$$(G_{({\bf i}})^A{}_{\vert B'\vert}(G_{{\bf j})})_C{}^{B'} = -g_{\bf ij}1_S\qquad\hbox{and}\qquad (G_{({\bf i}})_{\vert B\vert}{}^{A'}(G_{{\bf j})})^B{}_{C'} = -g_{\bf ij}1_{S'},\eqno(3.18)$$
where $g_{\bf ij}$ are the components of the scalar product on $V^a = {\bf R}^{2,2}$ with respect to a $\Psi$-ON basis. Equivalently, one can write
$$\gamma_{({\bf i}}\gamma_{{\bf j})} = -g_{\bf ij}1_4.\eqno(3.19)$$
Any such collection of endomorphisms $G_1,\ldots,G_4$ generates ${\bf R}_{2,2}$ as an algebra so that 
$$\gamma_1\gamma_2\gamma_3\gamma_4 = \pmatrix{1_2&0_2\cr 0_2&-1_2\cr} = \Lambda.\eqno(3.20)$$
Different choices of the $G_{\bf i}$ provide different copies of ${\bf R}^{2,2}$ in ${\bf R}(4)$, i.e., different models of ${\bf R}_{2,2}$ as ${\bf R}(4)$ but for which the main involution and Clifford conjugation take the same form as in (3.11) and the decomposition ${\bf P}^\alpha = S^A \oplus S^{A'}$ induced by $\Lambda$ is respected. Any model of ${\bf R}_{2,2}$ as ${\bf R}(4)$ which satisfies (3.11) represents $X$ as the elements of ${\bf R}(4)$ of the form
$$\pmatrix{0_2&Z\cr -{^*\! Z}&0_2\cr}.\eqno(3.21)$$
For such models, the vector representation of {\bf Spin(2,2)} is: for $g = \left({\alpha \atop 0_2}{0_2 \atop \delta}\right) \in {\bf R}^0_{2,2}$,
$$g\pmatrix{0_2&Z\cr -{^*\! Z}&0_2\cr}{\hat g}^{-1} = \pmatrix{0_2&\alpha Z\delta^{-1}\cr -{^*\!(}\alpha Z\delta^{-1})&0_2\cr}.\eqno(3.22)$$
Noting $Z \in \Hom(S^{A'},S^A)$, $Z^B{}_{B'} \mapsto \alpha^A{}_BZ^B{}_{B'}(\delta^{-1})^{B'}{}_{A'}$ suffices to describe the vector representation of ${\bf Spin(2,2)}$. Using the symplectic form on $S^{A'}$, one can rewrite this expression as:
$$\eqalignno{Z^{BB'} &\mapsto \alpha^A{}_BZ^B{}_{B'}(\delta^{-1})^{B'A'}\cr
& = -\alpha^A{}_BZ^{BB'}(\delta^{-1})_{B'}{}^{A'}\cr 
&= \alpha^A{}_BZ^{BB'}({^*\! \delta^{-1}})^{A'}{}_{B'}\qquad\hbox{using (2.21b)}\cr 
&= \alpha^A{}_BZ^{BB'}(\pm\delta^{A'}{}_{B'})\qquad\hbox{according as $gg^- = \pm 1$},\cr
&= \pm\alpha^A{}_B\delta^{A'}{}_{B'}Z^{BB'}\qquad\hbox{sign according as $gg^- = \pm 1$}.&(3.23)\cr}$$
Analogous to (3.12), write:
$$\eqalignno{{\bf Spin^\bfminus(2,2)} &:= \{\,g \in {\bf R}^0_{2,2}:gg^- = -1\,\} \cong \left\{\,\pmatrix{\alpha&0_2\cr 0_2&\delta\cr}:\alpha,\ \delta \in {\bf ASp(2;r)} = {\bf ASL(2;R)}\,\right\}\cr
\noalign{\vskip 6pt}
&\cong {\bf ASL(2;R)} \times {\bf ASL(2;R)},&(3.24)\cr}$$
where {\bf ASp(2;R)} is the set of anti-symplectic automorphisms of ${\bf R}^2_{\rm sp}$ and {\bf ASL(2;R)} is the subset of {\bf GL(2;R)} whose elements have determinant $-1$. For the choice (3.6--7), from (3.9),
$$(Z^{\bf BB'}) = \pmatrix{u+x&y-v\cr y+v&u-x\cr}.\eqno(3.25)$$
For $g \in {\bf ASpin(2,2)}$, $g$ is odd, hence of the form $\left({0_2 \atop \beta}{\gamma \atop 0_2}\right)$, with $\gamma{^*\! \gamma} = \pm 1 = \beta{^*\! \beta}$. Denote the subset for which $gg^- = 1$ by ${\bf ASpin^\bfplus(2,2)}$, and the subset for which $gg^- = -1$ by ${\bf ASpin^\bfminus(2,2)}$. The vector representation is
$$g\pmatrix{0_2&Z\cr -{^*\! Z}&0_2\cr}{\hat g}^{-1} = \pmatrix{0_2&\gamma\,{^*\! Z}\beta^{-1}\cr -{^*\! (}\gamma\,{^*\! Z}\beta^{-1})&0_2\cr},\eqno(3.26)$$
and thus may be represented by
$$\eqalignno{Z^{BB'} &\mapsto \gamma^A{}_{B'}({^*\! Z})^{B'}{}_B(\beta^{-1})^{BA'}\cr
&= \gamma^A{}_{B'}(-Z_B{}^{B'})(\beta^{-1})^{BA'}\qquad\hbox{using (2.21b)}\cr
&= \gamma^A{}_{B'}(\beta^{-1})_B{}^{A'}Z^{BB'}\cr
&= \gamma^A{}_{B'}(\pm{^*\! \beta})_B{}^{A'}Z^{BB'}\qquad\hbox{according as $gg^- = \pm 1$}\cr
&= \pm\gamma^A{}_{B'}(-\beta)^{A'}{}_BZ^{BB'}\qquad\hbox{using (2.21b)}\cr
&= \mp\gamma^A{}_{B'}\beta^{A'}{}_BZ^{BB'}\qquad\hbox{sign according as $gg^- = \pm 1$}.&(3.27)\cr}$$
Of particular interest is (3.27) when $g$ is a unit vector $u \in X$, whence by (3.21) $g$ takes the form $\left({0_2 \atop -{^*\! Z_u}} {Z_u \atop 0_2}\right)$, with ${^*\! Z}_u = \pm Z_u^{-1}$ according as $gg^- = \pm1$, i.e., according as $g^2 = \mp1$ (since $g^- = -g$), i.e., according as $u$ is time/space like. So, (3.27) becomes
$$Z^{BB'} \mapsto \mp(Z_u)^A{}_{B'}[-({^*\! Z_u})^{A'}{}_B]Z^{BB'} = \pm(Z_u)^A{}_{B'}[\pm(Z_u^{-1})^{A'}{}_B]Z^{BB'} = (Z_u)^A{}_{B'}(Z_u^{-1})^{A'}{}_BZ^{BB'},\eqno(3.28)$$
independently of whether $gg^- = \pm 1$, i.e., of whether $u$ is time/space like. For example, with $u = E_1$, $E_2$, $E_3$ and $E_4$, respectively, one obtains, as expected, the hyperplane reflexions
$$\vcenter{\openup2\jot \halign{$\hfil#$&&${}#\hfil$&\qquad$\hfil#$\cr
\pmatrix{u+x&y-v\cr y+v&u-x\cr} &\mapsto \pmatrix{-u+x&y-v\cr y+v&-(u+x)\cr}& \pmatrix{u+x&y-v\cr y+v&u-x\cr} &\mapsto \pmatrix{u+x&y+v\cr y-v&u-x\cr}\cr
\pmatrix{u+x&y-v\cr y+v&u-x\cr} &\mapsto \pmatrix{u-x&y-v\cr y+v&u+x\cr} & \pmatrix{u+x&y-v\cr y+v&u-x\cr} &\mapsto \pmatrix{u+x&-(y+v)\cr v-y&u-x\cr}\cr}}\eqno(3.29)$$
respectively.

One easily computes that an even product of unit vectors $u_i\leftrightarrow\ \left({0_2 \atop -{^*\! Z_i}}{Z_i \atop 0_2}\right)$ in $X$ takes the form,
$$u_1\ldots u_{2k}\ \leftrightarrow\ \pmatrix{\pm(Z_1{^*\! Z_2})\ldots(Z_{2k-1}{^*\! Z}_{2k})&0_2\cr 0_2&\pm({^*\! Z}_1Z_2)\ldots({^*\! Z}_{2k-1}Z_{2k})\cr},\eqno(3.30)$$
with sign according as $k$ is even/odd. Now $u_1\ldots u_{2k} \in {\bf Spin^\bfplus(2,2)}$ iff there is an even number of space-like unit vectors amongst the $u_i$ (equivalently, an even number of unit time-like vectors). Substituting into (3.23), the vector representation of $g := u_1\ldots u_{2k}$ is, with the sign according as $gg^- = \pm 1$,
$$\eqalignno{\pm\alpha^A{}_B\delta^{A'}{}_{B'}Z^{BB'} &= \pm[(Z_1{^*\! Z}_2)\ldots(Z_{2k-1}{^*\! Z}_{2k})]^A{}_B[({^*\! Z}_1Z_2)\ldots({^*\! Z}_{2k-1}Z_{2k}]^{A'}{}_{B'}Z^{BB'}\cr
&= [(Z_1Z_2^{-1})\ldots(Z_{2k-1}Z_{2k}^{-1})]^A{}_B[(Z_1^{-1}Z_2)\ldots(Z_{2k-1}^{-1}Z_{2k})]^{A'}{}_{B'}Z^{BB'}&(3.31)\cr}$$
where one uses ${^*\! Z_i} = \pm Z_i^{-1}$ once for each $i=1,\ldots,2k$, the sign according as $u_i$ is time/space like, so the cumulative negative signs from the space-like vectors exactly equals the sign at the front of the expression in the first line and so cancels it. Since every element $g$ of {\bf Spin(2,2)} is a product of an even number of unit vectors, (3.31) gives a form of the vector representation of {\bf Spin(2,2)} that is independent of whether $gg^- = \pm 1$. Of course that form can be obtained directly as $2k$ applications of (3.28). For an element $g = u_1\ldots u_{2k+1} \in {\bf ASpin(2,2)}$, one need only compose (3.31) with (3.28) (effectively $2k+1$ applications of (3.27)) to obtain, with $g\ \leftrightarrow\ \left({0_2 \atop \beta}{\gamma \atop 0_2}\right)$, the vector representation in a form independent of whether $gg^- = \pm1$:
$$\mp\gamma^A{}_{B'}\beta^{A'}{}_BZ^{BB'} = [(Z_1Z_2^{-1})\ldots(Z_{2k-1}Z_{2k}^{-1})Z_{2k+1}]^A{}_{B'}[(Z_1^{-1}Z_2)\ldots(Z_{2k-1}^{-1}Z_{2k})Z_{2k+1}^{-1}]^{A'}{}_BZ^{BB'}.\eqno(3.32)$$
To obtain an analogue of the spinor formalism of PRI, one must choose an explicit identification of $V^a = {\bf R}^{2,2}$ with $S^A \otimes S^{A'}$. In PRI, such a choice is avoided in so far as they regard space-time vectors as certain kinds of space-time spinors, i.e., on physical grounds, spinors are regarded as primitive and space-time vectors as derived objects, not an independent kind of object, the space of which happens to be isomorphic to a tensor product of spinor spaces by virtue of the geometry. Such considerations are irrelevant to the mathematical proceedings here. $V^a$ and $S^A \otimes S^{A'}$ are independent spaces, which happen to be isomorphic. PRI's interpretation permits them, in the abstract index notation, to regard abstract indices $a$ and $AA'$ as just different forms of the same index, so that they may write equations such as $\ell^a = o^Ao^{A'}$. Such an equation will not be technically correct here; but I will adopt a convention that permits one to write it, with an implicit meaning understood.

Choose {\sl spin frames} $\{o^A,\iota^{A'}\}$ for $S^A$ and $\{o^{A'},\iota^{A'}\}$ for $S^{A'}$, i.e., bases satisfying $\iota^A o_A = 1$ and $\iota^{A'}o_{A'} = 1$, and define the following elements of $\Hom(S^{A'},S^A)$:
$$\vcenter{\openup1\jot \halign{$\hfil#$&&${}#\hfil$&\qquad$\hfil#$\cr
(G_1)^A{}_{A'} &:= o^Ao_{A'} + \iota^A\iota_{A'}  &&:o^{A'} \mapsto -\iota^A\hbox{  \&  }:\iota^{A'} \mapsto o^A, & \hbox{a symplectomorphism}\cr
(G_2)^A{}_{A'} &:= \iota^Ao_{A'} - o^A\iota_{A'}  &&:o^{A'} \mapsto o^A\hbox{  \&  }:\iota^{A'} \mapsto \iota^A, & \hbox{a symplectomorphism}\cr
(G_3)^A{}_{A'} &:= o^Ao_{A'} - \iota^A\iota_{A'}  &&:o^{A'} \mapsto \iota^A\hbox{  \&  }:\iota^{A'} \mapsto o^A, & \hbox{an anti-symplectomorphism}\cr
(G_4)^A{}_{A'} &:= \iota^Ao_{A'} + o^A\iota_{A'}  &&:o^{A'} \mapsto -o^A\hbox{  \&  }:\iota^{A'} \mapsto \iota^A, & \hbox{an anti-symplectomorphism}\cr}}\eqno(3.33)$$
Now $G_1,\ldots,G_4$ satisfy (3.14--16); in particular, with respect to the chosen spin frames, they have matrix representations $K$, $1_2$, $M$ and $N$ respectively. Now define
$$\sigma_{\bf i}{}^A{}_{A'} := {1 \over \sqrt2}(G_{\bf i})^A{}_{A'},\qquad{\bf i}=1,\ldots,4.\eqno(3.34)$$
Let $\{\phi^1,\ldots,\phi^4\}$ denote the dual basis of the standard basis $\{e_1,\ldots,e_4\}$ of ${\bf R}^4$ and define
$$\sigma_a{}^A{}_{A'} := \phi^{\bf i}_a\sigma_{\bf i}{}^A{}_{A'} \in \Hom(S^{A'},S^A \otimes V_a),\eqno(3.35)$$
where, as usual, the Einstein summation convention is assumed for a pair of {\sl like} indices, one contravariant, one covariant (here the numerical index {\bf i}). Indices are raised/lowered on $\sigma_a{}^A{}_{A'}$, of course, using the relevant metric ($g_{ab}$, $\epsilon_{AB}$, or $\epsilon_{A'B'}$).
\vskip 24pt
\noindent {\bf 3.36 Lemma}\hfil\break
$$\sigma_{\bf i}{}^{AA'}\sigma_{\bf j}{}_{AA'} = g_{\bf ij},$$
where $g_{\bf ij} = g_{ab}e^a_{\bf i}e^b_{\bf j} =s_{2,2}(e_{\bf i},e_{\bf j})$.

Proof. From the first equation in (3.18), one deduces 
$$(G_{\bf i})^{AB'}(G_{\bf j})_{CB'} + (G_{\bf j})^{AB'}(G_{\bf i})_{CB'} = 2g_{\bf ij}1_S,$$
whence
$$2(G_{\bf i})^{AB'}(G_{\bf j})_{AB'} = 2g_{\bf ij}\tr(1_S) = 4g_{\bf ij},$$
and the assertion follows from the definition (3.34).\bull
\vskip 24pt
Now define the isomorphism
$$V^a \to S^A \otimes S^{A'} \hskip 1.25in v^a \mapsto v^a\sigma_a{}^{AA'} =: v^{AA'},\eqno(3.37)$$
observing that $v^a = v^{\bf i}e^a_{\bf i} \mapsto v^{\bf i}e^a_{\bf i}\phi^{\bf j}_a\sigma_{\bf j}{}^{AA'} = v^{\bf i}\sigma_{\bf i}{}^{AA'}$. With respect to the chosen spin frames, one computes that
$$\sigma_1{}^{\bf AA'} = {1 \over \sqrt2}\pmatrix{1&0\cr 0&1\cr} \qquad \sigma_2{}^{\bf AA'} = {1 \over \sqrt2}\pmatrix{0&-1\cr 1&0\cr} \qquad \sigma_3{}^{\bf AA'} = {1 \over \sqrt2}\pmatrix{1&0\cr 0&-1\cr} \qquad \sigma_4{}^{\bf AA'} = {1 \over \sqrt2}\pmatrix{0&1\cr 1&0\cr},\eqno(3.38)$$
whence the isomorphism takes the form, with respect to the chosen spin frames and standard basis of $V^a$:
$$(u,v,x,y) \mapsto {1 \over \sqrt2}\pmatrix{u+x&y-v\cr y+v&u-x\cr},\eqno(3.39)$$
consistent with (3.25). Note that
$$\sigma^a{}_{AA'} :=  g^{ab}\sigma_{bAA'} = g^{ab}\phi^{\bf k}_b\sigma_{{\bf k}AA'} = e^a_{\bf i}g^{\bf ik}\sigma_{{\bf k}AA'} =: e^a_{\bf i}\sigma^{\bf i}{}_{AA'},\eqno(3.40)$$
whence the induced isomorphism on $V_a$, preserving the natural pairing between $V^a$ and $V_a$, can be written
$$\psi_a \in V_a \mapsto \psi_a\sigma^a{}_{AA'};\eqno(3.41)$$
in particular,
$$\delta^{\bf j}_{\bf i} = e^a_{\bf i}\phi^{\bf j}_a \mapsto (e^a_{\bf i}\sigma_a{}^{AA'})(\phi^{\bf j}_b\sigma^b{}_{AA'}) = \sigma_{\bf i}{}^{AA'}\sigma^{\bf j}{}_{AA'} = \delta^{\bf j}_{\bf i}$$
by (3.36). The significance of the factor $\sqrt2$ incorporated into the identification of $V^a$ with $S^A \otimes S^{A'}$ is that, by (3.36),
$$\sigma_a{}^{AA'}\sigma_b{}^{BB'}\epsilon_{AB}\epsilon_{A'B'} = \phi^{\bf i}_a\phi^{\bf j}_b\sigma_{\bf i}{}^{AA'}\sigma_{{\bf j}AA'} = \phi^{\bf i}_a\phi^{\bf j}_bg_{\bf ij} = g_{ab}.\eqno(3.42)$$
Without the $\sqrt2$ factor, one would have a factor of 2 in the right-hand side. In PRI, the expression $\epsilon_{AB}\epsilon_{A'B'}$ is taken to {\sl define} the metric $g_{ab}$.

Note that for $\sigma_a{}^{AA'}$, there is no identification of the index $a$ with the index pair $AA'$, and no contraction across them implied, contra PRI (3.1.4). In order to identify indices $a$ and $AA'$, one must regard the isomorphism (3.37) as (non canonically) identifying $V^a$ with $S^A \otimes S^{A'}$; then, instead of writing $v^a\ \leftrightarrow v^{AA'}$ or $g_{ab}\ \leftrightarrow \epsilon_{AB}\epsilon_{A'B'}$, I will adopt the convention that one writes $v^a = v^{AA'}$ and $g_{ab} = \epsilon_{AB}\epsilon_{A'B'}$, with the object $\sigma_a{}^{AA'}$ that effects the identification implicitly understood. With this convention in place, one can then treat the indices $a$ and $AA'$ as (effectively) different forms of the same label as in PRI (though for a different reason). Because of this convention, when explicitly employing the sigma's, one must write them as $\sigma_a{}^{AA'}$ but the convention does not apply to their indices. In PRI, the sigma objects take a concrete form and an abstract form. First, the expressions in (3.38) are
$$\sigma_{\bf i}{}^{\bf AA'} = \sigma_{\bf i}{}^{AA'}\epsilon_A{}^{\bf A}\epsilon_{A'}{}^{\bf A'} = (e^a_{\bf i}\sigma_a{}^{AA'})\epsilon_A{}^{\bf A}\epsilon_{A'}{}^{\bf A'},\eqno(3.43)$$
which is the analogue of the first equation in PRI (3.1.37), i.e., $\sigma_{\bf i}{}^{\bf AA'}$ are the Infeld-van der Waerden symbols for the chosen spin frames and $\Psi$-ON basis. Second, in PRI, each $\phi^{\bf j}_a$ is in fact a spinorial object. Here, one can derive the equations
$$\phi^1_a = {o_Ao_{A'} + \iota_A\iota_{A'} \over \sqrt2} \qquad \phi^2_a = {\iota_Ao_{A'} - o_A\iota_{A'} \over \sqrt2} \qquad \phi^3_a = {o_Ao_{A'} - \iota_A\iota_{A'} \over \sqrt2} \qquad \phi^4_a = {\iota_Ao_{A'} + o_A\iota_{A'} \over \sqrt2},\eqno(3.44)$$
recalling that here the equality signs stand for the identification (3.37), but in PRI (3.44) are effectively definitions of the $\phi^{\bf i}_a$. Upon substituting these expressions into (3.35) one finds
$$\sigma_a{}^B{}_{B'} = \epsilon_A{}^B\epsilon_{A'B'},\eqno(3.45)$$
for which compare PRII (B.86).

Of course, any $\Psi$-ON basis of $V^a$ could be employed in place of the standard basis to define an identification of $V^a$ with $S^A \otimes S^{A'}$. I take (3.37) to be understood hereafter; the choice of spin frames underlying (3.37) is completely arbitrary, the spin frames chosen will be referred to as the {\sl distinguished spin frames} and hereafter denoted by $\{\check o^A,\check\iota^A\}$ and $\{\check o^{A'},\check\iota^{A'}\}$.

When working exclusively with the identification between $V^a$ and $S^A \otimes S^{A'}$, it is convenient to rewrite the vector representation (3.23) and (3.27) as follows: for $\alpha^A{}_B \in \End(S)$, $\delta^{A'}{}_{B'} \in  \End(S')$, $\gamma^A{}_{B'} \in \Hom(S',S)$, and $\beta^{A'}{}_B \in \Hom(S,S')$, the vector representations of {\bf Spin(2,2)} and {\bf ASpin(2,2)} are
$$v^{BB'}  \mapsto \alpha^A{}_B\delta^{A'}{}_{ B'}v^{BB'} \hskip 1.25in v^{BB'} \mapsto \beta^{A'}{}_B\gamma^A{}_{B'}v^{BB'},\eqno(3.46)$$
respectively. When $\alpha^A{}_B$ and $\delta^{A'}{}_{B'}$ are each symplectic automorphisms (equivalently, each belongs to {\bf SL(2;R)}), one has the vector representation of ${\bf Spin^\bfplus(2,2)}$ as ${\bf SO^\bfplus(2,2)}$ and the element represented is $\pm\left({\alpha \atop 0_2}{0_2 \atop \delta}\right)$; when each is an anti-symplectic automorphism (equivalently, each belongs to {\bf ASL(2;R)}), one has the vector representation of ${\bf Spin^\bfminus(2,2)}$ as ${\bf O^\bfminus_\bfminus(2,2)}$ and the element represented is $\pm\left({\alpha \atop 0_2}{0_2 \atop -\delta}\right)$. When $\beta^{A'}{}_B$ and $\gamma^A{}_{B'}$ are each symplectic isomorphisms, one has the vector representation of ${\bf ASpin^\bfplus(2,2)}$ as ${\bf O_\bfplus^\bfminus(2,2)}$ and the element represented is $\pm\left({0_2 \atop \beta}{\gamma \atop 0_2}\right)$; when each is an anti-symplectic isomorphism, one has the vector representation of ${\bf ASpin^\bfminus(2,2)}$ as ${\bf O_\bfminus^\bfplus(2,2)}$ and the element represented is $\pm\left({0_2\atop \beta}{-\gamma \atop 0_2}\right)$. The identification of the components of {\bf O(2,2)} in the last sentence will be justified below; otherwise all that has been done here is to absorb the minus sign in (3.23) and (3.27) into the element represented so as to gain more uniform expressions for the representation of {\bf O(2,2)} in terms of spinors.

Note that while a homomorphism does not in general possess a determinant, the determinant of a homomorphism $T$ between symplectic spaces of the same dimension does: one defines the determinant $\det(T)$ to be the determinant of any matrix representation of $T$ with respect to symplectic bases. It is easily seen that the determinant is indeed well defined. In the present context, it follows that for $T \in \Hom(S',S)$, say,
$$\epsilon_{AB}T^A{}_{C'}T^B{}_{D'} = \det(T)\epsilon_{C'D'},\eqno(3.47)$$
with an analogous equation if $T \in \Hom(S,S')$. In particular, $\beta^{A'}{}_B$ and $\gamma^A{}_{B'}$ are characterized as (anti-)symplectic isomorphisms iff their determinant is $(-)1$.

From (3.37), it follows that $v^a$ is null iff $v^{AA'}v^{BB'}\epsilon_{AB}\epsilon_{A'B'} = 0$; this condition is equivalent to $v^{AA'}$ being a decomposable element of $S^A \otimes S^{A'}$. Thus, given any pair of spin frames $\{o^A,\iota^A\}$ and $\{o^{A'},\iota^{A'}\}$, one can define a basis of null vectors for $V^a$:
$$\ell^a := o^Ao^{A'} \qquad \tilde m^a := \iota^Ao^{A'} \qquad n^a := \iota^A\iota^{A'} \qquad m^a := o^A \iota^{A'},\eqno(3.48)$$
which is called a {\sl null tetrad}; note that $\{\ell^a,\tilde m^a,n^a,-m^a\}$ is a Witt basis, but by null tetrad I will mean a basis of null vectors of the particular form specified in (3.48) and so ordered. Each null tetrad defines a $\Psi$-ON frame $\{U^a,V^a,X^a,Y^a\}$ according to
$$U^a := {\ell^a + n^a \over \sqrt2} \qquad V^a := {\tilde m^a - m^a \over \sqrt2} \qquad X^a := {\ell^a - n^a \over \sqrt2} \qquad Y^a := {\tilde m^a + m^a \over \sqrt2}.\eqno(3.49)$$
Given that ${\bf Spin^\bfplus(2,2)} \cong {\bf SL(2;R)} \times {\bf SL(2;R)}$ two-one covers ${\bf SO^\bfplus(2,2)}$, it will follow that the $\Psi$-ON frames constructed by means of (3.48--49) constitute an $\bf SO^\bfplus$-orientation class. It is easily confirmed that from the distinguished spin frames (3.48-49) constructs the standard basis (e.g., $\check\ell^a := \sigma^a{}_{AA'}\check o^A\check o^{A'} = e^a_{\bf j}\sigma^{\bf j}{}_{AA'}\check o^A\check o^{A'} = e^a_{\bf j}(g^{{\bf j}1} - g^{{\bf j}3})/\sqrt2 = (e^a_1+e^a_3)/\sqrt2$; similarly, $\check{\tilde m}^a := \sigma^a{}_{AA'}\check\iota^A\check o^{A'} = (e^a_2+e^a_4)/\sqrt2$, $\check n^a := \sigma^a{}_{A'}\check\iota^A\check\iota^{A'} = (e^a_1 - e^a_3)\sqrt2$, $\check m^a := \sigma^a{}_{AA'}\check o^a\check\iota^{A'} = (-e^a_2 + e^a_4)/\sqrt2$).

Now, $\Psi$-ON bases $\{U^a,V^a,X^a,Y^a\}$ correspond one-to-one with Witt bases $\{E_1,E_2,E_3,E_4\}$ according to
$$\displaylines{E^a_1 = {U^a + X^a \over \sqrt2} \qquad E^a_2 = {V^a + Y^a \over \sqrt2} \qquad E^a_3 = {U^a - X^a \over \sqrt2} \qquad E^a_4 = {V^a - Y^a \over \sqrt2}\cr
\hfill\llap(3.50)\cr
U^a = {E^a_1+E^a_3 \over \sqrt2} \qquad V^a = {E^a_2+E^a_4 \over \sqrt2} \qquad X^a = {E^a_1 - E^a_3 \over \sqrt2} \qquad Y^a = {E^a_2 - E^a_4 \over \sqrt2}\cr}$$
Write $E^a_1 = \alpha^A\pi^{A'}$. Since $E_1$ and $E_2$ are distinct, orthogonal, null vectors, either $E^a_2 = \beta^A\pi^{A'}$ with $\beta^A\alpha_A \not= 0$ or $E^a_2 = \alpha^A\eta^{A'}$ with $\eta^{A'}\pi_{A'} \not= 0$. 

Consider the first case. By scaling $\alpha^A$ and $\beta^A$ by a suitable $\lambda \in {\bf R}^+$ and $\pi^{A'}$ by $\lambda^{-1}$, one can suppose $\beta^A\alpha_A = \pm 1$. Writing $E^a_3 = \nu^A\zeta^{A'}$, $\zeta^{D'}\pi_{D'} \not= 0$ whence $E^a_2 \perp E^a_3$ entails $\nu^D\beta_D = 0$. So, one can write $E^a_3 = \beta^A\xi^{A'}$, where $\xi^{A'}\pi_{A'} \not= 0$, whence $g_{ab}E^a_1E^b_3 = 1$ implies $\xi^{A'}\pi_{A'} = \pm1$ (with signs corresponding to the earlier sign choice). Writing $E^4_a = \gamma^A\kappa^{A'}$, $g_{ab}E^a_4E^b_2 = 1$ and $g_{ab}E^a_4E^b_1 = 0$ imply $\gamma^A\alpha_A = 0$, whence $g_{ab}E^a_4E^b_3 = 0$ entails $\kappa^{D'}\xi_{D'} = 0$. So, $E^a_4 = \lambda\alpha^A\xi^{A'}$, for some $\lambda \in {\bf R}^*$. But $g_{ab}E^a_4E^b_2 = 1$ implies $\lambda = -1$, in which case
$$E^a_1 = \alpha^A\pi^{A'} \qquad E^a_2 = \beta^A\pi^{A'} \qquad E^a_3 = \beta^A\xi^{A'} \qquad E^a_4 = -\alpha^A\xi^{A'}.\eqno(3.51{\rm a})$$
When $\pi^{A'}\xi_{A'} = 1 = \beta^A\alpha_A$, $\{\alpha^A,\beta^A\}$ and $\{\pi^{A'},\xi^{A'}\}$ are spin frames and $\{E_1,E_2,E_3,-E_4\}$ is a null tetrad; when $\pi^{A'}\xi_{A'} = -1 = \beta^A\alpha_A$, $\{\beta^A,\alpha^A\}$ and $\{\xi^{A'},\pi^{A'}\}$ are spin frames and $\{E_3,-E_4,E_1,E_2\}$ is a null tetrad. Of course the transformation which maps $\alpha^A \mapsto \beta^A$ and $\beta^A \mapsto \alpha^A$ is an element of {\bf ASL(2;R)}.

In the second case, a similar argument shows that the Witt basis takes the form
$$E^a_1 = \alpha^A\pi^{A'} \qquad E^a_2 = \alpha^A\xi^{A'} \qquad E^a_3 = \beta^A\xi^{A'} \qquad E^a_4 = -\beta^A\pi^{A'},\eqno(3.51{\rm b})$$
either for spin frames $\{\alpha^A,\beta^A\}$ and $\{\pi^{A'},\xi^{A'}\}$ or spin frames $\{\beta^A,\alpha^A\}$ and $\{\xi^{A'},\pi^{A'}\}$, in which cases $\{E_1,-E_4,E_3,E_2\}$ and $\{E_3,E_2,E_1,-E_4\}$ are null tetrads respectively.

In all cases it is easy to check that the Witt basis determines the corresponding pair of spin frames uniquely up to an overall sign. 

These four possibilities have the following expected interpretation. Starting with the standard basis, the associated spin frames are the distinguished pair (and their negatives). Any other spin frame $\{o^A,\iota^A\}$ of $S$ can be obtained as the image of $\{\check o^A,\check\iota^A\}$ by a (unique) $\alpha \in {\bf Sp}(S) \cong {\bf SL(2;R)}$; similarly, any spin frame $\{o^{A'},\iota^{A'}\}$ of $S'$ is the image of $\{\check o^{A'},\check\iota^{A'}\}$ by a (unique) $\delta \in {\bf Sp}(S') \cong {\bf SL(2;R)}$ (equivalently, if $\underline\alpha$ is the matrix representation of $\alpha$ with respect to $\{\check o^A,\check\iota^A\}$, then $\{o^A,\iota^A\}$ is the image under the natural right action of {\bf SL(2;R)} by $\underline\alpha$ of $\{\check o^A,\check\iota^A\}$, viz., $\epsilon_{\bf B}{}^A = \alpha(\check\epsilon_{\bf B}{}^A) = \check\epsilon_{\bf A}{}^A\underline\alpha^{\bf A}{}_{\bf B}$, with an analogous statement for $\delta$ and the primed spin frames). By (3.12) and (3.46), the action of the pair $\alpha$ and $\delta$ on the distinguished spin frames corresponds via the vector representation of ${\bf Spin^\bfplus(2,2)}$ to the action of an element of ${\bf SO^\bfplus(2,2)}$ on the standard basis, e.g.,
$$\hat U^a := {\hat\ell^a + \hat n^a \over \sqrt2} = {\hat o^A\hat o^{A'} + \hat\iota^A\hat\iota^{A'} \over \sqrt2} = {\alpha^A{}_B\delta^{A'}{}_{B'}[\check o^B\check o^{B'} + \check\iota^B\check\iota^{B'}] \over \sqrt2} = \alpha^A{}_B\delta^{A'}{}_{B'}\check U^{BB'} = \alpha^A{}_B\delta^{A'}{}_{B'}e^{BB'}_1,$$
and one has the correspondence: for arbitrary spin frames $\{o^A,\iota^A\}$ and $\{o^{A'},\iota^{A'}\}$,
$$\eqalignno{\pm(\{o^A,\iota^A\},\{o^{A'},\iota^{A'}\})\ &\leftrightarrow\ \hbox{null tetrads    } \{\ell^a,\tilde m^a,n^a,m^a\}\cr
&\leftrightarrow\ \hbox{$\Psi$-ON frames in the ${\bf SO^\bfplus}$-orientation class of the standard basis},&(3.52)\cr}$$
as determined by (3.48--49).

Now consider the anti-symplectic automorphisms induced by:
$$\check o^A \mapsto \alpha^A := \check\iota^A \qquad \check\iota^A \mapsto \beta^A := \check o^A \qquad \check o^{A'} \mapsto \pi^{A'} := \check\iota^{A'} \qquad \check\iota^{A'} \mapsto \xi^{A'} := \check o^{A'},\eqno(3.53)$$
which induces
$$\vcenter{\openup1\jot \halign{$\hfil#$&&${}#\hfil$&\qquad$\hfil#$\cr
\check\ell^a &\mapsto E^a_1 := \alpha^A\pi^{A'} = \check n^a & \check{\tilde m}^a &\mapsto E^a_2 := \beta^A\pi^{A'} = \check m^a\cr
\check n^a &\mapsto E^a_3 := \beta^A\xi^{A'} = \check\ell^a & \check m^a &\mapsto -E^a_4 := \alpha^A\xi^{A'} = \check{\tilde m}^a,\cr}}\eqno(3.54)$$
i.e., $\{E_3,-E_4,E_1,E_2\}$ is a null tetrad, as in the second possibility in (3.51a), and
$$\vcenter{\openup1\jot \halign{$\hfil#$&&${}#\hfil$&\qquad$\hfil#$\cr
e^a_1 &\mapsto U^a := \displaystyle{E^a_1 + E^a_3 \over \sqrt2} = e^a_1 & e^a_2 &\mapsto V^a := \displaystyle{E^a_2 + E^a_4 \over \sqrt2} = -e^a_2\cr
e^a_3 &\mapsto X^a := \displaystyle{E^a_1 - E^a_3 \over \sqrt2} = -e^a_3 & e^a_4 &\mapsto Y^a := \displaystyle{E^a_2 - E^a_4 \over \sqrt2} = e^a_4,\cr}}\eqno(3.55)$$
which is, as expected, an ${\bf O^\bfminus_\bfminus(2,2)}$ transformation. Composing the anti-symplectomorphisms (3.53) each with symplectomorphisms (of $S$ and $S'$ respectively) is then equivalent to composing the ${\bf O^\bfminus_\bfminus(2,2)}$ transformation (3.55) with ${\bf SO^\bfplus(2,2)}$ transformations. Thus, one obtains the following correspondence: for arbitrary spin frames $\{o^A,\iota^A\}$ and $\{o^{A'},\iota^{A'}\}$,
$$\eqalignno{\pm(\{\iota^A,o^A\},\{\iota^{A'}o^{A'}\})\ &\leftrightarrow\ \hbox{Witt bases of the form    } \{E^a_1 = \iota^A\iota^{A'}, E^a_2 = o^A\iota^{A'}, E^a_3 = o^Ao^{A'}, E^a_4 = -\iota^Ao^{A'}\}\cr
&\leftrightarrow\ \hbox{$\Psi$-ON bases in the ${\bf SO^\bfplus}$-orientation class that is the image}\cr
&\qquad \hbox{ of the standard basis under ${\bf O^\bfminus_\bfminus(2,2)}$.}&(3.56)\cr}$$

Now consider the symplectic isomorphisms induced by
$$\check o^A \mapsto \pi^{A'} := \check o^{A'} \qquad \check\iota^A \mapsto \xi^{A'} := \check\iota^{A'} \qquad \check o^{A'} \mapsto \alpha^A := \check o^A \qquad \check\iota^{A'} \mapsto \beta^A := \check\iota^A.\eqno(3.57)$$
These isomorphisms induce
$$\vcenter{\openup1\jot \halign{$\hfil#$&&${}#\hfil$&\qquad$\hfil#$\cr
\check\ell^a &\mapsto E^a_1 := \alpha^A\pi^{A'} = \check\ell^a & \check{\tilde m}^a &\mapsto E^a_2 := \alpha^A\xi^{A'} = \check m^a\cr
\check n^a &\mapsto E^a_3 := \beta^A\xi^{A'} = \check n^a & \check m^a &\mapsto -E^a_4 := \beta^A\pi^{A'} = \check{\tilde m}^a,\cr}}\eqno(3.58)$$
whence
$$\vcenter{\openup1\jot \halign{$\hfil#$&&${}#\hfil$&\qquad$\hfil#$\cr
e^a_1 &\mapsto U^a := \displaystyle{E^a_1 + E^a_3 \over \sqrt2} = e^a_1 & e^a_2 &\mapsto V^a := \displaystyle{E^a_2 + E^a_4 \over \sqrt2} = -e^a_2\cr
e^a_3 &\mapsto X^a := \displaystyle{E^a_1 - E^a_3 \over \sqrt2} = e^a_3 & e^a_4 &\mapsto Y^a := \displaystyle{E^a_2 - E^a_4 \over \sqrt2} = e^a_4,\cr}}\eqno(3.59)$$
an ${\bf O^\bfminus_\bfplus(2,2)}$ transformation. Composing the symplectic isomorphisms (3.57) each with symplectomorphisms (of $S$ and $S'$ respectively) is then equivalent to composing the ${\bf O^\bfminus_\bfplus(2,2)}$ transformation (3.59) with ${\bf SO^\bfplus(2,2)}$ transformations. Hence, one obtains the following correspondence: for arbitrary spin frames $\{o^A,\iota^A\}$ and $\{o^{A'},\iota^{A'}\}$,
$$\eqalignno{\pm(\{o^{A'},\iota^{A'}\},\{o^A,\iota^A\})\ &\leftrightarrow\ \hbox{Witt bases of the form    } \{E^a_1 = o^Ao^{A'}, E^a_2 = o^A\iota^{A'}, E^a_3 = \iota^A\iota^{A'}, E^a_4 = -\iota^Ao^{A'}\}\cr
&\leftrightarrow\ \hbox{$\Psi$-ON bases in the ${\bf SO^\bfplus}$-orientation class that is the image}\cr
&\qquad \hbox{ of the standard basis under ${\bf O^\bfminus_\bfplus(2,2)}$,}&(3.60)\cr}$$
compare (3.51b).

Finally, consider the anti-symplectic isomorphisms induced by
$$\check o^A \mapsto \pi^{A'} := \check \iota^{A'} \qquad \check\iota^A \mapsto \xi^{A'} := \check o^{A'} \qquad \check o^{A'} \mapsto \alpha^A := \check \iota^A \qquad \check\iota^{A'} \mapsto \beta^A := \check o^A.\eqno(3.61)$$
These isomorphisms induce
$$\vcenter{\openup1\jot \halign{$\hfil#$&&${}#\hfil$&\qquad$\hfil#$\cr
\check\ell^a &\mapsto E^a_1 := \alpha^A\pi^{A'} = \check n^a & \check{\tilde m}^a &\mapsto E^a_2 := \alpha^A\xi^{A'} = \check{\tilde m}^a\cr
\check n^a &\mapsto E^a_3 := \beta^A\xi^{A'} = \check\ell^a & \check m^a &\mapsto -E^a_4 := \beta^A\pi^{A'} = \check m^a,\cr}}\eqno(3.62)$$
whence
$$\vcenter{\openup1\jot \halign{$\hfil#$&&${}#\hfil$&\qquad$\hfil#$\cr
e^a_1 &\mapsto U^a := \displaystyle{E^a_1 + E^a_3 \over \sqrt2} = e^a_1 & e^a_2 &\mapsto V^a := \displaystyle{E^a_2 + E^a_4 \over \sqrt2} = e^a_2\cr
e^a_3 &\mapsto X^a := \displaystyle{E^a_1 - E^a_3 \over \sqrt2} = -e^a_3 & e^a_4 &\mapsto Y^a := \displaystyle{E^a_2 - E^a_4 \over \sqrt2} = e^a_4,\cr}}\eqno(3.63)$$
an ${\bf O_\bfminus^\bfplus(2,2)}$ transformation. Composing the anti-symplectic isomorphisms (3.61) each with symplectomorphisms (of $S$ and $S'$ respectively) is then equivalent to composing the ${\bf O_\bfminus^\bfplus(2,2)}$ transformation (3.63) with ${\bf SO^\bfplus(2,2)}$ transformations. Hence, one obtains the following correspondence: for arbitrary spin frames $\{o^A,\iota^A\}$ and $\{o^{A'},\iota^{A'}\}$,
$$\eqalignno{\pm(\{\iota^{A'},o^{A'}\},\{\iota^A,o^A\})\ &\leftrightarrow\ \hbox{Witt bases of the form    } \{E^a_1 = \iota^A\iota^{A'}, E^a_2 = \iota^A o^{A'}, E^a_3 = o^Ao^{A'}, E^a_4 = -o^A\iota^{A'}\}\cr
&\leftrightarrow\ \hbox{$\Psi$-ON bases in the ${\bf SO^\bfplus}$-orientation class that is the image}\cr
&\qquad \hbox{ of the standard basis under ${\bf O_\bfminus^\bfplus(2,2)}$,}&(3.64)\cr}$$
compare (3.51b) again. In (3.52--64), Witt bases and $\Psi$-ON bases are related by (3.50). These results justify the identification of the vector representation with the action of the components of {\bf O(2,2)} asserted following (3.46).

Corresponding to PRI (3.3.31), define the four-form
$$e_{abcd} := \epsilon_{AC}\epsilon_{BD}\epsilon_{A'D'}\epsilon_{B'C'} - \epsilon_{AD}\epsilon_{BC}\epsilon_{A'C'}\epsilon_{B'D'},\eqno(3.65)$$
and  put $e^{abcd} := g^{af}g^{bg}g^{ch}g^{di}e_{fghi}$. Then, one computes
$$e_{abcd}e^{abcd} = 24,\eqno(3.66)$$
and for $\Psi$-ON bases $\{U^a,V^a,X^a,Y^a\}$ and their associated Witt bases $\{E_1,E_2,E_3,E_4\}$
$$e_{abcd}U^aV^bX^cY^d = e_{abcd}E^a_1E^b_2E^c_3E^d_4 = \cases{e_{abcd}\ell^an^b\tilde m^cm^d = 1,& for (3.52) and (3.56);\cr -e_{abcd}\ell^an^b\tilde m^cm^d = -1,& for (3.60) and (3.64);\cr}\eqno(3.67)$$
i.e., $e_{abcd}$ is the volume form for the orientation determined by the standard basis. If $\{g_a^{\bf a}:{\bf a} = 1,\ldots,4\,\}$ is the dual basis of $\{U^a,V^a,X^a,Y^a\}$, then
$$e_{abcd} = \pm24 g_{[a}^1g^2_bg^3_cg^4_{d]} = \pm24 U_{[a}V_bX_cY_{d]} \hskip 1.25in e^{abcd} = \pm24U^{[a}V^bX^cY^{d]},\eqno(3.68)$$
with signs according to the orientation class of $\{U^a,V^a,X^a,Y^a\}$.

PRI \S 3.6 show how to derive the vector representation from the fact of the isomorphism (3.37) together with (3.42). I will recapitulate this derivation for the more general case of the group {\bf NO(2)} of transformations of $V^a = {\bf R}^{2,2}$ that are either orthogonal or anti-orthogonal. Note that {\bf O(2,2)} is a subgroup of {\bf NO(2)} of index two, see Law (1992). In particular, if $F_2$ is the anti-orthogonal automorphism of ${\bf R}^{2,2}$ which interchanges $e_1$ with $e_3$ and $e_2$ with $e_4$, then ${\bf NO(2)} = {\bf O(2,2)} \amalg F_2{\bf O(2,2)}$ ($F_2$ is the standard paraHermitian structure on ${\bf R}^{2,2}$).

Each $L^a{}_b \in {\bf NO(2)}$ maps a null vector to a null vector. Thus, for all spinors $\kappa^A \in S^A$ and $\xi^{A'} \in S^{A'}$,
$$L^{AA'}{}_{BB'}\kappa^B\xi^{A'} = \tau^A\eta^{A'},\eqno(3.69)$$
for some spinors $\tau^A$ and $\eta^{A'}$. Exactly as in PRI, p. 169, one concludes that
$$L^a{}_b\hbox{  takes the form  } \theta^A{}_{B'}\mu^{A'}{}_B\hbox{  or  } \phi^A{}_B\mu^{A'}{}_{B'}.\eqno(3.70)$$
Now
$$g_{ab} = \pm g_{cd}L^c{}_aL^d{}_b.\eqno(3.71)$$
Upon substituting the forms (3.70) into (3.71), transvecting by $g^{ab}$, and recalling (3.47), one obtains
$$1 = \cases{\pm\left({1 \over 2}\epsilon^{A'B'}\epsilon_{CD}\theta^C{}_{A'}\theta^D{}_{B'}\right)\left({1 \over 2}\epsilon^{AB}\epsilon_{C'D'}\lambda^{C'}{}_A\lambda^{D'}{}_B\right)\cr \pm\left({1 \over 2}\epsilon^{AB}\epsilon_{CD}\phi^C{}_A\phi^D{}_B\right)\left({1 \over 2}\epsilon^{A'B'}\epsilon_{C'D'}\mu^{C'}{}_{A'}\mu^{D'}{}_{B'}\right)\cr} = \cases{\pm\det(\theta^C{}_{A'})\det(\lambda^{C'}{}_A)\cr \pm\det(\phi^C{}_A)\det(\mu^{C'}{}_{A'})\cr}.$$
If, for example, $\det(\theta^C{}_{A'}) = \rho$, then $\det(\lambda^{C'}{}_A) = \pm\rho^{-1}$, whence dividing $\theta^C{}_{A'}$ by $\sqrt{\vert \rho\vert}$ while multiplying $\lambda^{C'}{}_A$ by the same factor, one can arrange that $\theta^C{}_{A'}$ and $\lambda^{C'}{}_A$ are each of determinant $\pm 1$, and similarly for $\phi^C{}_A$ and $\mu^{C'}{}_{A'}$; assume this is done.
\vskip 24pt
\noindent {\bf 3.72 Proposition}\hfil\break
Let $T \in \End({\bf R}^{2,2})$ have spinor description $\alpha^A{}_B\delta^{A'}{}_{B'}$ or $\beta^{A'}{}_B\gamma^A{}_{B'}$, where $\alpha$, $\delta$, $\beta$, and $\gamma$ each have determinant $\pm 1$. Then $T$ {\sl is} an element of {\bf NO(2)}.

Proof. Choosing spin frames for $S$ and $S'$, the possible spinor forms for the matrix representation of $T$ all take the form $AVB$, $A$, $B$, $V \in {\bf R}(2)$ with $A$ and $B$ of determinant $\pm1$ and $V = (v^{\bf AA'})\ \leftrightarrow\ (v^{\bf a})$ under (3.37--38). It follows that $AVB = 0_2$ iff $V = 0_2$, i.e., $T$ is an automorphism. Also, $\det(AVB) = 0$ iff $\det(V) = 0$. From (3.39), $2\det(V) = s_{2,2}(v,v)$, so $T(v)$ is null iff $v$ is null. Hence, $g(v,w) := s_{2,2}\bigl(T(v),T(w)\bigr)$ is a scalar product on ${\bf R}^{2,2}$ with the same null cone as $s_{2,2}$ itself, whence $g(v,w) = ks_{2,2}(v,w)$ for some $k \in {\bf R}^*$ (e.g., see Beem et al. 1996, p. 26). Thus, $\det(AVB) = k\det(V)$, and $k=\pm1$ according as $\det(A)$ and $\det(B)$ have the same/opposite sign.\bull
\vskip 24pt 
According to whether each of $\theta$, $\lambda$, $\phi$ and $\mu$ has determinant $\pm 1$, eight possibilities occur in (3.70), which one expects to correspond to the eight connected components of {\bf NO(2)}. Indeed, the four components of {\bf O(2,2)} have already been identified above. The anti-orthogonal automorphism $F_2$ interchanges $e_1$ with $e_3$ and $e_2$ with $e_4$, whence fixes $\check\ell^a$ and $\check{\tilde m^a}$ and maps each of $\check n^a$ and $\check m^a$ to their negative. This transformation corresponds to the mapping that fixes $\{\check o^A,\check\iota^A\}$ and the mapping $\kappa^{A'}{}_{B'}$ that maps$\{\check o^{A'},\check\iota^{A'}\}$ to $\{\check o^{A'},-\check\iota^A{'}\}$ (and also corresponds to $(-\epsilon_B{}^A,-\kappa^{A'}{}_{B'})$ of course), i.e.,
$$F_2 = \epsilon_B{}^A\kappa^{A'}{}_{B'} = \epsilon^A{}_B(\check o^{A'}\check\iota_{B'} + \check\iota^{A'}\check o_{B'}).$$
Hence, given the structure of {\bf NO(2)} and (3.52--64), one deduces the following short exact sequences:
$$\vbox{\settabs \+${\bf 0}\ \to\ \{\pm(1,1)\}\ \to$\hskip .2in&$\{\,(\alpha^A{}_B,\delta^{A'}{}_{B'}):\bigl(\det(\alpha),\det(\delta)\bigr) = (1,1)\,\}$\hskip .4in&$\to$\hskip .2in&${\bf SO^\bfplus(2,2)}$\hskip.2in&$\to\ {\bf 0}$\cr
\+${\bf 0}\ \to\ \{\pm(1,1)\}\ \to$\hskip .2in&$\{\,(\alpha^A{}_B,\delta^{A'}{}_{B'}):\bigl(\det(\alpha),\det(\delta)\bigr) = (1,1)\,\}$\hskip .4in&$\to$\hskip .2in&${\bf SO^\bfplus(2,2)}$\hskip.2in&$\to\ {\bf 0}$\cr
\+${\bf 0}\ \to\ \{\pm(1,1)\}\ \to$\hskip .2in&$\{\,(\alpha^A{}_B,\delta^{A'}{}_{B'}):\bigl(\det(\alpha),\det(\delta)\bigr) = (1,-1)\,\}$\hskip .4in&$\to$\hskip .2in&$F_2{\bf SO^\bfplus(2,2)}$\hskip.2in&$\to\ {\bf 0}$\cr
\+${\bf 0}\ \to\ \{\pm(1,1)\}\ \to$\hskip .2in&$\{\,(\alpha^A{}_B,\delta^{A'}{}_{B'}):\bigl((\det(\alpha),\det(\delta)\bigr) = (-1,-1)\,\}$\hskip .4in&$\to$\hskip .2in&${\bf O^\bfminus_\bfminus(2,2)}$\hskip.2in&$\to\ {\bf 0}$\cr
\+${\bf 0}\ \to\ \{\pm(1,1)\}\ \to$\hskip .2in&$\{\,(\alpha^A{}_B,\delta^{A'}{}_{B'}):\bigl(\det(\alpha),\det(\delta)\bigr) = (-1,1)\,\}$\hskip .4in&$\to$\hskip .2in&$F_2{\bf O^\bfminus_\bfminus(2,2)}$\hskip.2in&$\to\ {\bf 0}$\cr
\+${\bf 0}\ \to\ \{\pm(1,1)\}\ \to$\hskip .2in&$\{\,(\beta^{A'}{}_B,\gamma^A{}_{B'}):\bigl(\det(\beta),\det(\gamma)\bigr) = (1,1)\,\}$\hskip .4in&$\to$\hskip .2in&${\bf O^\bfminus_\bfplus(2,2)}$\hskip.2in&$\to\ {\bf 0}$\cr
\+${\bf 0}\ \to\ \{\pm(1,1)\}\ \to$\hskip .2in&$\{\,(\beta^{A'}{}_B,\gamma^A{}_{B'}):\bigl(\det(\beta),\det(\gamma)\bigr) = (1,-1)\,\}$\hskip .4in&$\to$\hskip .2in&$F_2{\bf O^\bfminus_\bfplus(2,2)}$\hskip.2in&$\to\ {\bf 0}$\cr
\+${\bf 0}\ \to\ \{\pm(1,1)\}\ \to$\hskip .2in&$\{\,(\beta^{A'}{}_B,\gamma^A{}_{B'}):\bigl(\det(\beta),\det(\gamma)\bigr) = (-1,-1)\,\}$\hskip .4in&$\to$\hskip .2in&${\bf O_\bfminus^\bfplus(2,2)}$\hskip.2in&$\to\ {\bf 0}$\cr
\+${\bf 0}\ \to\ \{\pm(1,1)\}\ \to$\hskip .2in&$\{\,(\beta^{A'}{}_B,\gamma^A{}_{B'}):\bigl(\det(\beta),\det(\gamma)\bigr) = (-1,1)\,\}$\hskip .4in&$\to$\hskip .2in&$F_2{\bf O_\bfminus^\bfplus(2,2)}$\hskip.2in&$\to\ {\bf 0}$\cr}\eqno(3.73)$$
e.g., the elements of $F_2{\bf O^\bfminus_\bfplus(2,2)}$ are of the form $(\kappa^{A'}{}_{C'}\beta^{C'}{}_B,\epsilon_C{}^A\gamma^C{}_{B'}) = (\kappa^{A'}{}_{C'}\beta^{C'}{}_B,\gamma^A{}_{B'})$, where\hfil\break
$(\beta^{C'}{}_B,\gamma^C{}_{B'}) \in {\bf O^\bfminus_\bfplus(2,2)}$.

From $e_{fghi}L^f{}_aL^g{}_bL^h{}_cL^i{}_d = \det(L)e_{abcd}$, (3.65), and (3.47), one computes that: for $L^a{}_b = \alpha^A{}_B\delta^{A'}{}_{B'}$, $\det(L)e_{abcd} = [\det(\alpha)\det(\delta)]^2e_{abcd}$, i.e., these forms have determinant 1 and determine {\bf SNO(2)}; for $L^a{}_b = \beta^{A'}{}_B\gamma^A{}_{B'}$ one obtains $\det(L)e_{abcd} = -[\det(\beta)\det(\gamma)]^2e_{abcd}$, i.e., $\det(L) = -1$ and these forms determine ${\bf ASNO(2)} := {\bf NO(2)} \setminus {\bf SNO(2)}$. These calculations accord with (3.73). The argument on PRI, p. 170, that the form $\beta^{A'}{}_B\gamma^A{}_{B'}$ cannot be path-connected with the identity transformation $\epsilon_B{}^A\epsilon_{B'}{}^{A'}$ is valid in the present circumstances too.

These results show the spinor formalism is consistent in the sense that the spinorial object corresponding to an $L \in {\bf O(2,2)}$ under the correspondence induced by (3.37) and (3.41) is an element of ${\bf Pin(2,2)}$ and that its vector representation is precisely the action of $L$. 

Now consider elements $L$ of {\bf ASNO(2)}, i.e., {\sl improper} (anti-)orthogonal automorphisms. By (3.72), such $L$ correspond to pairs $\pm(\beta^{A'}{}_B,\gamma^A{}_{B'})$, each term of which is either a symplectic or anti-symplectic isomorphism of symplectic planes, i.e., satisfies
$$\gamma^A{}_{C'}\gamma^B{}_{D'}\epsilon^{C'D'}\epsilon_{AB} = \pm 2 = \beta^{A'}{}_C\beta^{B'}_D\epsilon^{CD}\epsilon_{A'B'}.\eqno(3.74)$$
Thus, $\beta^a := \beta^{A'}{}_C\epsilon^{AC}$ and $\gamma^a := \gamma^A{}_{C'}\epsilon^{A'C'}$ are non-null vectors, of like {\sl character}, i.e., both time like or both space like, iff $L$ is orthogonal, and of opposite character iff $L$ is anti-orthogonal.
\vskip 24pt
\noindent {\bf 3.75 Lemma}\hfil\break
$$\displaylines{L^a{}_b\beta^b = \pm \gamma^a,\qquad\hbox{according as $\beta$ satisfies (3.74) with $\pm $ signs};\cr
L^a{}_b\gamma^b = \pm \beta^a,\qquad\hbox{according as $\gamma$ satisfies (3.74) with $\pm $ signs};\cr
\gamma^a\beta_a = -\tr(\gamma \circ \beta).\cr}$$

Proof. Straightforward computations. For example, $L^a{}_b\gamma^b = \beta^{A'}{}_B\gamma^A{}_{B'}\gamma^{BB'} = \beta^{A'}{}_B\gamma^A{}_{B'}\gamma^B{}_{C'}\epsilon^{B'C'} = \pm\epsilon^{AB}\beta^{A'}{}_B$.\bull
\vskip 24pt
\noindent {\bf 3.76 Lemma}\hfil\break
For any rank-two covariant tensor $H_{ab}$:
$$\displaylines{H_{BAA'B'} = {1 \over 2}(H_{ab} + H_{ba} - H_c{}^cg_{ab} - e_{ab}{}^{cd}H_{cd})\cr
H_{ABB'A'} = {1 \over 2}(H_{ab} + H_{ba} - H_c{}^cg_{ab} + e_{ab}{}^{cd}H_{cd}).\cr}$$
It follows that interchange of $A$ and $B$ in any spinorial object is achieved by
$$U_{ab}{}^{cd} := \epsilon_A{}^D\epsilon_B{}^C\epsilon_{A'}{}^{C'}\epsilon_{B'}{}^{D'} = {1 \over 2}(g_a{}^cg_b{}^d + g_a{}^dg_b{}^c - g_{ab}g^{cd} - e_{ab}{}^{cd})$$
while interchange of $A'$ and $B'$ is achieved by
$$V_{ab}{}^{cd} := \epsilon_A{}^C\epsilon_B{}^D\epsilon_{A'}{}^{D'}\epsilon_{B'}{}^{C'} = {1 \over 2}(g_a{}^cg_b{}^d + g_a{}^dg_b{}^c - g_{ab}g^{cd} + e_{ab}{}^{cd})$$
It follows that, for example,
$$\displaylines{U^{ab}{}_{cd} = U_{cd}{}^{ab} \hskip 1.25in V^{ab}{}_{cd} = V_{cd}{}^{ab}\cr
U_{ab}{}^{cd} = V_{ab}{}^{dc} = V_{ba}{}^{cd}.\cr}$$

Proof. By exactly the same argument as on PRI, p. 153, these results are analogues of PRI (3.4.53--54), (3.4.57), (3.4.59--60).\bull
\vskip 24pt
\noindent {\bf 3.77 Lemma}\hfil\break
For $L \in {\bf ANO(2)}$,
$$L_{ab} = \beta_{(a}\gamma_{b)} - {1 \over 2}(\beta^c\gamma_c)g_{ab} + {1 \over 2}e_{abcd}\gamma^c\beta^d.$$

Proof. Since $L_{ab} = \gamma_{AB'}\beta_{A'B}$, then $L_{ab} = V_{ab}{}^{cd}\gamma_c\beta_d$ and the result follows.\bull
\vskip 24pt
\noindent {\bf 3.78 Observation}\hfil\break
It is well known that $L \in \End_\pm({\bf R}^{p,q})$, i.e., is self/skew-adjoint, iff $L_{ab}$ is symmetric/skew:
$$s_{p,q}(L\ ,\ ) =\pm s_{p,q}(\ ,L\ )\ \iff\ g_{cb}L^c{}_a = \pm g_{ac}L^c{}_b,\qquad\hbox{i.e., } L_{ba} = \pm L_{ab}.\eqno(3.78.1)$$
It is of interest to identify (anti-)orthogonal automorphisms $L^a{}_b$ for which $L_{ab}$ is symmetric or skew. By (3.78.1), if $L_{ab}$ is symmetric, then $L^a{}_b$ is orthogonal/anti-orthogonal iff $L^2 = \pm 1$ i.e., an involution/complex structure, while if $L_{ab}$ is skew, then $L^a{}_b$ is orthogonal/anti-orthogonal iff $L^2 = \mp 1$.

An involution $F$ of $V = {\bf R}^{p,q}$ induces a decomposition into its two eigenspaces $V = E_1(F) \oplus E_{-1}(V)$.  If $F$ is also orthogonal, then in fact $V = E_1(F) \operp E_{-1}(F)$. If $F$ is anti-orthogonal, then each of the two eigenspaces must be totally null, whence each is of dimension at most $\dim(V)/2$. Clearly, $\det(F) = (-1)^k$, where $k = \dim\bigl(E_{-1}(F)\bigr)$. If $F$ is improper, then $k$ is odd. Hence for $V = {\bf R}^{2,2}$, if $F$ were improper and anti-orthogonal, then one of $E_1(F)$ or $E_{-1}(F)$ would be a three-dimensional totally null subspace, which is impossible; thus ${\bf R}^{2,2}$ admits no improper, anti-orthogonal, involutions. Anti-orthogonal involutions for which the eigenspaces are of equal dimension are called {\sl paraHermitian structures}. 

All complex structures have positive determinant, so there are no improper, (anti)-orthogonal, anti-involutions. Only ${\bf R}^{2p,2q}$ admits orthogonal complex structures $J$; the $\pm i$-eigenspaces of $J_{\bf C}$ are totally null in ${\bf C}^{2p,2q}$. ${\bf R}^{n,n}$ admits anti-orthogonal complex structures $J$; the $\pm i$-eigenspaces in ${\bf C}^{n,n}$ are orthogonal to each other.
\vskip 24pt
\noindent {\bf 3.79 Proposition}\hfil\break
Let $L \in {\bf ASNO(2)}$. If $L_{ab}$ is symmetric, then $L^a{}_b$ is orthogonal, whence involutary, $\beta^a = \pm \gamma^a$, and $L$ takes one of two forms (signs according as $\beta^a = \pm \gamma^a$):\hfil\break
(i) $\gamma^a$ is time-like; $L_{ab} = \pm\gamma_a\gamma_b \mp g_{ab}$; $L:\gamma^a \mapsto \beta^a$ and $:\beta^a \mapsto \gamma^a$, whence is reflexion in the line $\langle \gamma^a \rangle_{\bf R}$/in the hyperplane $\langle \gamma^a \rangle_{\bf R}^\perp$ according as $\beta^a = \pm \gamma^a$;\hfil\break
(ii) $\gamma^a$ is space-like; $L_{ab} = \pm\gamma_a\gamma_b \pm g_{ab}$; $L:\gamma^a \mapsto -\beta^a$ and $:\beta^a \mapsto -\gamma^a$, whence is reflexion in the hyperplane $\langle \gamma^a \rangle_{\bf R}^\perp$/in the line $\langle \gamma^a \rangle_{\bf R}$ according as $\beta^a = \pm \gamma^a$.

As already noted in (3.78), there are no involutary, improper anti-orthogonal automorphisms of ${\bf R}^{2,2}$ nor improper, (anti-)orthogonal complex structures.

Proof. Suppose $L_{ab} = L_{(ab)}$. From the spinor form of improper (anti-)orthogonal automorphisms and PRI (3.5.2), or (3.77), it follows that $\beta^a \propto \gamma^a$. Given the normalization condition (3.74), $\beta^a = \pm \gamma^a$ (and are non-null). Substitution into (3.77) gives the two forms of $L_{ab}$ according to the character of $\gamma^a$. From these two forms, one easily deduces the action of $L$ on $\langle \gamma^a \rangle_{\bf R}^\perp$, while (3.75) gives the action on $\gamma^a$ itself. Since $\beta^a$ and $\gamma^a$ are of the same character, i.e., $\beta^{A'}{}_B$ and $\gamma^A{}_{B'}$ have the same determinant, $L$ is necessarily orthogonal. Hence, there are no improper, anti-orthogonal, anti-involutions.

When $L_{ab} = L_{[ab]}$, (3.77) gives
$$\gamma_{(a}\beta_{b)} = {1 \over 2}(\gamma^c\beta_c)g_{ab}.\eqno(*)$$
Transvection of this equation by $g^{ab}$ shows that $\gamma^c\beta_c=0$; whence $(*)$ becomes $\gamma_{(a}\beta_{b)} = 0$, and transvecting this equation by $\gamma^b$ yields $\beta_a(\gamma^b\gamma_b) = 0$. By (3.74), $\gamma^a$ is non-null, so $\beta^a = 0$, which is impossible, as expected.\bull
\vskip 24pt
Now consider proper (anti-)orthogonal automorphisms of ${\bf R}^{2,2}$, i.e., elements of {\bf SNO(2)}:
$$L^a{}_b = \alpha^A{}_B\delta^{A'}{}_{B'}, \hskip 1in \alpha^A{}_C\alpha^B{}_D\epsilon^{CD} = \pm\epsilon^{AB}, \hskip 1in \delta^{A'}{}_{C'}\delta^{B'}{}_{D'}\epsilon_{C'D'} = \pm\epsilon_{A'B'},\eqno(3.80)$$
whence
$$\alpha_{AB}\alpha^{AB} = \pm2,\hbox{  according as  } \det(\alpha) = \pm1\hskip .5in \delta_{A'B'}\delta^{A'B'} = \pm2.\hbox{  according as  } \det(\delta) = \pm1.\eqno(3.81)$$
Writing
$$\vcenter{\openup1\jot \halign{$\hfil#$&&${}#\hfil$&\qquad$\hfil#$\cr
\alpha_{AB} &= \psi_{AB} + \alpha\epsilon_{AB}; & \psi_{AB} &= \psi_{(AB)},&\alpha &:= \displaystyle{1 \over 2}\alpha_A{}^A,\cr
\delta_{A'B'} &= \chi_{A'B'} + \delta\epsilon_{A'B'}; & \chi_{A'B'} &= \chi_{(A'B')}, & \delta &:= \displaystyle{1 \over 2}\delta_{A'}{}^{A'}\cr}}\eqno(3.82)$$
then
$$L_{ab} = pg_{ab} + F_{ab} + T_{ab},\eqno(3.83)$$
where
$$p := \alpha\delta, \hskip .5in F_{ab} := \delta\psi_{AB}\epsilon_{A'B'} + \alpha\chi_{A'B'}\epsilon_{AB}, \hskip .5in T_{ab} := \psi_{AB}\chi_{A'B'},\eqno(3.84)$$
with $p \in {\bf R}$, $F_{ab}$ skew, and $T_{ab}$ symmetric and trace free. If one regards $F_{ab}$ as a Maxwell tensor, its energy tensor is $(p/2\pi)T_{ab}$.

Further analysis depends on the decomposition of the symmetric spinors $\psi_{AB}$ and $\chi_{A'B'}$. For $\psi_{AB}$ there are two possibilities:
$$\psi_{AB} = \mu_A\nu_B + \nu_A\mu_B\qquad\hbox{or}\qquad\psi_{AB} = \pm\kappa_{(A}\bar\kappa_{B)} = \pm(\rho_A\rho_B + \sigma_A\sigma_B),\eqno(3.85)$$
where $\mu_A$, $\nu_A$, $\rho_A$, and $\sigma_A$ are real spinors and $\kappa_A = \rho_A + i\sigma_A$ and $\sigma^A\rho_A \not= 0$ (which condition ensures that $\kappa_A$ is nontrivially complex, i.e., not just a complex scalar multiple of a real spinor),
with
$$\psi_{AB}\psi^{AB} = \cases{-2(\nu^A\mu_A)^2 \leq 0;&\cr 2(\sigma^A\rho_A)^2 > 0;&\cr}\eqno(3.86)$$
respectively. Substitution of these forms into (3.81) yields
$$\alpha^2 - (\nu^A\mu_A)^2 = \pm1\qquad\hbox{and}\qquad \alpha^2 + (\sigma^A\rho_A)^2 = \pm1,\eqno(3.87)$$
respectively. The same possibilities occur for $\chi_{A'B'}$, of course, for which there is an analogous development to that which follows.
\vskip 24pt
\noindent {\bf Lemma 3.88}\hfil\break
The quantity $\alpha_{AB}$ in (3.80) takes one of three forms.\hfil\break
$R_2$: with $\mu^D\nu_D \not=0$, $\det(\alpha) = \pm1$ (equivalently, $\alpha_{AB}\alpha^{AB} = \pm2$), and $\psi_{AB} = \mu_A\nu_B + \nu_A\mu_B$ then
$$\alpha_{AB} = \mu_A\nu_B\left(1 + {\alpha \over (\nu^D\mu_D)}\right) + \nu_A\mu_B\left(1 - {\alpha \over (\nu^D\mu_D)}\right).\eqno(3.88{\rm a})$$
$CC$: with $\sigma^D\rho_D \not = 0$, $\kappa_A = \rho_A + i\sigma_A$, $\det(\alpha) = 1$ (equivalently, $\alpha_{AB}\alpha^{AB} = 2$), and $\psi_{AB} = \pm\kappa_{(A}\bar\kappa_{B)}$ then
$$\alpha_{AB} = \rho_A\left(\pm\rho_B + {\alpha \over (\sigma^D\rho_D)}\sigma_B\right) + \sigma_A\left(\pm\sigma_B - {\alpha \over (\sigma^D\rho_D)}\rho_B\right) = \pm\kappa_{(A}\bar\kappa_{B)} + \alpha\epsilon_{AB}.\eqno(3.88{\rm b})$$
$R_1$: with $\mu^D\nu_D = 0$, then $\det(\alpha) = 1$ (equivalently, $\alpha_{AB}\alpha^{AB} = 2$), and $\psi_{AB} \propto \mu_A\mu_B$ then
$$\alpha_{AB} = \pm(\zeta\mu_A\mu_B + \epsilon_{AB}),\eqno(3.88{\rm c})$$
where $\mu_A$, $\nu_A$, $\rho_A$ and $\sigma_A$ are real spinors. The notation $R_2$ indicates that $\psi_{AB}$ has two real distinct principal spinors (PSs), $CC$ indicates that $\psi_{AB}$ has a pair of complex conjugate PSs, while $R_1$ indicates that $\psi_{AB}$ has a multiple real PS.

Proof. If $\nu^D\mu_D \not= 0$, then
$$\epsilon_{AB} = {1 \over \nu^D\mu_D}(\mu_A\nu_B - \nu_A\mu_B) = {1 \over \sigma^D\rho_D}(\rho_A\sigma_B - \sigma_A\rho_b).$$
Combining the relevant form with (3.85) and substituting into (3.82) yields the cases $R_2$ and $CC$. Note that in the case $CC$, (3.87) implies $\det(\alpha) = 1$. For the case $R_1$, write $\nu_A = (\zeta/2)\mu_A$, so that $\psi_{AB} = \zeta\mu_A\mu_B$. From (3.87), $\alpha^2 = \pm1$, but as $\alpha$ is real then $\alpha^2 = 1$, and $\alpha_{AB} = \zeta\mu_A\mu_B \pm \epsilon_{AB}$, whence $\alpha_{AB}\alpha^{AB} = 2$, i.e., $\det(\alpha) = 1$. By absorbing a sign into $\zeta$ if necessary, one can write $\alpha_{AB} = \pm(\zeta\mu_A\mu_B + \epsilon_{AB})$.\bull
\vskip 24pt
\noindent{\bf Lemma 3.89}\hfil\break
The endomorphism $\alpha^A{}_B$ has the following eigenvalues and eigenvectors according to the three cases of (3.88):\hfil\break
$R_2$: $\alpha^A{}_B$ has two distinct real eigenvalues each of geometric multiplicity one, $\lambda_1 := -(\nu^D\mu_D + \alpha)$ with real eigenvector $\mu^A$, and $\lambda_2 := \nu^D\mu_D - \alpha$ with eigenvector $\nu^A$;\hfil\break
$CC$: $\alpha^A{}_B$ has a complex conjugate pair of eigenvalues, each of geometric multiplicity one, $\lambda_1:= \pm i(\sigma^D\rho_D) - \alpha$ with complex eigenvector $\kappa^A$, and $\lambda_2 := \mp i(\sigma^D\rho_D) - \alpha$ with complex eigenvector $\bar\kappa^A$;\hfil\break
$R_1$: if $\zeta\not= 0$, $\alpha^A{}_B = \pm(\zeta\mu^A\mu_B + \epsilon^A{}_B)$ has a single real eigenvalue $\mp1$, of geometric multiplicity one, with real eigenvector $\mu_A$.

Thus, in each case, the eigenvectors are just the PSs of $\psi_{AB}$. Note that for cases $R_2$ and $CC$, (3.87) can be interpreted as the fact that the determinant is the product of the eigenvalues; in particular, $\lambda_2 = \pm\lambda^{-1}$ in these cases, with sign according as $\det(\alpha^A{}_B) = \pm 1$. In case $R_1$, the single eigenvalue $\mp1$ has algebraic multiplicity 2, implying $\det(\alpha^A{}_B) = 1$, consistent with (3.88).

Proof. Since $\alpha^A{}_B$ has at most two distinct eigenvalues, the results follow from simple calculations using the forms in (3.88); namely, compute the action of $\alpha^A{}_B$ on the PSs of $\psi_{AB}$. For case $R_1$, note that $\alpha^A{}_B\gamma^B = \lambda\gamma^A$ is equivalent to $\pm\zeta(\gamma^D\mu_D)\mu^A \mp \gamma^A = \lambda\gamma^A$, i.e., to $(\lambda \pm 1)\gamma^A = \pm\zeta(\gamma^D\mu_D)\mu^A$. If $\lambda \not= \mp1$, $\gamma^A \propto \mu^A$; if $\lambda = \mp1$, then $\gamma^D\mu_D = 0$.\bull
\vskip 24pt
\noindent {\bf 3.90 Notation}\hfil\break
For any endomorphism $T \in \End(V)$, $V \cong {\bf K}^n$ a {\bf K}-linear space, and (real or complex) eigenvalue $\lambda$, let $m_\lambda(T)$ denote the algebraic multiplicity of $T$, $E_\lambda(T)$ the eigenspace of $\lambda$ (over $V$), and $M_\lambda(T)$ the ({\bf K}-)dimension of $E_\lambda(T)$ (the {\sl geometric multiplicity}). When ${\bf K} = {\bf R}$, let $T_{\bf C}$ denote the complexification of $T$, i.e., the {\bf C}-linear extension to an element of $\End({\bf C}V)$; then $E_\lambda(T_{\bf C})$ is the eigenspace of $\lambda$ for $T_{\bf C}$ over ${\bf C}V$ and $M_\lambda(T_{\bf C})$ its (complex) dimension. By analogy with ${\bf R}^{p,q}$, introduce the notation ${\bf C}^{p,q}$, i.e., ${\bf C}^{p+q}$ equipped with the scalar product $s_{p,q}$ that makes the standard basis of ${\bf C}^{p+q}$ $\Psi$-ON with signature $(p,q)$. Similarly, $\overline{\bf C}^{p,q}$ is ${\bf C}^{p+q}$ equipped with the {\sl Hermitian} scalar product $h_{p,q}$ that makes the standard basis of ${\bf C}^{p+q}$ $\Psi$-unitary with signature $(p,q)$. Of course, all the ${\bf C}^{p,q}$ of a given dimension are isometric to each other and signature is not an invariant of these spaces; the notation ${\bf C}^{p,q}$ is for convenience only, as the complexification of ${\bf R}^{p,q}$.
\vskip 24pt
The following results are easy and no doubt well known.
\vskip 24pt
\noindent {\bf 3.91 Proposition}\hfil\break
Let $L \in {\bf O(p,q)}$.\hfil\break
(1) For a real eigenvalue $\lambda \not= \pm 1$, $E_\lambda(L)$ is totally null in ${\bf R}^{p,q}$.\hfil\break
(2) For an eigenvalue $\lambda$ of $L_{\bf C}$, $\lambda \not= \pm 1$ implies $E_\lambda(L_{\bf C})$ is totally null in ${\bf C}^{p,q}$ and $\vert \lambda \vert \not= 1$ implies $E_\lambda(L_{\bf C})$ is totally null in $\overline{\bf C}^{p,q}$.\hfil\break
(3) In particular, for $L \in {\bf O(n)}$, all eigenvalues are unimodular.\hfil\break
(4) Let $\lambda$ and $\mu$ be eigenvalues of $L$. If $\lambda\mu \not= 1$, then $E_\lambda(L_{\bf C})$ and $E_\mu(L_{\bf C})$ are orthogonal in ${\bf C}^{p,q}$ (if both $\lambda$ and $\mu$ are real then $E_\lambda(L)$ and $E_\mu(L)$ are orthogonal in ${\bf R}^{p,q}$), while if $\lambda\overline\mu \not= 1$ then $E_\lambda(L_{\bf C})$ and $E_\mu(L_{\bf C})$ are orthogonal in $\overline{\bf C}^{p,q}$.\hfil\break
(5) If $\lambda \not= \pm 1$ is an eigenvalue of $L$, then $\lambda^{-1}$ is a distinct eigenvalue and
$$M_{\lambda^{-1}}(L) = M_\lambda(L) \hskip 1.5in m_{\lambda^{-1}}(L) = m_\lambda(L).$$
Thus, the eigenvalues, other than $\pm 1$, of an orthogonal automorphism occur in pairs $\lambda$ and $\lambda^{-1}$ having the same multiplicities, while eigenspaces of distinct eigenvalues $\lambda$ and $\mu$ are orthogonal, except perhaps when $\mu = \lambda^{-1}$.
\vskip 24pt
Similarly, one can easily prove the following analogue.
\vskip 24pt
\noindent {\bf 3.92 Proposition}\hfil\break
Let $L \in {\bf AO(n,n)}$.\hfil\break 
(1) if $\lambda$ is a real eigenvalue of $L$ then $E_\lambda(L)$ is totally null in ${\bf R}^{n,n}$,\hfil\break
(2) if $\lambda$ is an eigenvalue of $L_{\bf C}$, then $E_\lambda(L_{\bf C})$ is totally null in $\overline{\bf C}^{n,n}$ and totally null in ${\bf C}^{n,n}$ provided $\lambda^2  \not= -1$,\hfil\break
(3) for eigenvalues $\lambda$ and $\mu$ of $L$, if $\lambda\mu \not = -1$ then $E_\lambda(L_{\bf C})$ and $E_\mu(L_{\bf C})$ are orthogonal in ${\bf C}^{n,n}$ (if both eigenvalues are real then $E_\lambda(L)$ and $E_\mu(L)$ are orthogonal in ${\bf R}^{n,n}$), while if $\lambda\overline\mu \not= -1$ then $E_\lambda(L_{\bf C})$ and $E_\mu(L_{\bf C})$ are orthogonal in $\overline{\bf C}^{n,n}$,\hfil\break
(4) if $\lambda \not= \pm i$ is an eigenvalue of $L$ then $-\lambda^{-1}$ is a distinct eigenvalue of $L$ and
$$M_{-\lambda^{-1}}(L) = M_\lambda(L) \hskip 1.5in m_{-\lambda^{-1}}(L) = m_\lambda(L).$$
\vskip 24pt
\noindent Thus, the eigenvalues, other than $\pm i$, of an anti-orthogonal automorphism occur in pairs $\lambda$ and $-\lambda^{-1}$ having the same multiplicities, while eigenspaces of distinct eigenvalues $\lambda$ and $\mu$ are orthogonal, except possibly when $\mu = -\lambda^{-1}$.
\vskip 24pt
It is easy to see that for $L \in {\bf SNO(2)}$, $L^a{}_b = \alpha^A{}_B\delta^{A'}{}_{B'}$ has a (real or complex) null eigenvector $\xi^A\chi^{A'}$ iff $\xi^A$ is an eigenvector of $\alpha^A{}_B$ and $\chi^{A'}$ of $\delta^{A'}{}_{B'}$ (where $\xi^A$ and $\chi^{A'}$ may be real or complex). Hence, by applying (3.89) and its analogue for $\delta^{A'}{}_{B'}$, one obtains a characterization of the null eigenvectors of $L \in {\bf SNO(2)}$; specifically, for each of the four components of {\bf SNO(2)}, according as $\det(\alpha^A{}_B) = \pm 1$ and $\delta^{A'}{}_{B'} = \pm 1$, one can apply the applicable cases $R_2$, $CC$, and $R_1$ to each of $\alpha^A{}_B$ and $\delta^{A'}{}_{B'}$ (cases $R_1$ and $CC$ not being applicable when the determinant is $-1$). Thus, there are nine cases for ${\bf SO^\bfplus(2,2)}$, one for ${\bf O^\bfminus_\bfminus(2,2)}$, and three each for $F_2{\bf SO^\bfplus(2,2)}$ and $F_2{\bf O^\bfminus_\bfminus(2,2)}$.
\vskip 24pt
\noindent {\bf 3.93 Remark}\hfil\break
In the case $CC$, as $\psi_{AB}$ is symmetric in $\kappa_A$ and $\bar\kappa_A$, one has the freedom to specify that $\kappa_A = \rho_A + i\sigma_A$ satisfies $\sigma^D\rho_D > 0$ (if not, interchange $\kappa_A$ and $\bar\kappa_A$). Similarly, in case $R_2$, one has the freedom to specify that $\nu^D \mu_D > 0$ (otherwise interchange $\mu_A$ and $\nu_A$).
\vskip 24pt
Now consider $L^a{}_b \in {\bf SNO(2)}$ for which $L_{ab}$ is symmetric or skew.
\vskip 24pt
\noindent {\bf 3.94 Proposition}\hfil\break
Other than multiples of the identity, $L^a{}_b = \alpha^A{}_B\delta^{A'}{}_{B'} \in {\bf SNO(2)}$ has $L_{ab}$ symmetric iff:\hfil\break
i)  $L^a{}_b \in {\bf SO^\bfplus(2,2)}$ with $L^a{}_b = \psi^A{}_B\chi^{A'}{}_{B'}$, where $\psi_{AB} = \pm\kappa_{(A}\bar\kappa_{B)}$ and $\chi_{A'B'} = \pm\lambda_{(A'}\bar\lambda_{B')}$ (signs not paired), $\kappa = \rho_A + i\sigma_A$ and $\lambda_{A'} = \eta_{A'} + i\theta_{A'}$ both nontrivially complex spinors such that $\{\rho^A,\sigma^A\}$ and $\{\eta^{A'},\theta^{A'}\}$ are spin frames. $L^a{}_b = \pm(\rho^A\rho_B + \sigma^A\sigma_B)(\eta^{A'}\eta_{B'} + \theta^{A'}\theta_{B'})$; for the given spin frames, if $\{U^a,V^a,X^a,Y^a\}$ is the associated $\Psi$-ON basis (3.48--49), then $\pm L = 1 \oplus (-1)$ with respect to $\langle U^a,V^a \rangle_{\bf R} \operp \langle X^a,Y^a \rangle_{\bf R}$, i.e., $\pm L$ is reflexion in $\langle U^a,V^a \rangle_{\bf R} \cong {\bf R}^{2,0}$; or\hfil\break
ii) $L^a{}_b \in {\bf O^\bfminus_\bfminus(2,2)}$ with $L^a{}_b = \psi^A{}_B\chi^{A'}{}_{B'}$, $\psi_{AB} = \mu_A\nu_B + \nu_A\mu_B$, $\chi_{A'B'} = \pi_{A'}\tau_{B'} + \tau_{A'}\pi_{B'}$, with $\{\mu^A,\nu^A\}$ and $\{\pi^{A'},\tau^{A'}\}$ spin frames. If $\{U^a,V^a,X^a,Y^a\}$ is the associated $\Psi$-ON basis (3.48--49), then $L = 1 \oplus (-1)$ with respect to $\langle U^a,X^a \rangle_{\bf R} \operp \langle V^a,Y^a \rangle_{\bf R}$, i.e., $L$ is reflexion in $\langle U^a,X^a \rangle_{\bf R} \cong {\bf R}^{1,1}$; or\hfil\break
iii) $\pm L^a{}_b \in F_2{\bf O^\bfminus_\bfminus(2,2)}$ with $L^a{}_b = \psi^A{}_B\chi^{A'}{}_{B'}$, where $\psi_{AB} = \mu_A\nu_B + \nu_A\mu_B$ and $\chi_{A'B'} = \lambda_{(A'}\bar\lambda_{B')} = (\eta_{A'}\theta_{B'} + \theta_{A'}\eta_{B'})$, $\lambda_{A'} = \eta_{A'} + i\theta_{A'}$, and $\{\mu^A\,\nu^A\}$ and $\{\eta^{A'},\theta^{A'}\}$ are spin frames. $\pm L^a{}_b$ are a complex-conjugate pair of anti-orthogonal complex structures; if $\{U^a,V^a,X^a,Y^a\}$ is the $\Psi$-ON basis associated to the spin frames by (3.48--49), $\{U^a,V^a,Y^a,X^a\}$ is a complex basis for $L^a{}_b$, i.e., $L(U) = Y$ and $L(V) = X$, whence $\langle U^a - iY^a, V^a - iX^a \rangle_{\bf C}$ is the $i$-eigenspace of $L_{\bf C}$, and is nondegenerate and orthogonal to the $-i$-eigenspace $\langle U^a + iY^a, V^a + iX^a \rangle_{\bf C}$ in ${\bf C}^{2,2}$; or\hfil\break
iv) $\pm L^a{}_b \in F_2{\bf SO^\bfplus(2,2)}$ with $L^a{}_b = \psi^A{}_B\chi^{A'}{}_{B'}$, where $\psi_{AB} = \kappa_{(A}\bar\kappa_{B)} = (\rho_A\rho_B + \sigma_A\sigma_B)$, $\kappa_A = \rho_A + i\sigma_A$, and $\chi_{A'B'} = \pi_{A'}\tau_{B'} + \tau_{A'}\pi_{B'}$, and $\{\rho^A,\sigma^A\}$ and $\{\eta^{A'},\theta^{A'}\}$ are spin frames. $\pm L^a{}_b$ are a complex-conjugate pair of anti-orthogonal complex structures; if $\{U^a,V^a,X^a,Y^a\}$ is the $\Psi$-ON basis associated to the spin frames by (3.48--49), then $\{U^a,X^a,Y^a,V^a\}$ is a complex basis for $L^a{}_b$, whence $\langle U^a - iY^a, X^a - iV^a\rangle_{\bf C}$ is the $i$-eigenspace of $L_{\bf C}$, and is nondegenerate and orthogonal to the $-i$-eigenspace $\langle U^a + iY^a,X^a + iV^a \rangle_{\bf C}$ in ${\bf C}^{2,2}$.

$L^a{}_b = \alpha^A{}_B\delta^{A'}{}_{B'} \in {\bf SNO(2)}$ has $L_{ab}$ skew iff:\hfil\break
a) $\pm L^a{}_b \in {\bf SO^\bfplus(2,2)}$, with $L^a{}_b = \psi^A{}_B\epsilon^{A'}{}_{B'}$, where $\psi_{AB} = \kappa_{(A}\bar\kappa_{B)} = \rho_A\sigma_B + \sigma_A\rho_B$, $\kappa_A = \rho_A + i\sigma_A$, and $\{\rho^A,\sigma^A\}$ is a spin frame. $\pm L^a{}_b$ is a pair of complex-conjugate orthogonal complex structures. The $i$-eigenspace of $L_{\bf C}$ is the complex $\beta$-plane ${\bf W}_{[\bar\kappa]} := \{\,\bar\kappa^A\zeta^{A'}:\zeta^{A'} \in {\bf C}S'\,\}$ in ${\bf C}^{2,2}$ and the $-i$-eigenspace is the complex $\beta$-plane ${\bf W}_{[\kappa]} := \{\,\kappa^A\zeta^{A'}:\zeta^{A'} \in {\bf C}S'\,\}$, which are each totally null but not orthogonal to each other and ${\bf C}^{2,2} = {\bf W}_{[\bar\kappa]} \oplus {\bf W}_{[\kappa]}$; or\hfil\break
b) $\pm L^a{}_b \in F_2{\bf O^\bfminus_\bfminus(2,2)}$, with $L^a{}_b = \psi^A{}_B\epsilon^{A'}{}_{B'}$, where $\psi_{AB} = \mu_A\nu_B + \nu_A\mu_B$ and $\{\mu^A,\nu^A\}$ is a spin frame. $\pm L^a{}_b$ are anti-orthogonal involutions; in particular, paraHermitian structures. The eigenspaces are real $\beta$-planes in ${\bf R}^{2,2}$, $E_1(L) = W_{[\mu]} := \{\,\mu^A\omega^{A'}:\omega^{A'} \in S'\,\}$ and $E_{-1}(L) = W_{[\nu]} := \{\,\nu^A\omega^{A'}:\omega^{A'} \in S'\,\}$, which are each totally null but not orthogonal to each other and ${\bf R}^{2,2} = W_{[\mu]} \oplus W_{[\nu]}$; or\hfil\break
c) $\pm L^a{}_b \in {\bf SO^\bfplus(2,2)}$, with $L^a{}_b = \epsilon^A{}_B\chi^{A'}{}_{B'}$, where $\chi_{A'B'} = \lambda_{(A'}\bar\lambda_{B')} = \eta_{A'}\eta_{B'} + \theta_{A'}\theta_{B'}$, $\lambda_{A'} = \eta_{A'} + i\theta_{A'}$, and $\{\eta^{A'},\theta^{A'}\}$ is a spin frame. $\pm L^A{}_B$ is a pair of complex-conjugate orthogonal complex structures. The $i$-eigenspace of $L_{\bf C}$ is the complex $\alpha$-plane ${\bf Z}_{[\bar\lambda]} := \{\,\zeta^A\bar\lambda^{A'}:\zeta^A \in {\bf C}S\,\}$ in ${\bf C}^{2,2}$ and the $-i$-eigenspace is the complex $\alpha$-plane ${\bf W}_{[\lambda]} := \{\,\zeta^A\lambda^{A'}:\zeta^A \in {\bf C}S\,\}$, which are each totally null but not orthogonal to each other and ${\bf C}^{2,2} = {\bf W}_{[\bar\lambda]} \oplus {\bf W}_{[\lambda]}$; or\hfil\break
d) $\pm L^a{}_b \in F_2{\bf SO^\bfplus(2,2)}$, with $L^a{}_b = \epsilon^A{}_B\chi^{A'}{}_{B'}$, where $\chi_{A'B'} = \pi_{A'}\tau_{B'} + \tau_{A'}\pi_{B'}$ and $\{\pi^{A'},\tau^{A'}\}$ is a spin frame. $\pm L^a{}_b$ are anti-orthogonal involutions; in particular paraHermitian structures. The eigenspaces are real $\alpha$-planes in ${\bf R}^{2,2}$, $E_1(L) = Z_{[\pi]} := \{\,\omega^A\pi^{A'}:\omega^A \in S\,\}$ and $E_{-1}(L) = Z_{[\tau]} := \{\,\omega^A\tau^{A'}:\omega^A \in S\,\}$, which are each totally null but not orthogonal to each other and ${\bf R}^{2,2} = Z_{[\pi]} \oplus Z_{[\tau]}$.

$F_2$, of course, must be an example of (d); $F_2 = \epsilon^A{}_B(\check o^{A'}\check\iota_{B'} + \check\iota^{A'}\check o_{B'})$, as noted in connexion with (3.73).

Note that all orthogonal complex structures belong to ${\bf SO^\bfplus(2,2)}$. In fact, for the standard basis $\{e_1,\ldots,e_n\}$ of ${\bf R}^{2p,2q}$, the mapping $J:e_i\mapsto e_{p+i}$, $e_{p+i} \mapsto -e_i$, for $i=1,\ldots,p$, and $e_{2p+j} \mapsto e_{2p+q+j}$, $e_{2p+q+j} \mapsto -e_{2p+j}$, for $j=1,\ldots,q$ is an orthogonal complex structure and an element of ${\bf SO^\bfplus(2p,2q)}$. Every other orthogonal complex structure of ${\bf R}^{2p,2q}$ is of the form $LJL^{-1}$, for some $L \in {\bf O(2p,2q)}$, whence in ${\bf SO^\bfplus(2p,2q)}$, as the latter is a normal subgroup.

Proof.  $L \in {\bf SNO(2)}$ has $L_{ab} = L_{(ab)}$ iff $F_{ab} = 0$ in (3.83), i.e., by (3.84), iff $\delta\psi_{AB} = 0$ and $\alpha\chi_{A'B'} = 0$. Now $\delta$ and $\chi_{A'B'}$ zero is equivalent to $\delta^{A'}{}_{B'} = 0$, while $\alpha$ and $\psi_{AB}$ zero is equivalent to $\alpha^A{}_B = 0$; in both cases $L = 0$, which is not (anti-)orthogonal. If $\psi_{AB} =0$ and $\chi_{A'B'} = 0$, then $L_{ab} = pg_{ab}$ in (3.83) and $L = \pm 1$. Thus, the only non-trivial possibility is $\alpha=\delta=0$. Substituting this condition into (3.87) gives four possibilities, as follows.

Case $CC \times CC$, i.e., $\psi_{AB} = \pm\kappa_{(A}\bar\kappa_{B)}$, $\kappa_A = \rho_A + i\sigma_A$, and $\chi_{A'B'} = \pm\lambda_{(A'}\bar\lambda_{B')}$, $\lambda_A = \eta_{A'} + i\theta_{B'}$, with $(\sigma^D\rho_D)^2 = (\theta^{D'}\eta_{D'})^2 = 1$ by (3.83), whence, by virtue of (3.93), one can suppose $\sigma^D\rho_D = 1$ and $\theta^{D'}\eta_{D'} = 1$, i.e., $\{\rho^A,\sigma^A\}$ and $\{\eta^{A'},\theta^{A'}\}$ are spin frames. In this case $L^a{}_b \in {\bf SO^\bfplus(2,2)}$, and one computes that, under $\pm L$:
$$\vcenter{\openup1\jot \halign{$\hfil#$&&${}#\hfil$&\qquad$\hfil#$\cr
\rho^A\eta^{A'} &\mapsto (\sigma^D\rho_D)(\theta^{D'}\eta_{D'})\sigma^A\theta^{A'} = \sigma^A\theta^{A'} & \sigma^A\theta^{A'} &\mapsto (\sigma^D\rho_D)(\theta^{D'}\eta_{D'})\rho^A\eta^{A'} = \rho^A\eta^{A'}\cr
\rho^A\theta^{A'} &\mapsto (\sigma^D\rho_D)(\eta^{D'}\theta_{D'})\sigma^A\eta^{A'} = -\sigma^A\eta^{A'} & \sigma^A\eta^{A'} &\mapsto (\sigma^D\rho_D)(\eta^{D'}\theta_{D'})\rho^A\theta^{A'} = -\rho^A\theta^{A'}.\cr}}$$
Thus, constructing the null tetrad associated to the spin frames $\{\rho^A,\sigma^A\}$ and $\{\eta^{A'},\theta^{A'}\}$, and then the associated $\Psi$-ON basis (3.49), one deduces that with respect to the decomposition ${\bf R}^{2,2} = \langle U^a,V^a \rangle_{\bf R} \operp \langle X^a,Y^a \rangle_{\bf R}$, $\pm L = 1 \oplus (-1)$.

Case $R_2 \times R_2$, i.e., $\psi_{AB} = 2\mu_{(A}\nu_{B)}$ and $\chi_{A'B'} = 2\pi_{(A'}\tau_{B')}$, with $(\nu^D\mu_D)^2 = (\tau^{D'}\pi_{D'})^2 = 1$, whence $L \in {\bf O^\bfminus_\bfminus(2,2)}$ by (3.87). By (3.93), one can suppose $\{\mu^A,\nu^A\}$ and $\{\pi^{A'},\tau^{A'}\}$ are spin frames. Constructing the associated null tetrad and $\Psi$-ON basis (3.48--49), one computes that $L = 1 \oplus (-1)$ with respect to the decomposition ${\bf R}^{2,2} = \langle U^a,X^a \rangle_{\bf R} \operp \langle V^a,Y^a \rangle_{\bf R}$. Note that $-L_{ab}$ can be written as, say, $[\nu_A(-\mu_B) + (-\mu_A)\nu_B][\pi_{A'}\tau_{B'}+\tau_{A'}\pi_{B'}]$, with associated spin frames $\{\nu^A,-\mu^A\}$ and $\{\pi^{A'},\tau^{A'}\}$. For the $\Psi$-ON basis associated to these spin frames, $-L$ does take the form specified in the statement of the proposition.

Case $R_2 \times CC$, i.e., $\psi_{AB} = 2\mu_{(A}\nu_{B)}$ and $\chi_{A'B'} = \pm\lambda_{(A'}\bar\lambda_{B')}$, $\lambda_A = \eta_{A'} + i\theta_{B'}$, whence $\pm L^a{}_b = (\mu^A \nu_B + \nu^A\mu_B)(\eta^{A'}\eta_{B'} + \theta^{A'}\theta_{B'}) \in F_2{\bf O^\bfminus_\bfminus(2,2)}$ and are a complex-conjugate pair of anti-orthogonal complex structures. By (3.93), one can suppose $\{\mu^A,\nu^A\}$ and $\{\eta^{A'},\theta^{A'}\}$ are spin frames. The associated $\Psi$-ON basis $\{U^a,V^a,X^a,Y^a\}$ is such that $\{U^a,V^a,Y^a,X^a\}$ is a complex basis for $L^a{}_b$.

Case $CC \times R_2$, i.e., $\psi_{AB} = \pm\kappa_{(A}\bar\kappa_{B)} = \pm(\rho_A\rho_B + \sigma_A\sigma_B)$ and $\chi_{A'B'} = \pi_{A'}\tau_{B'} + \tau_{A'}\pi_{B'}$, whence $\pm L^a{}_b =\psi^A{}_B\chi^{A'}{}_{B'} \in F_2{\bf SO^\bfplus(2,2)}$ are a complex-conjugate pair of anti-orthogonal complex structures. By (3.87) and (3.93), one can suppose $\{\rho^A,\sigma^A\}$ and $\{\pi^{A'},\tau^{A'}\}$ are spin frames, and for the associated $\Psi$-ON basis $\{U^a,V^a,X^a,Y^a\}$, $\{U^a,X^a,Y^a,V^a\}$ is a complex basis of $L^a{}_b$.

Now suppose $L \in {\bf SNO(2)}$ satisfies $L_{ab} = L_{[ab]}$, whence in (3.83) $pg_{ab} + T_{ab} = 0$. Since $T_{ab}$ is traceless, transvecting by $g^{ab}$ yields $p=0$ and $T_{ab} = 0$, i.e., $\alpha\delta=0$ and $\psi_{AB}\chi_{A'B'} = 0$, which give two possibilities: $\alpha =0$ and $\chi_{A'B'}=0$, i.e., $L_{ab} = \delta\psi_{AB}\epsilon_{A'B'}$, and $\delta=0$ and $\psi_{AB} = 0$, i.e., $L_{ab} = \alpha\epsilon_{AB}\chi_{A'B'}$. 

For the first possibility, (3.87) implies $\psi_{AB} = \pm\kappa_{(A}\bar\kappa_{B)} = \pm(\rho_A\rho_B+\sigma_A\sigma_B)$ iff $\det(\alpha^A{}_B) = 1$ and $\psi_{AB} = \mu_A\nu_B + \nu_A\mu_B$ iff $\det(\alpha^A{}_B) = -1$. As $\det(\delta\epsilon^{A'}{}_{B'}) = \delta^2 > 0$, $\delta^{A'}{}_{B'} \in {\bf SL(2;R)}$, $\delta = \pm 1$, and noting (3.93), one obtains the two possibilities (a) and (b), respectively. The further statements in (a--b) are confirmed by simple computations. 

For the second possibility, (3.87) implies $\chi_{A'B'} = \pm\lambda_{(A'}\bar\lambda_{B')} = \pm(\eta_{A'}\eta_{B'}+\theta_{A'}\theta_{B'})$ iff $\det(\delta^{A'}{}_{B'}) = 1$ and $\chi_{A'B'} = \pi_{A'}\tau_{B'} + \tau_{A'}\pi_{B'}$ iff $\det(\delta^{A'}{}_{B'}) = -1$. As $\det(\alpha\epsilon^A{}_B) = \alpha^2 > 0$, $\alpha^A{}_B \in {\bf SL(2;R)}$, $\alpha = \pm 1$, and noting (3.93), one obtains the two possibilities (c) and (d) respectively.\bull
\vskip 24pt
As is well known, the Lie algebra {\bf so(p,q)} of {\bf O(p,q)} is the linear space $\End_-({\bf R}^{p,q})$ of skew-adjoint endomorphisms of ${\bf R}^{p,q}$; consequently, all these Lie algebras are linearly isomorphic to the space of rank two, skew tensors on ${\bf R}^{p+q}$ by virtue of (3.78.1) In particular, if $S^a{}_b \in \End_-({\bf R}^{2,2})$, then
$$S_{ab} = \phi_{AB}\epsilon_{A'B'} + \epsilon_{AB}\xi_{A'B'},\eqno(3.95)$$
for some symmetric spinors $\phi_{AB}$ and $\xi_{A'B'}$. Now ${\bf Spin^\bfplus(2,2)} \cong {\bf SL(2;R)} \times {\bf SL(2;R)}$, whence ${\bf spin(2,2)} \cong {\bf sl(2;R)} \oplus {\bf sl(2;R)}$. Since the homomorphism $\Upsilon$ of (3.3) is an isomorphism on a neighbourhood of the identity, it induces an isomorphism of Lie algebras, which one can recognize in (3.95) as follows. The Lie algebra of {\bf SL(2;R)} is the linear subspace of traceless elements in ${\bf R}(2)$. Noting that $\phi^A{}_A = \phi_{BA}\epsilon^{AB}$, then $\phi_{AB}$ is symmetric iff $\phi^A{}_B$ is traceless as an endomorphism of $S \cong {\bf R}^2$.

The induced action of ${\bf SO^\bfplus(2,2)}$ on the spaces $\Lambda^2_+({\bf R}^{2,2})$/$\Lambda^2_-({\bf R}^{2,2})$ of self-dual(SD)/anti-self-dual (ASD) two-multivectors of ${\bf R}^{2,2}$ yields a smooth, surjective homomorphism $\rho^+$, with kernel $\{\pm 1\}$, from ${\bf SO^\bfplus(2,2)}$ to ${\bf SO^\bfplus(1,2)} \times {\bf SO^\bfplus(1,2)}$ (where the induced scalar product on the space of multivectors is employed). This homomorphism induces an isomorphism: ${\bf so(2,2)} \cong {\bf so(1,2)} \oplus {\bf so(1,2)}$. Also, $\rho^+$ lifts to spinors as follows. If $L \in {\bf SO^\bfplus(2,2)}$, and $L\ \leftrightarrow\ \pm(\alpha^A{}_B,\delta^{A'}{}_{B'}) \in {\bf SL(2;R)} \times {\bf SL(2;R)} \cong {\bf Spin^\bfplus(2,2)}$, the action of $L^a{}_b$ on $S^{ab}$ corresponds to $\alpha^A{}_C\alpha^B{}_D\phi^{CD}$ and $\delta^{A'}{}_{C'}\delta^{B'}{}_{D'}\xi^{C'D'}$, i.e., to the induced actions of the two factors of ${\bf Spin^\bfplus(2,2)}$ on $S \odot S$ and $S' \odot S'$, each of which is, in turn, realizations of the vector representation of ${\bf Spin^\bfplus(1,2)} \cong {\bf SL(2;R)}$. Thus, at the level of Lie algebras,
$${\bf spin(1,2)} \oplus {\bf spin(1,2)} \cong {\bf spin(2,2)} \cong {\bf so(2,2)} \cong {\bf so(1,2)} \oplus {\bf so(1,2)},$$
which is also evident in (3.95).

Consider $\phi_{AB}$ in (3.95). The following considerations will also apply to $\xi_{A'B'}$. First suppose $\phi_{AB} = \kappa_{(A}\bar\kappa_{B)}$, with $\kappa_A = \rho_A + i\sigma_A$, $\sigma^D\rho_D > 0$. By (3.47), (3.80) and (3.86), $2\det(\phi^A{}_B) = \phi_{AB}\phi^{AB} = 2(\sigma^D\rho_D)^2 > 0$. Let $\lambda$ be a square root of $\det(\phi^A{}_B)$. Then, as $\phi_{AC}\phi^C{}_B = (1/2)\phi_{DE}\phi^{DE}\epsilon_{AB}$, $J^A{}_B := \lambda^{-1}\phi^A{}_B$ is a complex structure, whence
$$\exp(sJ^A{}_B) = -\epsilon^A{}_B\cos s + J^A{}_B\sin s,$$
and putting $s = \lambda t$ yields
$$\alpha^A{}_B := \exp(t\phi^A{}_B) = \phi^A{}_B\left({\sin(\lambda t)\over \lambda}\right) -\cos(\lambda t)\epsilon^A{}_B,\eqno(3.96)$$
which is of the form (3.82) and of case $CC$ in (3.88).

If instead $\phi_{AB} = \mu_A\nu_B + \nu_A \mu_B$, then $2\det(\phi^A{}_B) = \phi_{CD}\phi^{CD} = -2(\nu^D\mu_D)^2 \leq 0$. Let $\lambda$ now be a square root of $(\nu^D\mu_D)^2$, whence $P^A{}_B := \lambda^{-1}\phi^A{}_B$ is an involution, so
$$\exp(sP^A{}_B) = -\epsilon^A{}_B\cosh s + P^A{}_B\sinh s.$$
Putting $s=\lambda t$ yields
$$\alpha^A{}_B := \exp(t\phi^A{}_B) = \phi^A{}_B\left({\sinh(\lambda t)\over \lambda}\right) - \cosh(\lambda t)\epsilon^A{}_B.\eqno(3.97)$$
Now $\lambda = \pm(\nu^D\mu_D)$, so when $\lambda \not=0$, (3.97) is case $R_2$ of (3.88). In the limit $\lambda \to 0$, one obtains
$$\alpha^A{}_B = t\phi^A{}_B -\epsilon^A{}_B,\eqno(3.98)$$
which is case $R_1$ in (3.88) because $\nu^D\mu_D = 0$ when $\lambda = 0$.

Note that (3.96) describes one-parameter subgroups of {\bf SL(2;R)} homeomorphic to ${\bf S}^1$; in particular, when $\lambda=1$, such subgroups contain $\epsilon ^A{}_B$, i.e., $-1$ (e.g., for $t=\pi$). Moreover, taking the trace of (3.96) gives $\tr(\alpha^A{}_B) = 2\cos(\lambda t)$, whence such elements are so-called {\sl elliptic} or {\sl parabolic} elements of {\bf SL(2;R)}, see Carter et al. (1995), pp. 56--58. On the other hand, (3.97) and (3.98) describe one-parameter subgroups homeomorphic to {\bf R}, and taking the trace of (3.97) gives $\tr(\alpha^A{}_B) = 2\cosh(\lambda t) \geq 2$, i.e., the {\sl hyperbolic}, and certain parabolic, elements of {\bf SL(2;R)} in the image of the exponential mapping.

It is useful to have matrix descriptions of orthogonal transformations of ${\bf R}^{2,2}$ with respect to Witt bases, in particular with respect to the standard Witt basis associated to the standard basis via (3.50); this matrix representation may be called the {\sl standard matrix representation} of {\bf O(4;hb)}. By the previous construction, the standard Witt basis is $\{\check \ell^a,\check{\tilde m}^a,\check n^a,-\check m^a\}$, where $\{\check \ell^a,\check{\tilde m}^a,\check n^a,\check m^a\}$ is the null tetrad associated via (3.48) to the spin frames $\{\check o^A,\check\iota^A\}$ and $\{\check o^{A'},\check\iota^{A'}\}$. If
$$o^A := a\check o^A + b\check\iota^A \qquad \iota^A := c\check o^A + d\iota^A \qquad o^{A'} := \alpha \check o^{A'} + \beta\check\iota^{A'} \qquad \iota^{A'} := \gamma\check o^{A'} + \delta\check\iota^{A'},\eqno(3.99)$$
put
$$E_1^a := o^Ao^{A'} \qquad E_2^a := \iota^Ao^{A'} \qquad E_3^a := \iota^A\iota^{A'} \qquad E_4^a := -o^A\iota^{A'}.\eqno(3.100)$$
The transformation $L$ mapping the standard Witt basis to $\{E_1^a,\ldots,E_4^a\}$ has matrix representation
$$\underline L = \pmatrix{a\alpha&c\alpha&c\gamma&-a\gamma\cr b\alpha&d\alpha&d\gamma&-b\gamma\cr b\beta&d\beta&d\delta&-b\delta\cr -a\beta&-c\beta&-c\delta&a\delta\cr} =: \pmatrix{A&C\cr B&D\cr}\eqno(3.101)$$
with respect to the standard Witt basis; here $A$, $B$, $C$ and $D \in {\bf R}(2)$. When
$${\bf A}_1 := \pmatrix{a&c\cr b&d\cr} \hskip 1.25in {\bf A}_2 := \pmatrix{\alpha&\gamma\cr \beta&\delta\cr},\eqno(3.102)$$
are both elements of {\bf SL(2;R)}: $\{o^A,\iota^{A'}$ and $\{o^{A'},\iota^{A'}\}$ are spin frames; $\{o^Ao^{A'}, \iota^Ao^{A'}, \iota^A\iota^{A'},o^A\iota^{A'}\}$ is a null tetrad; (3.100) is a Witt basis; and $L$ is orthogonal. One easily confirms that the inverse of $\underline L$ is given by (1.2). Matrices of the form (3.101--102) therefore describe the identity component ${\bf SO^\bfplus(4;hb)}$ of {\bf O(4;hb)}. When ${\bf A}_1$ and ${\bf A}_2$ both belong to {\bf ASL(2;R)}: $\{\iota^A,o^A\}$ and $\{\iota^{A'},o^{A'}\}$ are spin frames; $\{o^Ao^{A'}, \iota^Ao^{A'}, \iota^A\iota^{A'},o^A\iota^{A'}\}$ is again a null tetrad; (3.100) is again a Witt basis; (1.2) is again the inverse of (3.101); and, by (3.56), (3.101) now describes the other connected component of {\bf SO(4;hb)}. Note that the two forms of (3.101) can be distinguished by the determinants of the eight $2 \times 2$-minors that have a common factor; in particular
$$\det(A) = \alpha^2, \qquad \det(B) = \beta^2, \qquad \det(C) = \gamma^2, \qquad \det(D) = \delta^2,\eqno(3.103)$$
for ${\bf SO^\bfplus(4;hb)}$ and 
$$\det(A) = -\alpha^2, \qquad \det(B) = -\beta^2, \qquad \det(C) = -\gamma^2, \qquad \det(D) = -\delta^2,\eqno(3.104)$$
for ${\bf SO(4;hb)} \setminus {\bf SO^\bfplus(4;hb)}$. Similarly, the $2 \times 2$-minors that have $\alpha$, $\beta$, $\gamma$ or $\delta$ as common factor, have determinants $a^2$, $-b^2$, $-c^2$, $d^2$ and $-a^2$, $b^2$, $c^2$, $-d^2$, in the two cases respectively.

Now consider mappings
$$\check o^A \mapsto \alpha\check o^{A'} + \beta\check\iota^{A'} \qquad \check\iota^A \mapsto \gamma\check o^{A'} + \delta\check\iota^{A'} \qquad \check o^{A'} \mapsto a\check o^A + b\check\iota^A \qquad \iota^{A'} \mapsto c\check o^A + d\check\iota^A.\eqno(3.105)$$
With the definitions (3.99) but now
$$E_1^a := o^Ao^{A'} \qquad E_2^a := o^A\iota^{A'} \qquad E_3^a := \iota^A\iota^{A'} \qquad E_4^a := - \iota^Ao^{A'},\eqno(3.106)$$
in place of (3.100), see (3.51b), the transformation $T$ mapping the standard Witt basis to $\{E_1^a,\ldots,E_4^a\}$ has matrix representation
$$\underline T = \pmatrix{a\alpha&a\gamma&c\gamma&-c\alpha\cr b\alpha&b\gamma&d\gamma&-d\alpha\cr b\beta&b\delta&d\delta&-d\beta\cr -a\beta&-a\delta&-c\delta&c\beta\cr} =: \pmatrix{A&C\cr B&D\cr},\eqno(3.107)$$
with respect to the standard Witt basis.

When ${\bf A}_1$ and ${\bf A}_2$ are both elements of {\bf SL(2;R)} or both elements of {\bf ASL(2;R)} then (3.106) is a Witt basis, $T$ is orthogonal, the inverse of (3.107) is given by (1.2), and
$$\det(A) = \det(B) = \det(C) = \det(D) = 0.\eqno(3.108)$$
When ${\bf A}_1$ and ${\bf A}_2$ are both elements of {\bf SL(2;R)}, one has, by (3.60), a matrix description of that component of {\bf O(4;hb)} corresponding to ${\bf O^\bfminus_\bfplus(2,2)}$ and the eight $2 \times 2$-minors that have a common factor have determinants $-a^2$, $b^2$, $-c^2$, $d^2$, $-\alpha^2$, $-\beta^2$, $\gamma^2$, and $\delta^2$. When ${\bf A}_1$ and ${\bf A}_2$ are both elements of {\bf ASL(2;R)}, one has, by (3.64), a matrix description of that component of {\bf O(4;hb)} corresponding to ${\bf O_\bfminus^\bfplus(2,2)}$ and the eight $2 \times 2$-minors that have a common factor have determinants $a^2$, $-b^2$, $c^2$, $-d^2$, $\alpha^2$, $\beta^2$, $-\gamma^2$, and $-\delta^2$.

Alternatively, one can use (3.101) with ${\bf A}_1$ and ${\bf A}_2$ in {\bf SL(2;R)} for a matrix description of the identity component and the coset structure of ${\bf O(2,2)} = {\bf O(4;hb)}$ to obtain matrix descriptions of the other components. For example, let
$$I_1 := \pmatrix{0&1&0&0\cr 1&0&0&0\cr 0&0&0&1\cr 0&0&1&0\cr} \quad I_2 := \pmatrix{-1&0&0&0\cr 0&1&0&0\cr 0&0&-1&0\cr 0&0&0&1\cr} \quad F := \pmatrix{0&0&-1&0\cr 0&1&0&0\cr -1&0&0&0\cr 0&0&0&1\cr} \quad J := \pmatrix{0&0&1&0\cr 0&1&0&0\cr 1&0&0&0\cr 0&0&0&1\cr}.\eqno(3.109)$$
Then
$$I_1^2 = I_2^2 = F^2 = J^2 = 1 \hskip 1.25in FJ = I_2 = JF.\eqno(3.110)$$
Interpreting the matrices in (3.109) as matrix representations of automorphisms with respect to the standard Witt basis, $I_1$ represents the automorphism $L_1$: $\check\ell^a \leftrightarrow \check{\tilde m}^a$, $\check n^a \leftrightarrow -\check m^a$, i.e., $e_1 \leftrightarrow e_2$ and $e_3 \leftrightarrow e_4$, i.e., $L_1$ has matrix representation $I_1$ with respect to the standard basis itself and is obviously an element of ${\bf O_\bfminus^\bfminus(2,2)}$. Multiplying $\underline L$ in (3.101) by $I_1$ on the left interchanges the first and second rows and the third and fourth rows of (3.102), which is equivalent to the interchange $a \leftrightarrow b$ and $c \leftrightarrow d$ in ${\bf A}_1$ and replacing $\beta$ by $-\beta$ and $\delta$ by $-\delta$ in ${\bf A}_2$, thus converting these matrices into elements of {\bf ASL(2;R)} and converting (3.103) into (3.104). Similarly, the automorphism $L_2$ represented by $I_2$ has $I_2$ as matrix representation with respect to the standard basis and is an element of ${\bf O_\bfminus^\bfminus(2,2)}$ too. Multiplying $\underline L$ in (3.101) by $I_2$ on the left is equivalent to replacing $a$ by $-a$ and $c$ by $-c$ in ${\bf A}_1$ and $\beta$ by $-\beta$ and $\delta$ by $-\delta$ in ${\bf A}_2$.

The automorphism $L_3$ represented by $F$ maps $e_1 \mapsto -e_1$ and fixes the other elements of the standard basis and thus is an element of ${\bf O^\bfminus_\bfplus(2,2)}$. Multiplying $\underline L$ in (3.101) by $F$ on the left is equivalent to replacing ${\bf A}_1$ by $\left({\beta \atop -\alpha}{\delta \atop -\gamma}\right)$ and ${\bf A}_2$ by $\left({-b \atop a}{-d \atop c}\right)$, both of which are elements of {\bf SL(2;R)} and with this substitution $F.{\underline L}$ is of the form (3.107). The automorphism $L_4$ represented by $J$ maps $e_3 \mapsto -e_3$ and fixes the other elements of the standard basis and thus is an element of ${\bf O^\bfplus_\bfminus(2,2)}$. Multiplying $\underline L$ in (3.101) by $J$ on the left is equivalent to replacing ${\bf A}_1$ by $\left({\beta \atop \alpha}{\delta \atop \gamma}\right)$ and ${\bf A}_2$ by $\left({b \atop a}{d \atop c}\right)$, each of which is an element of {\bf ASL(2;R)} and with this substitution $J.{\underline L}$ is of the form (3.107).

For spinor analysis on four-dimensional manifolds with a neutral metric, see Law (2009).
\vskip 24pt
\noindent {\section 4. Classification of Self-Adjoint Endomorphisms}
\vskip 12pt
The algebraic classification of the Ricci tensor $R_{ab}$ is just a special case of the algebraic classification of symmetric, rank-two covariant tensors $S_{ab} = S_{(ab)}$. By (3.78) $S_{ab}$ is symmetric iff the associated endomorphism $S^a{}_b := g^{ac}S_{cb}$ is self-adjoint. So, in general, the space of interest is $\End_+({\bf R}^{p,q})$. There are the following simple results for self-adjoint endomorphisms.
\vskip 24pt
\noindent {\bf 4.1 Proposition}\hfil\break
Let $L \in \End_+({\bf R}^{p,q})$.\hfil\break
i) any eigenvalue of $L_{\bf C}$ with a nonnull eigenvector in $\overline{\bf C}^{p,q}$ (in particular, real eigenvector nonnull in ${\bf R}^{p,q}$) is real. In particular, for a complex eigenvalue $\lambda$, $E_\lambda(L_{\bf C})$ is totally null in $\overline{\bf C}^{p,q}$. Consequently, for $L \in \End_+({\bf R}^{n,0})$, every eigenvalue is real.\hfil\break
ii) if $\lambda$ and $\mu$ are two distinct eigenvalues of $L$, then $E_\lambda(L_{\bf C})$ and $E_\mu(L_{\bf C})$ are orthogonal in ${\bf C}^{p,q}$.
\vskip 24pt
A classification of self-adjoint endomorphisms can be based on Jordan canonical forms (JCF). First some standard facts. 
\vskip 24pt
\noindent{\bf 4.2 Facts}\hfil\break
Let $J_m(\lambda)$ denote the $m \times m$ matrix with $\lambda$'s down the diagonal, 1's down the superdiagonal, and 0's elsewhere. I refer to this matrix as the {\sl Jordan block} of size $m$ for the eigenvalue $\lambda$. It will prove convenient to let $J_{-m}(\lambda)$ denote the $m \times m$ matrix with $\lambda$'s down the diagonal, -1's down the superdiagonal, and 0's elsewhere. Let
\vskip 6pt
$$K_t(a,b) := \pmatrix{\matrix{a&b\cr -b&a\cr}&\matrix{1&0\cr 0&1\cr}\cr \hbox{{\blowup 0}}&\matrix{a&b\cr -b&a\cr}&\ddots&\hbox{{\blowup 0}}\cr \noalign{\smallskip}\hbox{{\blowup 0}}&&\ddots&\matrix{1&0\cr 0&1\cr}\cr &&&\matrix{a&b\cr -b&a\cr}\cr} \in {\bf R}(2t).\eqno(4.2.1)$$
\vskip 6pt
If an endomorphism $S$ of ${\bf R}^n$ has a complex eigenvalue $\lambda$, suppose the Jordan blocks in the JCF of $S_{\bf C}$ determined by $\lambda$ are $J_{m_1}(\lambda),\ldots,J_{m_r}(\lambda)$. Then $\bar\lambda$ is also an eigenvalue of $S_{\bf C}$ and its Jordan blocks in the JCF of $S_{\bf C}$ are $J_{m_1}(\bar\lambda),\ldots,J_{m_r}(\bar\lambda)$. If $W$ is the summand in the JCF decomposition of ${\bf C}^n$ determining the Jordan blocks of $\lambda$, and $\{w_1,\ldots,w_m\}$ the basis of $W$ with respect to which $S_{\bf C}\vert_W$ has JCF $J_{m_1}(\lambda) \oplus \cdots \oplus J_{m_r}(\lambda)$ as matrix representation, then $\{\overline{w_1},\ldots,\overline{w_m}\}$ is the basis for $\overline W$ with respect to which $S_{\bf C}\vert_{\overline W}$ has JCF $J_{m_1}(\bar\lambda) \oplus \ldots \oplus J_{m_r}(\bar\lambda)$ as matrix representation. Define
$$u_j := {w_j + \overline{w_j} \over \sqrt 2} \hskip 1.25in v_j := {w_j - \overline{w_j} \over i\sqrt 2},\eqno((4.2.2)$$
whence
$$w_j = {u_j + iv_j \over \sqrt 2} \hskip 1.25in \overline{w_j} = {u_j - iv_j \over \sqrt2}.\eqno(4.2.3)$$
The set $B := \{u_1,v_1,\ldots,u_m,v_m\}$ is a subset of ${\bf R}^n$; putting $V := \langle u_1,v_1,\ldots,u_m,v_m \rangle_{\bf R}$, then $B$ is a basis for $V$ and ${\bf C}(V) = W \oplus \overline{W}$.

Let $w_{-1}$, $u_{-1}$ and $v_{-1}$ each denote {\bf 0}. The first Jordan block of $S_{\bf C}\vert_{W}$ arises from the elements $w_1,\ldots,w_{m_1}$, say, $m_1 \leq m$. Then, $S_{\bf C}(w_j) = \lambda w_j + w_{j-1}$, $1 \leq j \leq m_1$. Writing $\lambda = a + bi$, $a$, $b \in {\bf R}$, then
$$S(u_j) = au_j - bv_j + u_{j-1} \hskip 1.25in S(v_j) = bu_j + av_j + v_{j-1}.$$ 
Hence, with respect to $\{u_1,v_1,\ldots,u_{m_1},v_{m_1}\}$ $S\vert_{\langle u_1,v_1,\ldots,u_{m_1},v_{m_1} \rangle_{\bf R}}$ has matrix representation $K_{m_1}(a,b)$. It follows that with respect to the basis $B$, $S\vert_V$ has matrix representation that is block diagonal, the blocks being $K_{m_1}(a,b),\ldots,K_{m_r}(a,b)$, for some $m_1,\ldots,m_r$ satisfying $\sum_i\,m_i = m$. Repeating this construction for each pair of complex conjugate eigenvalues of $S$ allows one to construct a basis for ${\bf R}^n$ with respect to which $S$ takes the {\sl realized} JCF, i.e., is block diagonal, containing the Jordan blocks of real eigenvalues and the $K_{m_j}(a,b)$ for complex conjugate pairs of eigenvalues.
\vskip 24pt
The following result is stated in O'Neill (1983), pp. 261--262. A proof may be constructed by examining the construction of the Jordan canonical form and making more specialized choices to inductively construct bases that achieve the desired result. The proof is somewhat lengthy and I will not give it here.
\vskip 24pt
\noindent {\bf 4.3 Proposition}\hfil\break
An endomorphism $S \in \End_{\bf R}({\bf R}^{p+q})$ is in $\End_+({\bf R}^{p,q})$ iff ${\bf R}^{p,q}$ admits an orthogonal, $S$-invariant decomposition, viz.,
$${\bf R}^{p,q} = U_1 \operp \ldots \operp U_k\eqno(4.3.1)$$
with $S_\vdash(U_\ell) \leq U_\ell$, $\ell=1,\ldots,k$, such that either:\hfil\break
(1) $U_\ell$ admits a basis $\{v_1,\ldots,v_n\}$ ($n$ dependent upon $\ell$), satisfying, with $\epsilon = \pm 1$ (fixed for each summand),
$$s_{p,q}(v_i,v_j) = \cases{\epsilon,&for $i+j=n+1$,\cr
0,&otherwise,\cr}$$
with respect to which $S \vert_{U_\ell}$ has matrix representation a Jordan block $J_n(\lambda)$, where $\lambda$ is a real eigenvalue of $S$, or\hfil\break
(2) $U_\ell$ admits a basis $\{u_1,v_1,\ldots,u_m,v_m\}$ ($m$ dependent upon $\ell$), satisfying
$$s_{p,q}(u_i,u_j) = 1 = -s_{p,q}(v_i,v_j)\qquad\hbox{if}\qquad i+j=m+1$$ 
and all other scalar products zero, and with respect to which $S\vert_{U_\ell}$ has matrix representation $K_m(a,b)$, where $a+ib$, $a$, $b \in {\bf R}$, is a complex eigenvalue of $S$.

The basis of ${\bf R}^{p+q}$ provided by this decomposition yields the realized JCF for $S$.

In case (1), with $j := n+1 - i$, define
$$x_i := {v_i + \epsilon v_j \over \sqrt2} \hskip 1.25in y_i := {v_i - \epsilon v_j \over \sqrt2}.\eqno(4.3.2)$$
For $n=2k$, $\{v_1,\ldots,v_k,\epsilon v_{2k},\ldots,\epsilon v_{k+1}\}$ is a Witt basis and $\{x_1,\ldots,x_k,y_1,\ldots,y_k\}$ a $\Psi$-ON basis of signature $(n,n)$ for $U_\ell$. For $n = 2k+1$, $\{v_1,\ldots,v_k,\epsilon v_{2k+1},\ldots,\epsilon v_{k+2}\}$ is a Witt basis for the space\hfil\break
$\langle v_1,\ldots,v_k,v_{k+2},\ldots,v_{2k+1}\rangle_{\bf R}$ and $v_{k+1}$ is a unit vector; $\{x_1,\ldots,x_k,v_{k+1},y_1,\ldots,y_k\}$ is a $\Psi$-ON basis for $U_\ell$ of signature $(k+1,k)$ when $\epsilon = 1$ (in which case $y_{k+1} = {\bf 0}$) or of signature $(k,k+1)$ when $\epsilon = -1$ (in which case $x_{k+1} = {\bf 0}$).

In case (2), with $j := m+1 - i$, define
$$X^\pm_i := {u_i \pm u_j \over \sqrt2} \hskip 1.25in Y^\pm_i := {v_i \pm v_j \over \sqrt2}.\eqno(4.3.3)$$
When $m=2k$, $\{u_1,\ldots,u_k,v_1,\ldots,v_k,u_{2k},\ldots,u_{k+1},-v_{2k},\ldots,-v_{k+1}\}$ is a Witt basis and\hfil\break
$\{X^+_1,\ldots,X^+_k,Y^-_1,\ldots,Y^-_k,X^-_1,\ldots,X^-_k,Y^+_1,\ldots,Y^+_k\}$ is a $\Psi$-ON basis of signature $(2k,2k)$ for $U_\ell$. When $m=2k+1$, $\{u_1,\ldots,u_k,v_1,\ldots,v_k,u_{2k+1},\ldots,u_{k+2},-v_{2k+1},\ldots,-v_{k+2}\}$ is a Witt basis for the space\hfil\break
 $\langle u_1,\ldots,u_k,u_{k+2}\ldots,u_{2k+1},v_1,\ldots,v_k,v_{k+2},\ldots,v_{2k+1}\rangle_{\bf R}$ and $u_{k+1}$ and $v_{k+1}$ are unit vectors ;$X^-_{k+1} = Y^-_{k+1} = {\bf 0}$ and $\{X^+_1,\ldots,X^+_k,Y^-_1,\ldots,Y^-_k,u_{k+1},X^-_1,\ldots,X^-_k,Y^+_1,\ldots,Y^+_k,v_{k+1}\}$ is a $\Psi$-ON basis of signature $(m,m)$ for $U_\ell$.

Hence, odd-dimensional summands are of signature $(r,s)$ with $r-s = \pm 1$, which case only occurs for real eigenvalues, while even dimensional summands are neutral.
\vskip 24pt
\noindent {\bf 4.4 Corollary}\hfil\break
$S \in \End_+({\bf R}^{2,2})$ admits an orthogonal, $S$-invariant decomposition ${\bf R}^{2,2} = U_1 \operp \cdots \operp U_k$, $1 \leq k \leq 4$, such that:\hfil\break
i) for a one-dimensional summand $U = \langle v \rangle_{\bf R}$, $v$ is a unit vector and $S\vert_U = \lambda1$, for some eigenvalue $\lambda \in {\bf R}$;\hfil\break
ii) for a two-dimensional summand $U$, $U \cong {\bf R}^{1,1}$ and either\hfil\break
case (1) of (4.3) pertains, i.e., $U$ admits a basis $\{v_1,v_2\}$, with $s_{2,2}(v_1,v_2) = \epsilon = \pm 1$ the only nonzero scalar products, with respect to which $S\vert_U$ has matrix representation $J_2(\lambda)$, i.e., $U$ admits a Witt basis $\{v_1,\epsilon v_2\}$ with respect to which $S\vert_U$ has matrix representation
$$J_{\pm2}(\lambda) = \pmatrix{\lambda&\pm 1\cr0&\lambda\cr}\qquad\hbox{(which is self-adjoint according to (1.2))};$$
equivalently, with respect to the $\Psi$-ON basis $\{x_1,y_1\}$ of signature $(1,1)$ of (4.3.2), $S\vert_U$ has matrix representation
$$\pmatrix{\lambda+{\epsilon \over 2}&-\epsilon/2\cr \epsilon/2&\lambda - {\epsilon \over 2}\cr}\qquad\hbox{(which is self-adjoint according to (1.1))}$$
or\hfil\break
case 2 of (4.3) pertains, i.e., $U$ admits a basis $\{u_1,v_1\}$, which is $\Psi$-ON of signature $(1,1)$, with respect to which $S\vert_U$ has matrix representation
$$K_1(a,b) = \pmatrix{a&b\cr-b&a\cr};$$
equivalently, with $w := u_1 + iv_1$, $\lambda = a + ib$, $\{w,\bar w\}$ is a basis for $W := {\bf C}U$ with respect to which $S_{\bf C}\vert_W$ has matrix representation $J_1(\lambda) \oplus J_1(\bar\lambda)$, i.e.,
$$\pmatrix{\lambda&0\cr 0&\bar\lambda\cr},$$
$w$ and $\bar w$ are each null but not orthogonal to each other in $\overline{\bf C}^{2,2}$, but are each non-null and orthogonal to each other in ${\bf C}^{2,2}$, in accord with (4.1);\hfil\break
iii) for a three-dimensional summand $U$ case (1) of (4.3) pertains and, either\hfil\break
$U \cong {\bf R}^{1,2}$ and admits a basis $\{v_1,v_2,v_3\}$, whose nonzero scalar products are $s_{2,2}(v_1,v_3) = -1 = s_{2,2}(v_2,v_2)$ and with respect to which $S\vert_U$ has matrix representation $J_3(\lambda)$, $\lambda \in {\bf R}$, whence $\{x_1,v_2,y_1\}$ (with $\epsilon = -1$ in (4.3.2)) is a $\Psi$-ON basis of signature $(1,2)$, with respect to which $S\vert_U$ has matrix representation
$$\pmatrix{\lambda&{1 \over\sqrt2}&0\cr -{1 \over \sqrt2}&\lambda&{1 \over \sqrt2}\cr 0&{1 \over \sqrt2}&\lambda\cr}\qquad\hbox{(which is self-adjoint according to (1.1))},\eqno(4.4.1)$$
or\hfil\break
$U \cong {\bf R}^{2,1}$ and admits a basis $\{v_1,v_2,v_3\}$ whose nonzero scalar products are $s_{2,2}(v_1,v_3) = 1 = s_{2,2}(v_2,v_2)$ and with respect to which $S\vert_U$ has matrix representation $J_3(\lambda)$, $\lambda \in {\bf R}$, whence $\{x_1,v_2,y_1\}$ (with $\epsilon = 1$ in (4.3.2)) is a $\Psi$-ON basis of signature $(2,1)$, with respect to which $S\vert_U$ has matrix representation
$$\pmatrix{\lambda&{1 \over \sqrt2}&0\cr {1 \over \sqrt2}&\lambda&-{1 \over \sqrt2}\cr 0&{1 \over \sqrt2}&\lambda\cr}\qquad\hbox{(which is self-adjoint according to (1.1))};\eqno(4.4.2)$$
iv) for a four-dimensional summand, $U = {\bf R}^{2,2}$, and either\hfil\break
case (1) of (4.3) pertains, i.e., $U = {\bf R}^{2,2}$ admits a basis $\{v_1,v_2,v_3,v_4\}$ whose nonzero scalar products are $s_{2,2}(v_1,v_4) = s_{2,2}(v_2,v_3) = \epsilon = \pm1$ and with respect to which $S$ has matrix representation $J_4(\lambda)$, $\lambda \in {\bf R}$; $\{v_1,v_2,\epsilon v_4,\epsilon v_3\}$ is a Witt basis with respect to which $S$ has matrix representation
$$\pmatrix{\lambda&1&0&0\cr 0&\lambda&0&\epsilon\cr 0&0&\lambda&0\cr 0&0&1&\lambda\cr}\qquad\hbox{(which is self-adjoint according to (1.2))},\eqno(4.4.3)$$
and from (4.3.2) $\{x_1,x_2,y_1,y_2\}$ is a $\Psi$-ON basis of signature $(2,2)$ with respect to which $S$ has matrix representation
$$\pmatrix{\lambda&1/2&0&1/2\cr 1/2&\lambda + {\epsilon \over 2}&-1/2&-\epsilon/2\cr 0&1/2&\lambda&1/2\cr -1/2&\epsilon/2&1/2&\lambda - {\epsilon \over 2}\cr}\qquad\hbox{(which is self-adjoint according to (1.1))};\eqno(4.4.4)$$
or\hfil\break
case (2) of (4.3) pertains, i.e., $U = {\bf R}^{2,2}$ has a basis $\{u_1,v_1,u_2,v_2\}$ whose nonzero scalar products are $s_{2,2}(u_1,u_2) = 1 = -s_{2,2}(v_1,v_2)$ and with respect to which $S$ has matrix representation $K_2(a,b)$, $\lambda = a +ib$ and $\bar\lambda$ are eigenvalues with eigenvectors $w = u_1 + iv_1$ and $\bar w$ respectively; thus,$\{u_1,v_1,u_2,-v_2\}$ is a Witt basis with respect to which $S$ has matrix representation
$$\pmatrix{a&b&1&0\cr -b&a&0&-1\cr 0&0&a&-b\cr 0&0&b&a\cr}\qquad\hbox{(which is self-adjoint according to (1.2))},\eqno(4.4.5)$$
while from (4.3.3) $\{X^+_1,Y^-_1,X^-_1,Y^+_1\}$ is a $\Psi$-ON basis of signature $(2,2)$ with respect to which $S$ has matrix representation
$$\pmatrix{a+{1 \over 2}&0&-1/2&b\cr 0&a - {1 \over 2}&-b&1/2\cr 1/2&b&a - {1\over 2}&0\cr -b&-1/2&0&a+{1 \over 2}\cr}\qquad\hbox{(which is self-adjoint according to (1.1))}.\eqno(4.4.6)$$
\vskip 24pt
\noindent {\bf 4.5 Classification of Self-adjoint Endomorphisms of ${\bf R}^{2,2}$}\hfil\break
Putting together the possible summands in (4.4) results in the following classification. I use an explicit notation that describes the summands and the associated JCF on that summand. I also quote the Segre characteristic, which lists the sizes of Jordan blocks in the decomposition and groups together between round brackets those belonging to the same eigenvalue, for the sake of comparison with the classification in the case of Lorentzian signature, e.g., in PRII and Kramer et al. (1980), Table 5.1. Note one immediate difference. As noted in (4.3), even dimensional summands are neutral, whence no self-adjoint endomorphism of ${\bf R}^{1,3}$ can have a trivial decomposition in (4.3), whereas this is possible for self-adjoint endomorphisms of ${\bf R}^{2,2}$.\hfil\break
{\bf Type I}: (4.4)(iv), case (1): the JCF of $S$ is $J_4(\lambda)$, $\lambda \in {\bf R}$ with respect to a basis $\{v_1,v_2,v_3,v_4\}$ such that $\{v_1,v_2,\epsilon v_4,\epsilon v_3\}$ is a Witt basis, with respect to which the matrix representation is given by (4.4.3), while (4.4.4) gives the matrix representation with respect to the $\Psi$-ON basis $\{x_1,x_2,y_1,y_2\}$ of (4.3.2). Thus, $S$ has a single eigenvalue $\lambda$ with $m_\lambda(S) = 4$ and $M_\lambda(S) = 1$ (recall notation from (3.90)). I denote the type by $\bigl({\bf R}^4_{\rm hb},J_4(\lambda)\bigr)$(noting that the basis giving the JCF is a nonstandard ordering of a Witt basis). The Segre characteristic is [4] and this type has no analogue in the Lorentzian case.\hfil\break
{\bf Type II}: (4.4)(iv), case (2): $S$ has realized JCF $K_2(a,b)$ with respect to a basis $\{u_1,v_1,u_2,v_2\}$ such that $\{u_1,v_1,u_2,-v_2\}$ is a Witt basis, with respect to which $S$ has matrix representation given by (4.4.5), while (4.4.6) gives the matrix representation with respect to the $\Psi$-ON basis $\{X^+_1,Y^-_1,X^-_1,Y^+_1\}$ of (4.3.3). Thus, $S$ has a complex conjugate pair of eigenvalues $(\lambda,\bar\lambda)$, $\lambda = a + ib$, with $m_\lambda(S_{\bf C}) = 2$, $M_\lambda(S_{\bf C}) = 1$. With $w_1 := u_1 + iv_1$, $w_2 = u_2 + iv_2$, the basis $\{w_1,w_2,\bar w_1,\bar w_2\}$ of ${\bf C}^4$ gives the JCF $J_2(\lambda) \oplus J_2(\bar\lambda)$ of $S_{\bf C}$; each element of the basis is null in both ${\bf C}^{2,2} \cong {\bf C}^4_{\rm E}$ and $\overline{\bf C}^{2,2}$. Putting $U := \langle w_1,w_2 \rangle_{\bf C}$, then $U \cong \bar U \cong {\bf C}^2_{\rm E}$ in ${\bf C}^{2,2} \cong {\bf C}^4_{\rm E}$, while $U \oplus \bar U$ is a Witt decomposition in $\overline{\bf C}^{2,2}$. I denote the type by $\bigl({\bf R}^4_{\rm hb},K_2(a,b)\bigr)$ (the realized JCF being obtained with respect to a null tetrad), but one could also write $\bigl({\bf C}^2_{\rm E},J_2(\lambda)\bigr) \operp \bigl({\bf C}^2_{\rm E},J_2(\bar\lambda)\bigr)$ (for $S_{\bf C}$). The Segre characteristic is $[2\bar{2}]$ and this type has no analogue in the Lorentzian case.\hfil\break
{\bf Type IIIa}: the decomposition of (4.4) is ${\bf R}^{2,2} = U_1 \operp U_2$, $U_1 \cong {\bf R}^{1,0}$ and $U_2 \cong {\bf R}^{1,2}$. $U_1 = \langle v_1 \rangle$, $s_{2,2}(v_1,v_1) = 1$ (an instance of (4.4)(i)) and $U_2 = \langle v_2,v_3,v_4 \rangle_{\bf R}$, with nonzero scalar products $s_{2,2}(v_2,v_4) = -1 = s_{2,2}(v_3,v_3)$ (an instance of (4.4)(iii)). The basis $\{v_1,v_2,v_3,v_4\}$ gives the JCF $J_1(\lambda) \oplus J_3(\mu)$, $\lambda$, $\mu \in {\bf R}$. With $x := (v_2 - v_4)/\sqrt2$ and $y := (v_2 + v_4)/\sqrt2$, $\{v_1,x,v_3,y\}$ is a $\Psi$-ON basis with respect to which $S$ has matrix representation
$$\pmatrix{\lambda&0&0&0\cr 0&\mu&{1 \over\sqrt2}&0\cr 0&-{1 \over \sqrt2}&\mu&{1 \over \sqrt2}\cr 0&0&{1 \over \sqrt2}&\mu\cr}\qquad\hbox{(which is self-adjoint according to (1.1))}.\eqno(4.5.1)$$
Thus, $S$ has two real eigenvalues, with $m_\lambda(S) = 1 = M_\lambda(S)$ and $m_\mu(S) = 3$, $M_\mu(S) = 1$. I denote the type by $\bigl({\bf R}^{1,0},J_1(\lambda)\bigr) \operp \bigl({\bf R}^{1,2},J_3(\mu)\bigr)$.\hfil\break
{\bf Type IIIb}: the decomposition of (4.4) is ${\bf R}^{2,2} = U_1 \operp U_2$, $U_1 \cong {\bf R}^{2,1}$ and $U_2 \cong {\bf R}^{0,1}$. $U_2 = \langle v_4 \rangle$, $s_{2,2}(v_4,v_4) = -1$ (an instance of (4.4)(i)) and $U_1 = \langle v_1,v_2,v_3 \rangle_{\bf R}$, with nonzero scalar products $s_{2,2}(v_1,v_3) = 1 = s_{2,2}(v_2,v_2)$ (an instance of (4.4)(iii)). The basis $\{v_1,v_2,v_3,v_4\}$ gives the JCF $J_3(\lambda) \oplus J_1(\mu)$, $\lambda$, $\mu \in {\bf R}$. With $x := (v_1 + v_3)/\sqrt2$ and $y := (v_1 - v_3)/\sqrt2$, $\{x,v_2,y,v_4\}$ is a $\Psi$-ON basis with respect to which $S$ has matrix representation
$$\pmatrix{\lambda&{1 \over \sqrt2}&0&0\cr {1 \over \sqrt2}&\lambda&-{1 \over \sqrt2}&0\cr 0&{1 \over \sqrt2}&\lambda&0\cr \noalign{\smallskip}0&0&0&\mu}\qquad\hbox{(which is self-adjoint according to (1.1))}.\eqno(4.5.2)$$
Thus, $S$ has two real eigenvalues, with $m_\lambda(S) = 3$, $M_\lambda(S) = 1$ and $m_\mu(S) = 1 = M_\mu(S)$. I denote the type by $\bigl({\bf R}^{2,1},J_3(\lambda)\bigr) \operp \bigl({\bf R}^{0,1},J_1(\mu)\bigr)$. 

The Segre characteristic for Type III is [13] and this type is the analogue of the Lorentzian type $B$ in Kramer et al. (1980), Table 5.1, with a decomposition of the form ${\bf R}^{1,3} = {\bf R}^{1,2} \operp {\bf R}^{0,1}$. The coincidence $\lambda = \mu$ yields Segre characteristic [(13)].\hfil\break
{\bf Type IV}: the decomposition of (4.4) is ${\bf R}^{2,2} = U_1 \operp U_2 \operp U_3$, $U_1 \cong {\bf R}^{1,0}$, $U_2 \cong {\bf R}^2_{\rm hb}$, $U_3 \cong {\bf R}^{0,1}$; $U_1 = \langle v_1 \rangle_{\bf R}$ with $s_{2,2}(v_1,v_1) = 1$ (an instance of (4.4)(i)), $U_2 = \langle v_2,v_3 \rangle_{\bf R}$, with $\{v_2,\epsilon v_3\}$ ($\epsilon = \pm 1$) a Witt basis for $U_2$ (an instance of (4.4)(ii), case (1)), $U_3 = \langle v_4 \rangle_{\bf R}$, $s_{2,2}(v_4,v_4) = -1$ (an instance of (4.4)(i)). The basis $\{v_1,v_2,v_3,v_4\}$ gives the JCF $J_1(\lambda) \oplus J_2(\mu) \oplus J_1(\nu)$, with $\lambda$, $\mu$, $\nu \in {\bf R}$, of $S$. With
$$x := {v_2 + \epsilon v_3 \over \sqrt 2} \hskip 1.25in y := {v_2 - \epsilon v_3\over \sqrt2},\eqno(4.5.3)$$
as in (4.3.2), $\{v_1,x,y,v_4\}$ is a $\Psi$-ON basis with respect to which $S$ has matrix representation
$$\pmatrix{\lambda&0&0&0\cr 0&\mu+{\epsilon \over 2}&-\epsilon/2&0\cr 0&\epsilon/2&\mu - {\epsilon \over 2}&0\cr 0&0&0&\nu\cr}\qquad\hbox{(which is self-adoint according to (1.1))}.\eqno(4.5.4)$$
Thus, $S$ has three real eigenvalues, with $m_\lambda(S) = M_\lambda(S) = 1$, $m_\mu(S) = 2$, $M_\mu(S) = 1$, and $m_\nu(S) = M_\nu(S) = 1$. I denote the type by 
$$\bigl({\bf R}^{1,0},J_1(\lambda)\bigr) \operp \bigl({\bf R}^2_{\rm hb},J_{\pm2}(\mu)\bigr) \operp \bigl({\bf R}^{0,1},J_1(\nu)\bigr).$$
The Segre characteristic is [121] and this type is the analogue of the Lorentzian type $A3$ in Kramer et al. (1980), Table 5.1, with a decomposition of the form ${\bf R}^{1,3} = {\bf R}^{0,1} \operp {\bf R}^{1,1} \operp {\bf R}^{0,1}$. The possible coincidences of eigenvalues, with their Segre characteristics, are: $\lambda = \mu$ [(12)1]; $\mu = \nu$ [1(21)]; $\lambda = \nu$ $[(1\vert 2 \vert 1)]$ (where the vertical lines about 2 indicate it is not included in the coincidence of eigenvalues); and $\lambda = \mu = \nu$ [(121)].\hfil\break
{\bf Type V}: the decomposition of (4.4) is ${\bf R}^{2,2} = U_1 \operp U_2 \operp U_3$, $U_1 \cong {\bf R}^{1,0}$, $U_2 \cong {\bf R}^{1,1}$, $U_3 \cong {\bf R}^{0,1}$; $U_1 = \langle v_1 \rangle_{\bf R}$ with $s_{2,2}(v_1,v_1) = 1$ (an instance of (4.4)(i)), $U_2 = \langle v_2,v_3 \rangle_{\bf R}$, with $\{v_2,v_3\}$ a $\Psi$-ON basis for $U_2$ (an instance of (4.4)(ii), case (2)), $U_3 = \langle v_4 \rangle_{\bf R}$,$s_{2,2}(v_4,v_4) = -1$ (an instance of (4.4)(i)). The $\Psi$-ON basis $\{v_1,v_2,v_3,v_4\}$ gives the realized JCF $J_1(\lambda) \oplus K_1(a,b) \oplus J_1(\nu)$, with $\lambda$, $a$, $b$, $\nu \in {\bf R}$, of $S$:
$$\pmatrix{\lambda&0&0&0\cr 0&a&b&0\cr 0&-b&a&0\cr 0&0&0&\nu\cr}\qquad\hbox{(which is self-adoint according to (1.1))}.\eqno(4.5.5)$$
Hence, $S$ has two real eigenvalues with $m_\lambda(S) = M_\lambda(S) = 1 = m_\nu(S) = M_\nu(S)$ and a pair of complex conjugate eigenvalues $\mu$ and $\bar\mu$, $\mu := a + ib$, with $m_\mu(S_{\bf C}) = M_\mu(S_{\bf C}) = 1$. I denote the type by
$$\bigl({\bf R}^{1,0},J_1(\lambda)\bigr) \operp \bigl({\bf R}^{1,1},K_1(a,b)\bigr) \operp \bigl({\bf R}^{0,1},J_1(\nu)\bigr).$$
The Segre characteristic is $[11 \bar1 1]$ and this type is the analogue of the Lorentzian type $A2$ in Kramer et al. (1980), Table 5.1, with a decomposition of the form ${\bf R}^{1,3} = {\bf R}^{0,1} \operp {\bf R}^{1,1} \operp {\bf R}^{0,1}$. The only coincidence of eigenvalues that remains within type is: $\lambda = \nu$ $[(1\vert 1 \bar1 \vert 1)]$.

There are three possible decompositions into two summands each of dimension two. As each summand is even dimensional and therefore neutral, there is no Lorentzian analogue of these types.\hfil\break
{\bf Type VI}: the decomposition of (4.4) is ${\bf R}^{2,2} = U_1 \operp U_2$, $U_1 = \langle v_1,v_2 \rangle _{\bf R} \cong {\bf R}^2_{\rm hb}$ with $\{v_1,\epsilon v_2\}$ a Witt basis ($\epsilon = \pm 1$) (an instance of (4.4)(ii), case (1)) and $U_2 = \langle v_3,v_4\rangle_{\bf R} \cong {\bf R}^{1,1}$ with $\{v_3,v_4\}$ a $\Psi$-ON basis. The basis $\{v_1,v_2,v_3,v_4\}$ gives the realized JCF $J_2(\lambda) \oplus K_1(a,b)$ of $S$. With
$$x := {v_1 + \epsilon v_2 \over \sqrt 2} \hskip 1.25in y := {v_1 - \epsilon v_2\over \sqrt2},\eqno(4.5.6)$$
$\{x,v_3,v_4,y\}$ is a $\Psi$-ON basis with respect to which $S$ has matrix representation
$$\pmatrix{\lambda+{\epsilon \over 2}&0&0&-\epsilon/2\cr 0&a&b&0\cr 0&-b&a&0\cr \epsilon/2&0&0&\lambda - {\epsilon \over 2}\cr}\qquad\hbox{(which is self-adoint according to (1.1))}.\eqno(4.5.7)$$
Hence, $S$ has a real eigenvalue with $m_\lambda(S) = 2$, $M_\lambda(S) = 1$ and a pair of complex conjugate eigenvalues $\mu := a + ib$ and $\bar\mu$ with $m_\mu(S_{\bf C}) = 1 = M_\mu(S_{\bf C})$. I denote the type by
$$\bigl({\bf R}^2_{\rm hb},J_{\pm2}(\lambda)\bigr) \operp \bigl({\bf R}^{1,1},K_1(a,b)\bigr).$$
The Segre characteristic is $[21\bar1]$.\hfil\break
{\bf Type VII}: the decomposition of (4.4) is ${\bf R}^{2,2} = U_1 \operp U_2$, $U_1 = \langle v_1,v_2 \rangle _{\bf R} \cong {\bf R}^2_{\rm hb}$ with $\{v_1,\epsilon v_2\}$ a Witt basis ($\epsilon = \pm 1$) and $U_2 = \langle v_3,v_4 \rangle _{\bf R} \cong {\bf R}^2_{\rm hb}$ with $\{v_3,\omega v_4\}$ a Witt basis ($\omega = \pm 1$) (each an instance of (4.4)(ii), case (1)). The basis $\{v_1,v_2,v_3,v_4\}$ gives the JCF $J_2(\lambda) \oplus J_2(\mu)$ of $S$. The basis $\{v_1,v_3,\epsilon v_2,\omega v_4\}$ is a Witt basis, with respect to which $S$ has matrix representation
$$\pmatrix{\lambda&0&\epsilon&0\cr 0&\mu&0&\omega\cr 0&0&\lambda&0\cr 0&0&0&\mu\cr}\qquad\hbox{(which is self-adoint according to (1.2))}.\eqno(4.5.8)$$
With
$$x_1 := {v_1 + \epsilon v_2 \over \sqrt 2} \hskip .5in y_1 := {v_1 - \epsilon v_2 \over \sqrt2} \hskip .5in x_2 := {v_3 + \omega v_4 \over \sqrt 2} \hskip .5in y_2 := {v_3 - \omega v_4 \over \sqrt2},\eqno(4.5.9)$$
$\{x_1,x_2,y_1,y_2\}$ is a $\Psi$-ON basis with respect to which $S$ has matrix representation
$$\pmatrix{\lambda+{\epsilon \over 2}&0&-\epsilon/2&0\cr 0&\mu+{\omega \over 2}&0&-\omega/2\cr \epsilon/2&0&\lambda-{\epsilon \over 2}&0\cr 0&\omega/2&0&\mu-{\omega \over 2}\cr}\qquad\hbox{(which is self-adoint according to (1.1))}.\eqno(4.5.10)$$
Hence, $S$ has two real eigenvalues satisfying $m_\lambda(S) = 2 = m_\mu(S)$ and $M_\lambda(S) = 1 = M_\mu(S)$. I denote the type by
$$\bigl({\bf R}^2_{\rm hb},J_{\pm2}(\lambda)\bigr) \operp \bigl({\bf R}^2_{\rm hb},J_{\pm2}(\mu)\bigr).$$
The Segre characteristic is [22]. Coincidence of eigenvalues gives [(22)].\hfil\break
{\bf Type VIII}: the decomposition of (4.4) is ${\bf R}^{2,2} = U_1 \operp U_2$, $U_1 = \langle v_1,v_2 \rangle _{\bf R} \cong {\bf R}^{1,1}$, with $\{v_1,v_2\}$ a $\Psi$-ON basis, and $U_2 = \langle v_3,v_4 \rangle _{\bf R} \cong {\bf R}^{1,1}$, with $\{v_3,v_4\}$ a $\Psi$-ON basis (each an instance of (4.4)(ii), case (2)). The basis $\{v_1,v_2,v_3,v_4\}$ gives the realized JCF $K_1(a,b) \oplus K_1(c,d)$ of $S$. The basis $\{v_1,v_3,v_2,v_4\}$ is $\Psi$-ON and $S$ has matrix representation 
$$\pmatrix{a&0&b&0\cr 0&c&0&d\cr -b&0&a&0\cr 0&-d&0&c\cr}\qquad\hbox{(which is self-adoint according to (1.1))}.\eqno(4.5.11)$$
with respect to it. Hence, $S$ has two pairs of complex conjugate eigenvalues $\lambda := a+ ib$, $\bar\lambda$ and $\mu := c + id$, $\bar\mu$, satisfying $m_\lambda(S_{\bf C}) = 1 = m_\mu(S_{\bf C})$ and $M_\lambda(S_{\bf C}) = 1 = M_\mu(S_{\bf C})$. I denote this type by
$$\bigl({\bf R}^{1,1},K_1(a,b)\bigr) \operp \bigl({\bf R}^{1,1},K_1(c,d)\bigr).$$
The Segre characteristic is $[1\bar11\bar1]$ and the coincidence $\lambda=\mu$ gives $[(1\bar11\bar1)]$.\hfil\break
{\bf Type IX}: the decomposition of (4.4) is ${\bf R}^{2,2} = U_1 \operp U_2 \operp U_3 \operp U_4$, $U_i = \langle v_i \rangle_{\bf R}$ (each an instance of (4.4)(i)), with $\{v_1,v_2,v_3,v_4\}$ a $\Psi$-ON basis giving the JCF $J_1(\lambda) \oplus J_1(\mu) \oplus J_1(\nu) \oplus J_1(\sigma)$ of $S$. Thus, $S$ has four real eigenvalues, each of algebraic and geometric multiplicity one. I denote the type by
$$\bigl({\bf R}^{1,0},J_1(\lambda)\bigr) \operp \bigl({\bf R}^{1,0},J_1(\mu)\bigr) \operp \bigl({\bf R}^{0,1},J_1(\nu)\bigr) \operp \bigl({\bf R}^{0,1},J_1(\sigma)\bigr).$$
The Segre characteristic is [1111] and this type is the analogue of the Lorentzian type $A1$ of in Kramer et al. (1980), Table 5.1, with a decomposition of the form ${\bf R}^{1,3} = {\bf R}^{1,0} \operp {\bf R}^{0,1} \operp {\bf R}^{0,1} \operp {\bf R}^{0,1}$. Coincidences of eigenvalues gives: $\lambda = \mu$ [(11)11]; $\lambda = \nu$ $[(1\vert 1\vert 1)1]$; $\lambda = \sigma$ $[(1\vert 11\vert 1)]$; $\mu = \nu$ [1(11)1]; $\mu = \sigma$ $[1(1\vert 1 \vert 1)]$; $\nu = \sigma$ [11(11)]; $\lambda = \mu$ and $\nu =\sigma$ [(11)(11)]; $\lambda = \nu$ and $\mu = \sigma$ $[(1 \vert (1 \vert 1)\vert1)]$; $\lambda = \sigma$ and $\mu = \nu$ $[(1\vert (11) \vert 1)]$; $\lambda = \mu =\nu$ [(111)1]; $\lambda = \mu =\sigma$ $[(11\vert 1 \vert 1)]$; $\lambda = \nu =\sigma$ $[(1\vert 1 \vert 11)]$; $\mu = \nu =\sigma$ [1(111)]; $\lambda = \mu = \nu =\sigma$ [(1111)].

In addition to the coincidences within type, there are degenerations from one type to another when complex conjugate eigenvalues coincide, i.e., become real: when $\lambda = \bar\lambda$ ($b=0$) in type II, $[2\bar2]$ degenerates to [(22)] in type VII; when $\mu = \bar\mu$ in type V, $[11\bar11]$ degenerates to [1(11)1] in type IX; when $\mu = \bar\mu$ in type VI, $[21\bar1]$ degenerates to $[(1\vert 2 \vert 1)]$ in type IV (with a re-arrangement of basis); when $\lambda = \bar\lambda$ or $\mu = \bar\mu$ in type VIII, $[1\bar11\bar1]$ degenerates to $[(11)1\bar1]$ in type V, (with a re-arrangement of basis); a further coincidence of the remaining complex conjugate pair reduces the original form to $[(11)(11)]$ in type IX (without a re-arrangement of basis).\bull
\vskip 24pt
For each type there is a specific matrix $\underline M$, say, representing an endomorphism of that type with respect to some $\Psi$-ON (or Witt) basis, e.g., (4.5.2) for type IIIb, and these forms are distinct between types. This fact is obvious from the explicit forms given for $\underline M$ in (4.5) and also from the fact that the different forms are characterized by different JCFs. Now observe that there is a smooth left action of {\bf O(p,q)} on $\End_+({\bf R}^{p,q})$ given by
$$T \mapsto LTL^{-1} = LT\,{^*\! L}.\eqno(4.6)$$

The following result is stated for the purposes of this paper but could be given for $\End_+({\bf R}^{p,q})$. Note that $\End_+({\bf R}^{p,q})$ is an {\bf R}-linear space of dimension $n(n+1)/2$.
\vskip 24pt
\noindent {\bf 4.7 Proposition}\hfil\break
The action (4.6) of {\bf O(2,2)} on $\End_+({\bf R}^{2,2})$ is not transitive; each orbit lies within one of the 10 (distinguishing IIIa and IIIb) types given in (4.5). Fix a type, say IIIb. Then its $\underline M$ has $2$ free parameters, $\lambda$ and $\mu$, so write $\underline M(\lambda,\mu)$.Suppose the self-adjoint endomorphisms are represented by matrix representations with respect to the standard basis.  Write $T(\lambda,\mu)$ for that whose standard matrix representation is $\underline M(\lambda,\mu)$. Then the orbits of the $T(\lambda,\mu)$ for the different values of $(\lambda,\mu)$ are distinct and the disjoint union of these orbits comprises the subset of type IIIb self-adjoint endomorphisms. Let $K(\lambda,\mu)$ denote the subgroup of {\bf O(2,2)} whose elements commute with $T(\lambda,\mu)$. Then the orbit containing $T(\lambda,\mu)$ is diffeomorphic to ${\bf O(2,2)}/K(\lambda,\mu)$ and the number of degrees of freedom in type IIIb is $\bigl(6 - \dim(K(\lambda,\mu))\bigr) + 2$. For a generic orbit, $K(\lambda,\mu)$ is independent of $(\lambda,\mu)$. An analogous statement applies to each of the 10 types. One can equally well restate this result in terms of the matrix forms each type takes with respect to Witt bases rather than $\Psi$-ON bases.

Proof. Since the elements of an orbit are conjugates, they have the same eigenvalues and algebraic and geometric multiplicities. It follows that there is more than one orbit, that the orbits lie within types, and that $\underline M(\lambda',\mu') = {\underline L}\underline M(\lambda,\mu){\underline L}^{-1}$ iff $(\lambda',\mu') = (\lambda,\mu)$ since $(\lambda,\mu)$ determines the eigenvalue structure of $\underline M(\lambda,\mu)$. If $S$ and $T$ are two self-adjoint endomorphisms of the same type and with the same parameter values $(\lambda_0,\mu_0)$, say, in $\underline M$, then there are $\Psi$-ON bases $\{u_1,\ldots,u_4\}$ and $\{v_1,\ldots,v_4\}$ such that $S$ and $T$ each have matrix representation $\underline M(\lambda_0,\mu_0)$ with respect to $\{u_1,\ldots,u_4\}$ and $\{v_1,\ldots,v_4\}$, respectively. Let $L \in {\bf O(2,2)}$ be the element such that $L(u_i) = v_i$, $i=1,\ldots,4$. Then, one easily checks that $T = LSL^{-1}$. Thus, the orbits are as claimed.

Since {\bf O(2,2)} is a Lie group, each orbit is diffeomorphic to the quotient of {\bf O(2,2)} by the isotropy subgroup of the action on that orbit (e.g., see Warner 1983, 3.62); in particular has a manifold structure. Generically, the isotropy subgroup $K(\lambda,\mu)$ is independent of $(\lambda,\mu)$ (see details in (4.9) below) but will be larger when there are coincidences between the eigenvalues.\bull
\vskip 24pt
\noindent {\bf 4.8 Observation}\hfil\break
Consider the mapping $\chi:{\bf R}(n) \to {\bf R}^{n+1}$ that maps a matrix $M$ to the vector of coefficients of the characteristic polynomial $\det(M - \lambda1)$. Consider the left action of {\bf GL(n;R)} on ${\bf R}(n)$ generalizing (4.6): $M \mapsto LML^{-1}$, $M \in {\bf R}(n)$ and $L \in {\bf GL(n;R)}$. Then, $\chi$ is constant on the orbits of this action but takes different values on different orbits. Hence, $\chi$ provides a model of the projection onto the space of orbits; moreover, as $\chi$ is a continuous mapping to ${\bf R}^{n+1}$, one concludes the orbits are closed subspaces of ${\bf R}(n)$.

When the action is restricted to {\bf O(p,q)} in {\bf GL(p+q)} and the closed subspace $\End_+({\bf R}^{p,q})$ of $\End({\bf R}^{p+q})$ is viewed, via matrix representation, as a closed subspace of ${\bf R}(n)$, $\chi: \End_+({\bf R}^{p,q}) \to {\bf R}^{p+q+1}$ is still well defined and models the projection onto the space of orbits of the action (4.6) and one still concludes that the orbits of the action (4.6) on $\End_+({\bf R}^{p,q})$ are closed subspaces. Thus, if $m$ is the number of free parameters in the matrix form $\underline M$ of a given type of self-adjoint endomorphism, the subspace of $\End_+({\bf R}^{p,q})$ consisting of those elements of that type consists of a collection of closed subspaces (namely orbits) parametrized by ${\bf R}^m$. For example, for type IIIb for $(p,q) = (2,2)$, the closed orbits are parametrized by ${\bf R}^2$. Recall, however, that for a type for which coincidences of eigenvalues can occur, not all orbits have the same dimension.
\vskip 24pt
In practice, for a given type, I will choose the matrix form $\underline M$ of that type with respect to either $\Psi$-ON or Witt bases, whichever is more convenient, and then determine the subgroup of (the relevant matrix representation of ) {\bf O(2,2)} that commutes with $\underline M$. This computation will reveal the degrees of freedom in the type and something of the topology of the generic orbit. I will also consider the special cases of coincidences of eigenvalues. Note that if $T$ is an orthogonal automorphism commuting with $S$, and $U$ is $S$ invariant, then $T_\vdash(U)$ is also $S$ invariant. Moreover, if $u \in E_\lambda(S)$, then $T(u) \in E_\lambda(S)$ and is of the same character as $u$. These simple facts, bearing (4.3.1) in mind, yield expectations as to the form of the relevant isotropy subgroups. I won't elaborate these expectations; they will be evident from the explicit form of the isotropy subgroups computed below.
\vskip 24pt
\noindent {\bf 4.9 Observation}\hfil\break
Since the aim is to relate the classification of Ricci tensors to that of Ricci spinors, and the latter is the spinor version of the traceless Ricci tensor, for any self-adjoint endomorphism $R^a{}_b$, define $\Phi^a{}_b := R^a{}_b - (\tr(R)/4)1$. Of course, $\tr(R) = R^a{}_a = R_{ab}g^{ab} =: S$ and $\Phi^a{}_b$ is traceless. The classification of (4.5) induces a classification of the codimension-one linear subspace $\End^0_+({\bf R}^{2,2})$, of traceless self-adjoint endomorphisms, of $\End_+({\bf R}^{2,2})$. Moreover, the action (4.6) restricts to an action on $\End^0_+({\bf R}^{2,2})$ and the obvious analogue of (4.7) is valid. I therefore now restrict attention to $\End^0_+({\bf R}^{2,2})$. Note that $R(v) = \lambda v$ iff $\Phi(v) = \mu v$, where $\mu := \lambda - S/4$ and that
$$\det(R - \lambda1) = \det(\Phi + (S/4)1 - \lambda 1) = \det(\Phi - \mu 1),$$
i.e., the characteristic equation for $R$ in $\lambda$ equals the characteristic equation for $\Phi$ in $\mu$, whence the eigenvalue $\mu$ of $\Phi$ has the same algebraic and geometric multiplicities as the corresponding eigenvalue $\lambda$ of $R$.
\vskip 24pt
\noindent {\bf 4.10 Computation of Degrees of Freedom}\hfil\break
{\bf Type I}: for each self-adjoint endomorphism $S$ of type I, there is a Witt basis with respect to which $S$ has matrix representation (4.4.3); $\tr(S) = 0$ is equivalent to $\lambda = 0$, so for type I, let $\underline M$ denote the matrix of (4.4.3) with $\lambda = 0$, whence $\underline M$ itself has no free parameters. One computes that the elements of $\underline {\bf O(4;hb)}$ (i.e., the matrix representation of {\bf O(4;hb)} with respect to any Witt basis) that commute with $\underline M$ in this case are $\pm 1$ only. Thus, type I consists of a single orbit under the action of (4.6) on $\End_+^0({\bf R}^{2,2})$, and that orbit is diffeomorphic to ${\bf O(4;hb)}/\{\pm 1\}$. Since $\pm 1$ both lie in the identity component, the orbit is six-dimensional and has four connected components.\hfil\break
{\bf Type II}: for each $S$ of Type II, there is a Witt basis with respect to which $S$ has matrix representation (4.4.5); $\tr(S) = 0$ is equivalent to $a=0$ (and $b \not= 0$ is required to remain within type II), so let $\underline M(b)$ denote the matrix in (4.4.5) with $a=0$, which has one free parameter itself. One computes that the elements of $\underline {\bf O(4;hb)}$ that commute with $\underline M(b)$ are $\pm 1$ only. Thus, type II consists of two ($b > 0$ and $b < 0$) one-parameter families of orbits under the action of (4.6) on $\End_+^0({\bf R}^{2,2})$, and these orbits are diffeomorphic to ${\bf O(4;hb)}/\{\pm 1\}$, whence each orbit is again six-dimensional and has four connected components.\hfil\break
{\bf Type IIIa}: for each $S$ of Type IIIa, there is a $\Psi$-ON basis with respect to which $S$ has matrix representation (4.5.1); $\tr(S) = 0$ is equivalent to $\lambda + 3\mu = 0$. Let $\underline M(\mu)$ denote the matrix in (4.5.1) with $\lambda = -3\mu$. The elements of $\underline {\bf O(2,2)}$ that commute with $\underline M(\mu)$, $\mu \not= 0$, are $\pm 1$ and $\pm Z_1$, where
$$Z_1 := \pmatrix{1&0&0&0\cr0&-1&0&0\cr 0&0&-1&0\cr 0&0&0&-1\cr}.\eqno(4.10.1)$$
Relative to a $\Psi$-ON basis, $\pm Z_1$ represents an element of ${\bf O^\bfminus_\bfplus(2,2)}$. Thus, the isotropy subgroup $K(\mu)$, $\mu \not= 0$, is discrete so there are two ($\mu\not=0$) one-parameter families of orbits, each of which is six dimensional, but has only two connected components.

The coincidence $\lambda=\mu$ yields $\mu=0$ in the traceless case. If
$$T_1(h) := \pmatrix{1&h&0&-h\cr -h&{2-h^2 \over 2}&0&h^2/2\cr 0&0&1&0\cr -h&-h^2/2&0&{2+h^2 \over 2}\cr},\eqno(4.10.2)$$
for $h \in {\bf R}$, then $T_1(0) = 1$, and $T_1(h)T_1(k) = T_1(h+k)$, i.e., $\{\,T_1(h):h \in {\bf R}\,\}$ is a subgroup of $\underline {\bf O(2,2)}$; in fact, by (1.5), of $\underline {\bf SO^\bfplus(2,2)}$. The isotropy subgroup $K(0)$ is generated by this subgroup together with the elements $-1$ and $Z_1$. Hence, the orbit for $\lambda = \mu=0$ also has only two connected components but is only five dimensional.\hfil\break
{\bf Type IIIb}: for each $S$ of Type IIIb, there is a $\Psi$-ON basis with respect to which $S$ has matrix representation (4.5.2); $\tr(S) = 0$ is equivalent to $3\lambda + \mu = 0$. Let $\underline M(\lambda)$ denote the matrix in (4.5.2) with $\mu = -3\lambda$. The elements of $\underline {\bf O(2,2)}$ that commute with $\underline M(\lambda)$, $\lambda \not= 0$, are $\pm 1$ and $\pm Z_2$, where
$$Z_2 := \pmatrix{1&0&0&0\cr 0&1&0&0\cr 0&0&1&0\cr 0&0&0&-1\cr}.\eqno(4.10.3)$$
With respect to a $\Psi$-ON basis $\pm Z_2$ represents an element of ${\bf O^\bfplus_\bfminus(2,2)}$. Hence, the isotropy subgroup $K(\lambda)$ ($\lambda \not= 0$) is discrete and there are two ($\lambda \not= 0$) one-parameter families of orbits, each of which is six dimensional and has two connected components.

The coincidence $\lambda = \mu$ yields $\lambda = 0$ in the traceless case. If 
$$T_2(g) := \pmatrix{{2+g^2 \over 2}&0&-g^2/2&g\cr 0&1&0&0\cr g^2/2&0&{2-g^2 \over 2}&g\cr g&0&-g&1\cr},\eqno(4.10.4)$$
for $g \in {\bf R}$, then $T_2(0) = 1$ and $T_2(g)T_2(h) = T_2(g+h)$, i.e., $\{\,T_2(g):g \in {\bf R}\,\}$ is a subgroup of $\underline {\bf O(2,2)}$, in fact, by (1.5), of $\underline {\bf O^\bfplus_\bfplus(2,2)}$. The isotropy subgroup $K(0)$ is generated by this subgroup together with the elements $-1$ and $Z_2$. Hence the orbit for $\lambda=0$ also has only two connected components but is only of dimension five. The analogy between IIIa and IIIb is clear.\hfil\break
{\bf Type IV}: for each $S$ of type IV, there is a $\Psi$-ON basis with respect to which $S$ has matrix representation (4.5.4); $\tr(S) = 0$ is equivalent to $-2\mu = \lambda + \nu$. Let $\underline M(\lambda,\nu)$ denote the matrix in (4.5.4) with $-(\lambda + \nu)/2$ substituted for $\mu$. The elements of $\underline {\bf O(2,2)}$ that commute with generic $\underline M(\lambda,\nu)$ (i.e., excluding the cases of coincidence treated below) are $\pm 1$, $\pm Z_1$, $\pm Z_2$ and $\pm Z_3$, where
$$Z_3 := \pmatrix{1&0&0&0\cr 0&-1&0&0\cr 0&0&-1&0\cr 0&0&0&1\cr}.\eqno(4.10.5)$$
With respect to a $\Psi$-ON basis, $Z_3$ represents an element of ${\bf O^\bfminus_\bfminus(2,2)}$. Thus, the isotropy subgroup is discrete, but intersects each component of {\bf O(2,2)}. Hence, when $\lambda + 3\nu \not= 0$, $3\lambda + \nu \not =0$, and $\lambda \not= \nu$, (the generic condition with no coincidences amongst eigenvalues), there are six two-parameter families of orbits, each of which is six dimensional and has a single connected component.

The coincidence $\mu=\nu$ is equivalent to $\lambda + 3\nu = 0$ under the trace-free condition. With $3\lambda + \nu \not= 0$ and $\lambda \not= \nu$, the isotropy group is now generated by elements with matrix representations $Z_1$, $Z_2$, $Z_3$ and $H := \{\,T_3(k): k \in {\bf R}\,\}$, where
$$T_3(k) = \pmatrix{1&0&0&0\cr 0&{2+k^2 \over 2}&-k^2/2&k\cr 0&k^2/2&{2-k^2 \over 2}&k\cr 0&k&-k&1\cr}.\eqno(4.10.6)$$
$T_3(0) = 1$, $T_3(k)T_3(h) = T_3(k+h)$, and $H$ is a subgroup of $\underline {\bf SO^\bfplus(2,2)}$. Hence, there are two one-parameter families of orbits for this coincidence, each orbit being five dimensional and having a single connected component.

The coincidence $\lambda = \mu$ is equivalent, under the trace-free condition, to $3\lambda + \nu = 0$. With $\lambda + 3\nu \not= 0$ and $\lambda \not= \nu$, the isotropy group is generated by $Z_1$, $Z_2$, $Z_3$ and $H := \{\,T_4(k): k \in {\bf R}\,\}$, where
$$T_4(g) := \pmatrix{1&-g&g&0\cr g&{2-g^2 \over 2}&g^2/2&0\cr g&-g^2/2&{2+g^2 \over 2}&0\cr 0&0&0&1\cr}.\eqno(4.10.7)$$
The situation is analogous to that of the previous paragraph.

The coincidence $\lambda=\nu$ is equivalent, under the trace-free condition, to $\lambda=\nu=-\mu$. With $3\lambda+\nu \not= 0$ and $\lambda + 3\nu \not= 0$, the isotropy subgroup consists of elements with matrix representations of the form
$$T_5 := \pmatrix{a&0&0&c\cr 0&\pm 1&0&0\cr 0&0&\pm 1&0\cr b&0&0&d\cr},\qquad \pmatrix{a&c\cr b&d\cr} \in {\underline {\bf O(2,2)}}\eqno(4.10.8)$$
and where the sign choices are paired. This one-parameter subgroup contains $Z_1$, $Z_2$ and $Z_3$. Hence, there are again two one-parameter families of orbits, with each orbit five dimensional and with a single connected component.

Finally, the coincidence $\lambda = \mu = \nu$, under the trace-free condition, is equivalent to $\lambda = \mu = \nu = 0$. The isotropy subgroup consists of elements with matrix representations of the form
$$\pmatrix{a&\mp F&\pm F&c\cr e&\pm\left[{2+k^2-e^2 \over 2}\right]&\pm\left[{e^2-k^2 \over 2}\right]&k\cr e&\pm\left[{k^2-e^2 \over 2}\right]&\pm\left[{2+e^2-k^2 \over 2}\right]&k\cr b&\mp G&\pm G&d\cr},\eqno(4.10.9{\rm a})$$
where all the signs choices are correlated and
$${\underline L} := \pmatrix{a&c\cr b&d\cr} \in {\underline {\bf O(2,2)}},\qquad\hbox{and}\qquad\pmatrix{F\cr G\cr} = \pmatrix{a&c\cr b&d\cr}\pmatrix{e\cr -k\cr}.\eqno(4.10.9{\rm b})$$
For $e$, and $k \in {\bf R}$, and $\underline L$, $F$ and $G$ as in (4.10.9b), with
$$T_6(\underline L,e,k) := \pmatrix{a&-F&F&c\cr e&{2+k^2-e^2 \over 2}&{e^2-k^2 \over 2}&k\cr e&{k^2-e^2 \over 2}&{2+e^2-k^2 \over 2}&k\cr b&-G&G&d\cr},\eqno(4.10.10{\rm a})$$
then $T_6(0,0,0) = 1$, and
$$T_6(\underline L,e,k)T_6(\hat{\underline L},\hat e,\hat k) = T_6(L\hat{\underline L},e\hat a + k\hat b + \hat e,e\hat c + k\hat d + \hat k).\eqno(4.10.11)$$
When $\underline L \in \underline{\bf SO^\bfplus(1,1)}$, by (1.5) one may confirm that $T_6(\underline L,e,k)$ represents, with respect to $\Psi$-ON bases, elements of ${\bf SO^\bfplus(2,2)}$. The elements $-Z_1$, $Z_2$ and $-Z_3$ take the form of $T_6(\underline L,0,0)$ for obvious choices of $\underline L$. The isotropy subgroup is therefore generated by elements of the form $T_6(\underline L,e,k)$ together with $-1$.

Hence, the isotropy subgroup has three parameters and intersects each component of {\bf O(2,2)}. Thus, for the coincidence of all eigenvalues, there is a single orbit, which is connected and of dimension three.\hfil\break
{\bf Type V}: for each $S$ of type V, there is a $\Psi$-ON basis with respect to which $S$ has matrix representation (4.5.5); $\tr(S) = 0$ is equivalent to $-2a = \lambda + \nu$. Let $\underline M(\lambda,\nu,b)$, $\lambda$, $\nu \in {\bf R}$ and $b \in {\bf R}^*$, denote the matrix in (4.5.5) with $-(\lambda + \nu)/2$ substituted for $a$. The elements of $\underline {\bf O(2,2)}$ that commute with (generic) $\underline M(\lambda,\nu,b)$ are $\pm 1$, $\pm Z_1$, $\pm Z_2$ and $\pm Z_3$. Thus, generically, there are four ($b \not= 0$, $\lambda \not= \nu$)) three-parameter families of orbits, each orbit being 6 dimensional and consisting of a single connected component.

The only coincidence possible in type V is $\lambda= \nu$, which under the trace-free condition is equivalent to $-a = \lambda = \nu$. The isotropy subgroup $K(\lambda,b)$ consists of elements with matrix representations of the form $T_5$ in (4.10.8). Hence, there are two ($b \not= 0$) two-parameter families of orbits, each orbit being five dimensional and consisting of a single connected component.\hfil\break
{\bf Type VI}: for each $S$ of type VI, there is a $\Psi$-ON basis with respect to which $S$ has matrix representation (4.5.7); $\tr(S) = 0$ is equivalent to $\lambda + a = 0$. Let $\underline M(\lambda,b)$ denote the matrix in (4.5.7) with $a = -\lambda$. The elements of $\underline {\bf O(2,2)}$ that commute with $\underline M(\lambda,b)$ are $\pm 1$ and $\pm Z_3$. Hence, there are two ($b \not= 0$) two-parameter families of orbits, each orbit of dimension six and with two connected components.\hfil\break
{\bf Type VII}: for each $S$ of type VII, there is a Witt basis with respect to which $S$ has matrix representation (4.5.8); $\tr(S) = 0$ is equivalent to $\lambda+\mu = 0$. Let $\underline M(\lambda)$ denote the matrix in (4.5.8) with $\mu = -\lambda$. The elements of $\underline {\bf O(4;hb)}$ that commute with $\underline M(\lambda)$, $\lambda \not= 0$, are $\pm 1$ and $\pm Z_4$, where
$$Z_4 := \pmatrix{1&0&0&0\cr 0&-1&0&0\cr 0&0&1&0\cr 0&0&0&-1\cr}.\eqno(4.10.12)$$
With respect to the Witt basis, $Z_4$ represents an orthogonal automorphism $T$, the matrix representation of which with respect to the $\Psi$-ON basis (4.5.9) is also $Z_4$, whence is an element of ${\bf O^\bfminus_\bfminus(2,2)}$. Hence, for $\lambda \not= 0$, there are two ($\lambda \not= 0$) one-parameter families of orbits, each orbit six dimensional and with two connected components.

Coincidence of eigenvalues $\lambda=\mu$ in the trace-free case amounts to $\lambda = \mu = 0$. In this case, the isotropy subgroup consists of elements of the form
$$T(A,a) := \pmatrix{A&aAJ\cr {\bf 0}_2&{^\tau\! A}^{-1}\cr},\eqno(4.10.13)$$
where $a \in {\bf R}$, $J := \left({0 \atop 1}{-1 \atop 0}\right)$, and
$$A \in \cases{{\underline {\bf O(2)}},&$\epsilon=\omega$;\cr \cr {\underline {\bf O(1,1)}},&$\epsilon=-\omega$.\cr}\eqno(4.10.14)$$
One computes that
$$T(A,a)T(D,d) = T\bigl(AD,d+a\det(D)\bigr)$$
whence $T(1,0) = 1$ and $[T(A,a)]^{-1} = T\bigl(A^{-1},-\det(A)a\bigr)$. Noting that $\det\bigl(T(A,a)\bigr) = \det(A)\det(A^{-1}) = 1$, then the isotropy subgroup is a subgroup of {\bf SO(4;hb)}. Since it contains the element of ${\bf O^\bfminus_\bfminus(2,2)}$ represented by $Z_4$ (by taking $A = \left({1 \atop 0}{0 \atop -1}\right)$ and $a=0$), then it is not a subgroup of the identity component of {\bf O(4;hb)}. Hence, there is a single orbit, of dimension four and with two components.\hfil\break
{\bf Type VIII}: for each $S$ of type VIII, there is a $\Psi$-ON basis with respect to which $S$ has matrix representation (4.5.11); $\tr(S) = 0$ is equivalent to $a+c = 0$. Note that $b \not= 0$, $d \not= 0$ in case VIII. Let $\underline M(b,d,a)$ denote the matrix in (4.5.11) with $c = -a$. The elements of $\underline {\bf O(2,2)}$ that commute with (generic) $\underline M(b,d,a)$, $a \not= 0$, are $\pm 1$ and $\pm Z_4$. As noted above, the matrix $Z_4$ represents, with respect to a $\Psi$-ON basis, an element of ${\bf O^\bfminus_\bfminus(2,2)}$. Since $b$, $d$, and $a$ are each nonzero, there are eight three-parameter families of orbits, each orbit of dimension six and with two connected components.

The eigenvalue coincidence in this case is $a=c$ and $b=d$, which in the trace-free case imposes $a=c=0$. The elements of $\underline {\bf O(2,2)}$ that commute with $\underline M(b,b,0)$ are precisely those of the form
$$\pmatrix{A&-B\cr B&A\cr},\eqno(4.10.15)$$
i.e., the isotropy subgroup is
$${\bf C}(2) \cap {\bf O(2,2)} = {\bf GL(2;C)} \cap {\bf O(2,2)} = {\bf O(2;C)}.$$
The group {\bf O(2;C)} has two components and is of real dimension two. It is a subgroup of {\bf SO(2,2)}. In particular, any element of {\bf SO(2;C)} can be written in the form $\left({\cos z \atop \sin z}{-\sin z \atop \cos z}\right)$, for which the corresponding form in (4.10.15) has, with $z = x+iy$,
$$A = \pmatrix{\cos x\cosh y&-\sin x\cosh y\cr \sin x\cosh y&\cos x\cosh y\cr} \hskip 1.25in B = \pmatrix{-\sin x\sinh y& -\cos x\sinh y\cr \cos x\sinh y&-\sin x\sinh y\cr},$$
whence ${\bf SO(2,C)} \leq {\bf SO^\bfplus(2,2)}$. Elements of ${\bf ASO(2;C)} := {\bf O(2;C)} \setminus {\bf SO(2;C)}$ can be written in the form $\left({\cos z \atop -\sin z}{-\sin z \atop -\cos z}\right)$, for which the corresponding form in (4.10.15) has, with $z = x+iy$,
$$A = \pmatrix{\cos x\cosh y&-\sin x\cosh y\cr -\sin x\cosh y&-\cos x\cosh y\cr} \hskip 1.25in B = \pmatrix{-\sin x\sinh y& -\cos x\sinh y\cr -\cos x\sinh y&\sin x\sinh y\cr},$$
whence ${\bf ASO(2,C)} \subset {\bf O^\bfminus_\bfminus(2,2)}$. Since $b \not= 0$, there are two one-parameter families of orbits, each orbit of dimension four and with two connected components.\hfil\break
{\bf Type IX}: for each $S$ of type IX, there is a $\Psi$-ON basis with respect to which $S$ has matrix representation
$$\underline M(\lambda,\mu,\nu,\sigma) := \pmatrix{\lambda&0&0&0\cr 0&\mu&0&0\cr 0&0&\nu&0\cr 0&0&0&\sigma\cr}.\eqno(4.10.16)$$
All descriptions to follow are with respect to this basis. Obviously $\tr(S) = 0$ is $\lambda+\mu+\nu+\sigma = 0$. Generically, no two eigenvalues are equal, which gives six conditions that partition the hyperplane $\lambda + \mu + \nu +\sigma = 0$ in ${\bf R}^4$ into 12 disconnected regions. Assuming these six inequalities and the trace-free condition, the elements of $\underline {\bf O(2,2)}$ that commute with (generic) $\underline M(\lambda,\mu,\nu,\sigma)$ form a discrete subgroup of order 16 generated by -1, $K_1$, $K_2$ and $K_4$. Hence, there are 24 (=4!)three-parameter families of orbits, each orbit of dimension six and connected.

There are six different possibilities of two coincident eigenvalues, which fall into two cases. In the first case, let $e^2 = f^2 = 1$ and $\left({a \atop b}{c \atop d}\right) \in {\underline  {\bf O(2)}}$. When $\lambda = \mu$, elements of the isotropy subgroup take the form
$$\pmatrix{a&c&0&0\cr b&d&0&0\cr 0&0&e&0\cr 0&0&0&f\cr}$$
and when $\nu = \sigma$, elements of the isotropy subgroup take the form
$$\pmatrix{e&0&0&0\cr 0&f&0&0\cr 0&0&a&c\cr 0&0&b&d\cr}.$$
In the second case, let $e^2=f^2=1$ and $\left({a \atop b}{c \atop d}\right) \in {\underline {\bf O(1,1)}}$. When $\lambda = \nu$, elements of the isotropy subgroup take the form
$$\pmatrix{a&0&c&0\cr 0&e&0&0\cr b&0&d&0\cr 0&0&0&f\cr};$$
when $\lambda = \sigma$, elements of the isotropy subgroup take the form
$$\pmatrix{a&0&0&c\cr 0&e&0&0\cr 0&0&f&0\cr b&0&0&d\cr};$$
when $\mu = \nu$, elements of the isotropy subgroup take the form
$$\pmatrix{e&0&0&0\cr 0&a&c&0\cr 0&b&d&0\cr 0&0&0&f\cr};$$
when $\mu = \sigma$, elements of the isotropy subgroup take the form
$$\pmatrix{e&0&0&0\cr 0&a&0&c\cr 0&0&f&0\cr 0&b&0&d\cr}.$$
In each case, one has the condition $\lambda + \mu +\nu +\sigma = 0$ and two coincident eigenvalues, so effectively three unequal parameters restricted to a hyperplane in ${\bf R}^3$, i.e., there are six two-parameter families of orbits, each orbit of dimension five and connected.

There are three ways of having exactly two pairs of coincident eigenvalues. When $\lambda=\mu$ and $\nu=\sigma$, elements of the isotropy subgroup take the form $\left({A \atop {\bf 0}_2}{{\bf 0}_2 \atop B}\right)$, with $A$, $B \in {\underline {\bf O(2)}}$. With $\left({a_1 \atop b_1}{c_1 \atop d_1}\right)$ and $\left({a_2 \atop b_2}{c_2 \atop d_2}\right)$ in ${\underline{\bf O(1,1)}}$, when $\lambda = \nu$ and $\mu = \sigma$, elements of the isotropy subgroup take the form
$$\pmatrix{a_1&0&c_1&0\cr 0&a_2&0&c_2\cr b_1&0&d_1&0\cr 0&b_2&0&d_2\cr},$$
while when $\lambda=\sigma$ and $\mu=\nu$ elements of the isotropy subgroup take the form
$$\pmatrix{a_1&0&0&c_1\cr 0&a_2&c_2&0\cr 0&b_2&d_2&0\cr b_1&0&0&d_1\cr}.$$
In effect, there are two non-coincident parameters subject to one linear condition; hence, in each case there are two one-parameter families of orbits, each orbit being four dimensional and connected.

There are four cases of exactly three coincident eigenvalues. Let $A := (a_{ij}) \in {\bf R}(3)$ and $f^2=1$. When $\lambda=\mu=\nu$, elements of the isotropy subgroup take the form
$$\pmatrix{a_{11}&a_{12}&a_{13}&0\cr a_{21}&a_{22}&a_{23}&0\cr a_{31}&a_{32}&a_{33}&0\cr 0&0&0&f\cr},$$
and when $\lambda=\mu=\sigma$, elements of the isotropy subgroup take the form
$$\pmatrix{a_{11}&a_{12}&0&a_{13}\cr a_{21}&a_{22}&0&a_{23}\cr 0&0&f&0\cr a_{31}&a_{32}&0&a_{33}\cr},$$
with $A \in {\underline{\bf O(2,1)}}$ in each case. When $\lambda=\nu=\sigma$, elements of the isotropy subgroup take the form
$$\pmatrix{a_{11}&0&a_{12}&a_{13}\cr 0&f&0&0\cr a_{21}&0&a_{22}&a_{23}\cr a_{31}&0&a_{32}&a_{33}\cr},$$
while when $\mu=\nu=\sigma$, elements of the isotropy subgroup take the form
$$\pmatrix{f&0&0&0\cr 0&a_{11}&a_{12}&a_{13}\cr
0&a_{21}&a_{22}&a_{23}\cr 0&a_{31}&a_{32}&a_{33}\cr},$$
with $A \in {\underline{\bf O(1,2)}}$ in both cases. Hence, in each case, there are two one-parameter families of orbits, each orbit of dimension three and connected.

Finally, when all four eigenvalues coincide in the trace-free case, the zero endomorphism results.
\vskip 24pt
\noindent {\section 5. Algebraic Classification of the Ricci Spinor}
\vskip 12pt
For any tensor $S_{ab}$ over ${\bf R}^{2,2}$, by PRI (3.3.54)
$$S_{ab} = S_{(AB)(A'B')} + \lambda_{AB}\epsilon_{A'B'} + \epsilon_{AB}\mu_{A'B'} + \rho\epsilon_{AB}\epsilon_{A'B'},$$
where $\lambda_{AB}$ and $\mu_{A'B'}$ are symmetric. Then, $S_{ab}$ is traceless iff $\rho=0$ and $S_{ab}$ is symmetric iff
$$\lambda_{AB}\epsilon_{A'B'} + \epsilon_{AB}\mu_{A'B'} = \lambda_{BA}\epsilon_{B'A'} + \epsilon_{BA}\mu_{B'A'},$$
i.e., $\lambda_{AB}\epsilon_{A'B'} + \epsilon_{AB}\mu_{A'B'} = 0$. Hence, for a traceless, symmetric tensor $\Phi_{ab}$ the spinor form $\Phi_{ABA'B'}$ satisfies 
$$\Phi_{ABA'B'} = \Phi_{(AB)(A'B')}.\eqno(5.1)$$
A spinor $\Phi_{A\ldots LM'\ldots V'} = \Phi_{(A\ldots L)(M'\ldots V')}$, with $p$ unprimed indices, $q$ primed indices, will be called a $(p,q)$-spinor. The expression
$$P_\Phi(\xi^A,\zeta^{M'}) := \Phi_{A\ldots LM'\ldots V'}\xi^A\ldots\xi^L\zeta^{M'}\ldots\zeta^{V'},\eqno(5.2)$$
defines a polynomial, homogeneous of degrees $p$ and $q$ in $\xi^A$ and $\zeta^{A'}$, respectively. The zero locus $\omega$ lies on ${\bf RP}^1 \times {\bf RP}^1$, which space may be viewed as the quadric surface that is the projective null cone of ${\bf R}^{2,2}$ in ${\bf RP}^3$. As usual, it is convenient to complexify, allowing the spinors to come from ${\bf C}S$ and ${\bf C}S'$. Then (5.2) defines a zero locus $\Omega$ on ${\bf CP}^1 \times {\bf CP}^1$, viewed as the quadric surface that is the projective null cone of ${\bf C}^4_{\rm E}$ in ${\bf CP}^3$.

If $q=0$, say, $P_\Phi$ reduces to a degree $p$ homogeneous polynomial and defines a zero locus on ${\bf KP}^3$. By the fundamental theorem of algebra, this zero locus in ${\bf CP}^3$ consists of $p$ points (counted with multiplicity). This case (and the alternative $p=0$) is the basis of the algebraic classification of the Weyl spinors (PRII, Law 2006). The generalization to $(p,q)$-spinors is given in PRII, \S 8.7, and amounts to classifying the irreducible factors of the zero locus of (5.2). If the polynomial $P_\Phi$ factors, $P_\Phi=QR$ say, then each of $R$ and $Q$ must be independently homogeneous in each of $\xi^A$ and $\zeta^{M'}$ to respect the fact that $P_\Phi$ is. Thus, $Q$ and $R$ determine $(r,s)$-spinors $\Psi_{A\ldots DM' \ldots R'}$ and $\chi_{E\ldots LS'\ldots V'}$ such that
$$\Phi_{A\ldots L}^{M'\ldots V'} = \Psi_{(A\ldots D}^{(M'\ldots R'}\chi_{E\ldots L)}^{S'\ldots V')}.$$
Now restrict attention to real (2,2)-spinors, i.e., (5.1). Their algebraic classification is based on the possible $(r,s)$-spinors into which $\Phi_{ABA'B'}$ can factor. Note that while $\Phi_{ABA'B'}$ is itself real, its factors need not be, but if complex must occur as complex conjugate pairs. Table One lists the possibilities for $\Phi_{ABA'B'}$ for ${\bf R}^{2,2}$.
\vskip 24pt
$$\vbox{\offinterlineskip
\halign{&\strut\ #\ \cr
\multispan{3}\hfil\bf Table 1\hfil\cr
\noalign{\medskip}
\noalign{\hrule}
\noalign{\vskip 2pt}
\hfil\bf Spinor Factorization\hfil&\hfil\bf Notation\hfil&\hfil\bf DF\hfil\cr
\noalign{\vskip 2pt}
\noalign{\hrule}
\noalign{\vskip 2pt}
\hfil$\Phi_{ABA'B'}$\hfil&\hfil$(2,2)$\hfil&\hfil$9$\hfil\cr
\noalign{\vskip 2pt}
\noalign{\hrule}
\noalign{\vskip 2pt}
\hfil$\Psi^{(A'}_{AB}\Upsilon^{B')}$\hfil&\hfil$(2,1)(0,1)$\hfil&\hfil$7$\hfil\cr
\noalign{\vskip 2pt}
\noalign{\hrule}
\noalign{\vskip 2pt}
\hfil$\Psi_{(A}\Upsilon^{A'B'}_{B)}$\hfil&\hfil$(1,0)(1,2)$\hfil&\hfil$7$\hfil\cr
\noalign{\vskip 2pt}
\noalign{\hrule}
\noalign{\vskip 2pt}
\hfil$\Psi^{(A'}_{(A}\Upsilon^{B')}_{B)}$\hfil&\hfil$(1,1)(1,1)$\hfil&\hfil$7$\hfil\cr
\noalign{\vskip 2pt}
\noalign{\hrule}
\noalign{\vskip 2pt}
\hfil$\pm\Gamma^{(A'}_{(A}\bar\Gamma^{B')}_{B)}$\hfil&\hfil$(1,1)\overline{(1,1)}$\hfil&\hfil$7$\hfil\cr
\noalign{\vskip 2pt}
\noalign{\hrule}
\noalign{\vskip 2pt}
\hfil$\Psi_{(A}\Lambda_{B)}^{(A'}\Upsilon^{B')}$\hfil&\hfil$(1,0)(1,1)(0,1)$\hfil&\hfil$6$\hfil\cr
\noalign{\vskip 2pt}
\noalign{\hrule}
\noalign{\vskip 2pt}
\hfil$\Psi_{(A}\Sigma_{B)}\Lambda^{(A'}\Upsilon^{B')}$\hfil&\hfil$(1,0)(1,0)(0,1)(0,1)$\hfil&\hfil$5$\hfil\cr
\noalign{\vskip 2pt}
\noalign{\hrule}
\noalign{\vskip 2pt}
\hfil$\Gamma_{(A}\bar\Gamma_{B)}\Lambda^{(A'}\Upsilon^{B')}$\hfil&\hfil$(1,0)\overline{(1,0)}(0,1)(0,1)$\hfil&\hfil$5$\hfil\cr
\noalign{\vskip 2pt}
\noalign{\hrule}
\noalign{\vskip 2pt}
\hfil$\Psi_{(A}\Sigma_{B)}\Xi^{(A'}\bar\Xi^{B')}$\hfil&\hfil$(1,0)(1,0)(0,1)\overline{(0,1)}$\hfil&\hfil$5$\hfil\cr
\noalign{\vskip 2pt}
\noalign{\hrule}
\noalign{\vskip 2pt}
\hfil$\pm\Gamma_{(A}\bar\Gamma_{B)}\Xi^{(A'}\bar\Xi^{B')}$\hfil&\hfil$(1,0)\overline{(1,0)}(0,1)\overline{(0,1)}$\hfil&\hfil$5$\hfil\cr
\noalign{\vskip 2pt}
\noalign{\hrule}
\noalign{\vskip 2pt}
\hfil$\pm\Gamma_{(A}\bar\Gamma_{B)}\Upsilon^{A'}\Upsilon^{B'}$\hfil&\hfil$(1,0)\overline{(1,0)}(0,1)^2$\hfil&\hfil$4$\hfil\cr
\noalign{\vskip 2pt}
\noalign{\hrule}
\noalign{\vskip 2pt}
\hfil$\pm\Psi_A\Psi_B\Xi^{(A'}\bar\Xi^{B')}$\hfil&\hfil$(1,0)^2(0,1)\overline{(0,1)}$\hfil&\hfil$4$\hfil\cr
\noalign{\vskip 2pt}
\noalign{\hrule}
\noalign{\vskip 2pt}
\hfil$\Psi_{(A}\Sigma_{B)}\Upsilon^{A'}\Upsilon^{B'}$\hfil&\hfil$(1,0)(1,0)(0,1)^2$\hfil&\hfil$4$\hfil\cr
\noalign{\vskip 2pt}
\noalign{\hrule}
\noalign{\vskip 2pt}
\hfil$\Psi_A\Psi_B\Lambda^{(A'}\Upsilon^{B')}$\hfil&\hfil$(1,0)^2(0,1)(0,1)$\hfil&\hfil$4$\hfil\cr
\noalign{\vskip 2pt}
\noalign{\hrule}
\noalign{\vskip 2pt}
\hfil$\pm\Psi^{(A'}_{(A}\Psi^{B')}_{B)}$\hfil&\hfil$(1,1)^2$\hfil&\hfil$4$\hfil\cr
\noalign{\vskip 2pt}
\noalign{\hrule}
\noalign{\vskip 2pt}
\hfil$\pm\Psi_A\Psi_B\Upsilon^{A'}\Upsilon^{B'}$\hfil&\hfil$(1,0)^2(0,1)^2$\hfil&\hfil$3$\hfil\cr
\noalign{\vskip 2pt}
\noalign{\hrule}
\cr}}$$
\vskip 24pt
The first column lists the spinor form of the factorization; the second column the notation, adapted from PRII, that I will employ to describe that form, viz., denoting an $(r,s)$-spinor by $(r,s)$ and a complex conjugate pair of $(r,s)$-spinors by $(r,s)\overline{(r,s)}$; the third column lists the degrees of freedom in the specified form. The kernel symbols $\Psi$, $\Sigma$, $\Lambda$, and $\Upsilon$ here all denote real spinors (of varying rank) while $\Gamma$ and $\Xi$ here denote nontrivially complex spinors (of varying rank).

In Figures One and Two on the next page, the specializations that are possible from one factorization to another are diagrammed. To obtain the correct global topology in a single diagram, the two figures must be glued together by identifying the labels that occur in both figures (which, for convenience, occur in bold type). The degrees of freedom are listed on the far left and apply to all types at that horizontal level. In Figure One, the implication arrows do not denote specialization; rather, they indicate that the form $(2,0)(0,1)(0,1)$, say, is actually one of two possibilities, either $(1,0)\overline{(1,0)}(0,1)(0,1)$ or $(1,0)(1,0)(0,1)(0,1)$, according as the $(2,0)$-spinor factors into real $(1,0)$-spinors or a complex conjugate pair; either way the factorization is immediate, not a further specialization.

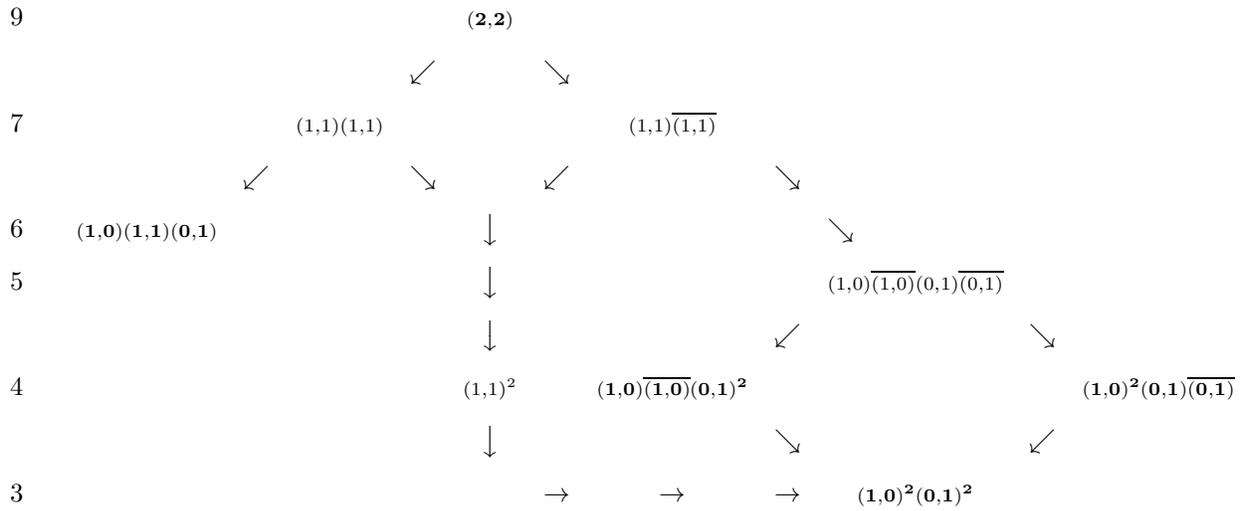
\vfill\eject
\vskip 12pt
\centerline{\bf Figure 1}
\vskip 12pt
$$
\matrix{9&&&&&&{\scriptstyle{\bf (2,2)}}&&&&\cr
&&&&&\swarrow&&\searrow&&&\cr
7&&&&{\scriptstyle(2,1)(0,1)}&&&&{\scriptstyle(1,0)(1,2)}&&\cr
&&&&\big\downarrow&\searrow&&\swarrow&\big\downarrow&&\cr
6&&&&\big\downarrow&&{\scriptstyle{\bf (1,0)(1,1)(0,1)}}&&\big\downarrow&&\cr
&&&&\big\downarrow&&\big\downarrow&&\big\downarrow&&\cr
5&&{\scriptstyle(1,0)\overline{(1,0)}(0,1)(0,1)}&\Leftarrow&{\scriptstyle(2,0)(0,1)(0,1)}&\Rightarrow&{\scriptstyle(1,0)(1,0)(0,1)(0,1)}&\Leftarrow&{\scriptstyle(1,0)(1,0)(0,2)}&\Rightarrow&{\scriptstyle(1,0)(1,0)(0,1)\overline{(0,1)}}\cr
&&\big\downarrow&&&\swarrow&&\searrow&&&\big\downarrow\cr
4&&{\scriptstyle {\bf (1,0)\overline{(1,0)}(0,1)^2}}&&{\scriptstyle (1,0)(1,0)(0,1)^2}&&&&{\scriptstyle (1,0)^2(0,1)(0,1)}&&{\scriptstyle {\bf (1,0)^2(0,1)\overline{(0,1)}}}\cr
&&\big\downarrow&&&\searrow&&\swarrow&&&\big\downarrow\cr
3&&&\rightarrow&&&{\scriptstyle{\bf (1.0)^2(0,1)^2}}&&&\leftarrow&\cr}$$
\vskip 60pt
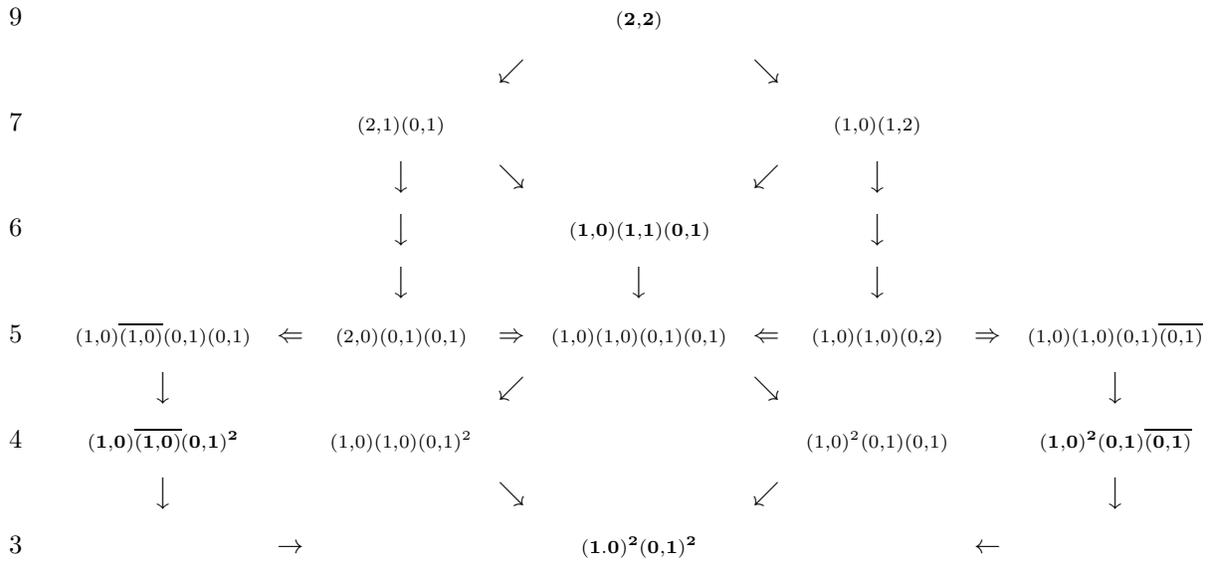
\centerline{\bf Figure 2}
\vskip 12pt
$$
\matrix{9&&&&&&$${\scriptstyle{\bf (2,2)}}&&&&&&\cr
&&&&&\swarrow&&\searrow&&&&&\cr
7&&&&{\scriptstyle(1,1)(1,1)}&&&&{\scriptstyle(1,1)\overline{(1,1)}}&&&&\cr
&&&\swarrow&&\searrow&&\swarrow&&\searrow&&&\cr
6&&{\scriptstyle {\bf(1,0)(1,1)(0,1)}}&&&&\big\downarrow&&&&\searrow\hfill&&\cr
5&&&&&&\big\downarrow&&&&{\scriptstyle (1,0)\overline{(1,0)}(0,1)\overline{(0,1)}}&&\cr
&&&&&&\big\downarrow&&&\swarrow&&\searrow&\cr
4&&&&&&{\scriptstyle (1,1)^2}&&{\scriptstyle {\bf (1,0)\overline{(1,0)}(0,1)^2}}&&&&{\scriptstyle {\bf (1,0)^2(0,1)\overline{(0,1)}}}\cr
&&&&&&\big\downarrow&&&\searrow&&\swarrow&\cr
3&&&&&&&\rightarrow&\rightarrow&\rightarrow&{\scriptstyle {\bf (1,0)^2(0,1)^2}}&&\cr}$$
\vfill\eject

This diagram is somewhat more complicated than the analogue in the case of ${\bf R}^{1,3}$; see PRII, Table (8.3.3) and PRII (8.8.1) for the analogue of Table One. As noted in PRII, p. 277, this classification of Ricci spinor types based on factorization into $(r,s)$-spinors (reducibility), is coarser than the classification of the corresponding traceless self-adjoint endomorphisms described in \S 4. Relating the classification of \S4 to that of \S5 is of obvious interest, but will also help refine the classification of $\Phi_{ABA'B'}$ to a comparable degree as that of \S 4. As Penrose has shown, e.g., PRII \S 8.8, this refinement involves the algebraic geometry of the zero locus defined by (5.2); in particular, null eigenvectors of $\Phi^a{}_b$ define points on the zero locus of (5.2) and the behaviour of the zero locus at these points is related to the classification of \S 4.
\vskip 24pt
\noindent {\section 6. Relating the Classifications}
\vskip 12pt
For traceless self-adjoint endomorphisms $\Phi^a{}_b$, I obtained in \S 4 (see also Table Two below) 10 types (counting IIIa and IIIb as distinct) and various cases of coincidence of eigenvalues; each of these distinct cases will hereafter be called subtypes. For each type, there is fixed matrix form that represents elements of $\End^0_+({\bf R}^{2,2})$ of that type with respect to some $\Psi$-ON or Witt basis. This form and the associated basis can also be used to describe the subtypes within a type. In this section, I take, for each type, the fixed matrix form $\Phi^{\bf a}{}_{\bf b}$, say, use the relevant scalar product to convert this to a fully covariant form, and contract over the dual basis to obtain the tensor form $\Phi_{ab}$. By using (3.48--50) for the appropriate kind of basis ($\Psi$-ON or Witt), one can express the resulting object in terms of spinors and obtain first $\Phi_{ABA'B'}$ and then the {\sl Ricci polynomial} (5.2). In this way, the reducibility of $\Phi_{ABA'B'}$ will be obtained.

When the Ricci spinor is reducible, the {\sl Ricci locus} $\Omega$, defined by the vanishing of the Ricci polynomial $P_\Phi$ on ${\bf CP}^1 \times {\bf CP}^1$, has multiple components, the intersections of which are singular points of the Ricci locus, leading to the expectation that the Ricci spinor type is related to the nature of the singularities of the Ricci locus. Of course, singularities may also arise from self intersections. The Lorentzian case is presented in PRII \S\S 8.7--8.8. Suppose $Q = \mu^A\eta^{A'}$ lies on the Ricci locus $\Omega$. Choosing spin frames $\{\mu^A,\nu^A\}$ and $\{\eta^{A'},\theta^{A'}\}$ for ${\bf C}S$ and ${\bf C}S'$, respectively, write 
$$\xi^A = X\mu^A + Y\nu^A \hskip 1.25in\zeta^{A'} = U\eta^{A'} + V\theta^{A'}.\eqno(6.1)$$
Introducing affine coordinates $y := Y/X$ and $v := V/U$, then $Q = (0,0)$ and
$$f(y,v) :=  P_\Phi(1,y,1,v),\eqno(6.2)$$
is an affine description of the Ricci polynomial whose zero locus coincides with $\Omega$ on a neighbourhood of $X = (0,0)$. For any line $L$, $y = \sigma t$, $v = \tau t$, passing through $(0,0)$, the intersections of $L$ and $\Omega$ (near $(0,0)$) are given by
$$\eqalignno{0 = F(t) &:= f(\sigma t,\tau t)&(6.3)\cr
&= (f_y\sigma+f_v\tau)t + {1 \over 2!}(f_{yy}\sigma^2 + 2f_{yv}\sigma\tau + f_{vv}\tau^2)t^2 + \cdots + {1 \over n!}\left(\sum_{r=0}^n\,{n \choose r}{\partial^nf \over \partial y^{n-r}\partial v^r}\sigma^{n-r}\tau^r\right)t^n + \cdots.\cr}$$
where all partial derivatives are evaluated at $(0,0)$.
The point $Q = (0,0)$ on $\Omega$ is nonsingular iff $(f_y,f_v) \not= (0,0)$, in which case the tangent at $Q$ is given, in this affine picture, by the the unique line through $(0,0)$ with at least second-order contact with $f^{\dashv}(\{0\})$, i.e., such that $F'(0) = 0$ (with $F''(0) \not= 0$ for precisely second-order contact), i.e., the line $L$ solving
$$F'(0) = f_y\sigma + f_y\tau = 0.\eqno(6.4)$$

On the other hand, $Q$ is singular iff $(f_y,f_v) = (0,0)$. In general, if all partial derivatives of $f$ at $(0,0)$ up to and including order $n-1$, but not all $n$'th order partial derivatives, are zero, then $(0,0)$ is called a {\sl singular point of order $n$} (or $n$-{\sl fold singularity}). In this case, all lines through $(0,0)$ have at least $n$'th-order contact with $\Omega$ at $(0,0)$, i.e., $F(0) = F'(0) = \cdots = F^{(n-1)}(0) = 0$, while the lines $L$ satisfying
$$F^{(n)}(0) = \sum_{r=0}^n\,{n \choose r}{\partial^nf \over \partial y^{n-r}\partial v^r}\sigma^{n-r}\tau^r = 0,\eqno(6.5)$$
are called {\sl tangents} to $\Omega$ at $(0,0)$ and have at least $(n+1)$'th order contact with $\Omega$ at $(0,0)$.

In particular, when $Q$ is singular but $(f_{yy},f_{yv},f_{vv}) \not= (0,0,0)$, $Q$ is called a {\sl double point}. The lines $L$ solving
$$F''(0) = f_{yy}\sigma^2 + 2f_{yv}\sigma\tau + f_{vv}\tau^2 = 0,\eqno(6.6)$$
are the tangents to $\Omega$ at $Q$ and have at least third-order contact with $\Omega$ at $Q$. A double point is called a {\sl node} if there are two distinct such tangents. Note that if $Q$ is a real point of $\Omega$, i.e., lies on $\omega$, then the tangents are either real, or complex conjugate and $Q$ is an isolated point of $\omega$. When (6.6) has a repeated root (necessarily real for the Ricci polynomial as it is real), i.e., the two tangents coincide, the double point is called a {\sl cusp}. If the tangent to a cusp has at least fourth-order contact, i.e., solves  (6.5) for $n=3$, then $X$ is called a {\sl tacnode}. See Walker (1962) for further details.
\vskip 24pt
\noindent {\bf 6.7 Lemma}\hfil\break
$Q$ is an $n$-fold singularity of $\Omega$ iff all the $(n-1)$'st, but not all the $n$'th, order partial derivatives of $P_\Phi$ vanish at $Q$.

Proof. A straightforward adaptation of Walker (1962), Theorem 2.4, p. 55. One uses the definition of $f$ in terms of $P_\Phi$ and Euler equations for $P_\Phi$ and its derivatives:
$$\xi^A{\partial P_\Phi \over \partial\xi^A} = 2P_\Phi = \zeta^{A'}{\partial P_\Phi \over \partial\zeta^{A'}},\qquad  \xi^A{\partial \over \partial\xi^A}\left({\partial P_\Phi \over \partial\xi^B}\right) = {\partial P_\Phi \over \partial\xi^B} \qquad  \xi^A{\partial \over \partial\xi^A}\left({\partial P_\Phi \over \partial\zeta^{B'}}\right) = 2{\partial P_\Phi \over \partial\zeta^{B'}},\eqno(6.8)$$
analogous equations for derivatives with respect to $\zeta^{A'}$ and similar equations for higher order partial derivatives, each evaluated at $X = (0,0)$, to relate the partial derivatives of $f$ at $(0,0)$ and those of $P_\Phi$ at $X$.\bull
\vskip 24pt
This result is useful, but also ensures that the nature of singularities is geometrical, being independent of the choice of affine representation and the choice of homogeneous coordinates. Since $P_\Phi$ is homogeneous of degree two in each of $\xi^A$ and $\zeta^{A'}$, (6.7) also indicates the orders of singularities expected since partial derivatives of order three or more in either $\xi^A$ or $\zeta^{A'}$ vanish automatically. Note that the space on which the Ricci polynomial is defined is a product of projective spaces; its `homogeneous coordinates are of the form $([X,Y],[U,V])$.

From (6.7), the condition for a singular point $Q$ on $\Omega$ can be written
$$\left({\partial P_\Phi \over \partial\xi^A}\right)\vert_{Q = \mu^A\nu^{A'}} = 0 \hskip 1.25in \left({\partial P_\Phi \over \partial\zeta^{A'}}\right)\vert_{Q = \mu^A\nu^{A'}} = 0,\eqno(6.9)$$
i.e.,
$$0 = \left({\partial P_\Phi \over \partial\xi^A}\right)\vert_{Q = \mu^A\nu^{A'}} = 2\Phi_{ABA'B'}\mu^B\nu^{A'}\nu^{B'} \hskip .5in 0 = \left({\partial P_\Phi \over \partial\zeta^{A'}}\right)\vert_{Q = \mu^A\nu^{A'}} = 2\Phi_{ABA'B'}\mu^A\mu^B\nu^{B'}.\eqno(6.10)$$
By two applications of PRI (3.5.16), (6.10) implies
$$\mu_A\sigma_{A'} = \Phi_{ABA'B'}\mu^B\nu^{B'} = \gamma_A\nu_{A'},$$
for some $\gamma_A$ and $\sigma_{A'}$, whence, by PRI (3.5.2),
$$\Phi_{ABA'B'}\mu^B\nu^{B'} = \chi\mu_A\nu_{A'},\eqno(6.11)$$
for some $\chi$, i.e., $\mu^A\nu^{A'}$ is a null eigenvector of $\Phi^a{}_b$. As (6.11) implies (6.10), it follows that singular points of $\Omega$ correspond to null eigenvectors of $\Phi^a{}_b$. This fact ensures the algebraic geometry of $\Omega$ is indeed relevant to the classification. From \S 4, one observes that generic null eigenvectors (null in ${\bf C}^{2,2} \cong {\bf C}^4_{\rm E} = {\bf C}({\bf R}^{2,2})$) occur only when $m > 1$ and $M < m$ for the eigenvalue. Additionally, when two distinct eigenvalues with eigenvectors of opposite character in the generic situation come into coincidence, the resulting eigenspace must contain null vectors.

From (6.10), e.g., by applying PRI (3.5.27), one deduces
$$\Phi_{ABA'B'}\mu^A\mu^B = \beta\nu_{A'}\nu_{B'} \hskip 1.25in \Phi_{ABA'B'}\nu^{A'}\nu^{B'} = \alpha\mu_A\mu_B,\eqno(6.12)$$
for some $\alpha$ and $\beta$. As
$$\displaylines{\left({\partial^2 P_\Phi \over \partial\xi^A\partial\xi^B}\right)\vert_{Q=\mu^A\nu^{A'}} = 2\Phi_{ABA'B'}\nu^{A'}\nu^{B'} = 2\alpha\mu_A\mu_B\cr
\noalign{\vskip 6pt}
\hfill \left({\partial^2 P_\Phi \over \partial\xi^A\partial\zeta^{A'}}\right)\vert_{Q=\mu^A\nu^{A'}} = 4\Phi_{ABA'B'}\mu^B\nu^{B'} = 4\chi\mu_A\nu_{A'}\hfill\llap(6.13)\cr
\noalign{\vskip 6pt}
\left({\partial^2 P_\Phi \over \partial\zeta^{A'}\partial\zeta^{B'}}\right)\vert_{Q=\mu^A\nu^{A'}} = 2\Phi_{ABA'B'}\mu^A\mu^B = 2\beta\nu_{A'}\nu_{B'},\cr}$$
one computes that $f_{yy} = 2\alpha$, $f_{yv} = 4\chi$, and $f_{vv} = 2\beta$, whence the discriminant of (6.6) is
$$4\chi^2 -\alpha\beta.\eqno(6.14)$$
Hence, the nature of a double point can be determined from the purely algebraic equations in (6.13) via (6.14).

Now, assuming $\Phi^a{}_b$ admits a null eigenvector $\mu^A\nu^{A'}$ with eigenvalue $\chi$, by (4.1)(ii) the eigenvectors of a distinct eigenvalue must lie in $\langle \mu^A\nu^{A'} \rangle^\perp = \langle \mu^A\nu^{A'},\mu^A\sigma^{A'},\rho^A\nu^{A'} \rangle_{\bf C}$ (within ${\bf C}^{2,2}$), where $\rho^A\mu_A = 1$, and $\sigma^{A'}\nu_{A'} = 1$. As $W := \langle \mu^A\nu^{A'} \rangle_{\bf C}$ is $\Phi^a{}_b$ invariant, since $\Phi^a{}_b$ is self-adjoint, $W^\perp$ is invariant too. Hence, using (6.13), one computes that, with respect to the null tetrad $\{\mu^A\nu^{A'},\rho^A\nu^{A'},\rho^A\sigma^{A'}\mu^A\sigma^{A'}\}$, $\Phi^a{}_b$ has matrix representation
$$\pmatrix{\chi&f&a&e\cr 0&-\chi&-e&-\beta\cr 0&0&\chi&0\cr 0&-\alpha&-f&-\chi\cr},$$
for some $a$, $e$ and $f$, for which the eigenvalue equation is $(\chi-\kappa)^2[(\chi+\kappa)^2-\alpha\beta] =0$.
Hence, in addition to $\chi$ itself (which is of algebraic multiplicity at least two, as noted previously), the other eigenvalues are
$$\kappa = -\chi \pm \sqrt{\alpha\beta}.\eqno(6.15)$$
One sees that these eigenvalues are in fact distinct to $\chi$ iff $4\chi^2 \not= \alpha\beta$, i.e., by (6.14), iff the null eigenvector $\mu^A\nu^{A'}$ does not define a cusp (or higher singularity) on $\Omega$. Conversely, a cusp (or higher order singularity) requires that $\chi$, the eigenvalue of the null eigenvector at which the singularity occurs, must have algebraic multiplicity at least three.

The preceding results provide some information on the relation between singularities and their nature and the eigenvalue and eigenvector structure of $\Phi^a{}_b$. I now turn to explicit descriptions of each (sub)type. The following observation will be useful.
\vskip 24pt
\noindent {\bf 6.16 Observation}\hfil\break
For a given type of (traceless) self-adjoint endomorphism $\Phi^a{}_b$, one particular example has the standard matrix representation of that type with respect to the standard $\Psi$-ON (or Witt) basis. All others in that type have the same matrix representation with respect to another $\Psi$-ON (or Witt) basis. When recasting $\Phi_{ab}$ as $\Phi_{ABA'B'}$ by substituting null tetrads, one must bear in mind the different relationships between null tetrads and $\Psi$-ON (and Witt) bases, depending on which coset of ${\bf SO^\bfplus(2,2)}$ $L$ belongs to, where $L$ maps the standard $\Psi$-ON (or Witt) basis to the basis in question. These relationships were given in (3.56), (3.60), and (3.64). Here I elaborate their consequences for the form of the Ricci polynomial. If a $\Psi$-ON basis is the image of the standard basis by $L \in {\bf O(2,2)}$, then of course the Witt basis associated to the former is the image of the standard Witt basis by $L$. If a $\Psi$-ON basis is the image of the standard basis by an element of ${\bf SO^\bfplus(2,2)}$, then the associated Witt basis $\{E_1,E_2,E_3,E_4\}$ is such that $\{E_1,E_2,E_3,-E_4\}$ is a null tetrad (with spin frames related to the standard spin frames by an element of ${\bf SL(2;R)} \times {\bf SL(2;R)}$).

Suppose $\{U^a,V^a,X^a,Y^a\}$ is a $\Psi$-ON basis that is the image of the standard basis by an element $L$ of ${\bf O^\bfminus_\bfminus(2,2)}$ (I must rely on context to prevent confusion of the notation for the elements $U^a$, $V^a$, $X^a$, and $Y^a$ of a $\Psi$-ON basis with that for homogeneous coordinates $([X,y],[U,V])$ for ${\bf CP}^1 \times {\bf CP}^1$); then the associated Witt basis $\{E_1,E_2,E_3,E_4\}$ is such that $\{E_3,-E_4,E_1,E_2\}$ is a null tetrad (with spin frames related to the standard spin frames by the composition of an element of ${\bf SL(2;R)} \times {\bf SL(2;R)}$ with $o^A \leftrightarrow \iota^A$ and $o^{A'} \leftrightarrow \iota^{A'}$). Hence, if $\check{\Phi}_{ABA'B'}$ corresponds to the $\Phi^a{}_b$ that has the type's form as matrix representation with respect to the standard $\Psi$-ON or Witt basis, as appropriate, to obtain the form of $\Phi_{ABA'B'}$ corresponding to the $\Phi^a{}_b$ that has the type's form as matrix representation with respect to either $\{U^a,V^a,X^a,Y^a\}$ or $\{E_1,E_2,E_3,E_4\}$, as appropriate, one can apply the transformation $o^A \leftrightarrow \iota^A$ and $o^{A'} \leftrightarrow \iota^{A'}$ to $\check{\Phi}_{ABA'B'}$ (and drop the checks over the elements of the standard spin frames). 

With 
$$\xi^A = Xo^A + Y\iota^A \hskip 1.25in \zeta^{A'} = Uo^{A'} + V\iota^{A'},\eqno(6.16.1))$$
the Ricci polynomial is $P_\Phi(X,Y,U,V) = \Phi_{ABA'B'}\xi^A\xi^B\zeta^{A'}\zeta^{B'}$. To obtain the Ricci polynomial of the $\Phi_{ABA'B'}$ of the previous paragraph from that of $\check{\Phi}_{ABA'B'}$, note that (6.16.1) is left unchanged, whence $\xi^Ao_A = Y$ in $\check{\Phi}_P$ is replaced by $\xi^A\iota_A = -X$, $\xi^A\iota_A = -X$ is replaced by $\xi^Ao_A = Y$, $\zeta^{A'}o_{A'} = V$ is replaced by $\zeta^{A'}\iota_{A'} = -U$ and $\zeta^{A'}\iota_{A'} = -U$ is replaced by $\zeta^{A'}o_{A'} = V$, i.e., one obtains $P_\Phi$ from $\check{P}_\Phi$ by
$$Y\ \leftrightarrow\ -X \hskip 1.25in V\ \leftrightarrow\ -U.\eqno(6.16.2)$$

Suppose now $L \in {\bf O^\bfminus_\bfplus(2,2)}$; then $\{E_1,-E_4,E_3,E_2\}$ is a null tetrad (with spin frames related to the standard spin frames by the composition of an element of ${\bf SL(2;R)} \times {\bf SL(2;R)}$ with $o^A \leftrightarrow o^{A'}$ and $\iota^A \leftrightarrow \iota^{A'}$). By the same reasoning as in the previous case, one can obtain the Ricci polynomial in this case by making the substitutions
$$Y\ \leftrightarrow\ V \hskip 1.25in X\ \leftrightarrow\ U,\eqno(6.16.3)$$
in $\check{P}_\Phi$.

Similarly, if $L \in {\bf O^\bfplus_\bfminus(2,2)}$, then $\{E_3,E_2,E_1,-E_4\}$ is a null tetrad (with spin frames related to the standard spin frames by the composition of an element of ${\bf SL(2;R)} \times {\bf SL(2;R)}$ with $o^A \leftrightarrow \iota^{A'}$ and $\iota^A \leftrightarrow o^{A'}$). One can obtain the Ricci polynomial in this case by making the substitutions
$$Y\ \leftrightarrow\ -U \hskip 1.25in X\ \leftrightarrow\ -V,\eqno(6.16.4)$$
in $\check{P}_\Phi$.

It will prove convenient, in the following, to refer to $\Psi$-ON (Witt) bases that are the image of the standard (Witt) basis by an element $L$ in a given component of {\bf O(2,2)} as being of {\sl the O-type} of that component, e.g., if $L \in {\bf O^\bfminus_\bfminus(2,2)}$, I shall say the basis is of O-type ${\bf O^\bfminus_\bfminus}$. Since the different Ricci polynomials that occur for endomorphisms of a given (sub)type vary only with the O-type of the basis with respect to which the endomorphism has the standard matrix form of the (sub)type, and these different forms are obtained from each other by (simple linear) changes of homogeneous coordinates, all Ricci polynomials of a given(sub)type have, geometrically, the same singularity structure. It will therefore suffice to consider the form of the Ricci polynomial that occurs for those endomorphisms of the (sub)type that take the (sub)type's standard matrix form with respect to bases of O-type ${\bf O^\bfplus_\bfplus}$.
\vskip 24pt
\noindent {\bf 6.17 Type I}\hfil\break
For each traceless, self-adjoint endomorphism $\Phi^a{}_b$ of type I, there is a Witt basis $\{E_1,E_2,E_3,E_4\}$ with respect to which $\Phi^a{}_b$ has matrix representation (4.4.3) with $\lambda=0$. One readily computes that
$$\Phi^{ab} = \epsilon E^a_2E^b_2 + 2E^a_{(1}E^b_{4)}.\eqno(6.17.1)$$
From (4.10), the isotropy subgroup is $\{\pm1\}$, so there are distinct elements of type I for each component of {\bf O(2,2)}. For those with matrix representation (4.4.3) ($\lambda=0$) with respect to Witt bases of O-type ${\bf O^\bfplus_\bfplus}$,
$$\Phi_{ABA'B'} = \epsilon\iota_A\iota_Ao_{A'}o_{B'} - o_Ao_B(o_{A'}\iota_{B'} + \iota_{A'}o_{B'}),\eqno(6.17.2)$$
and 
$$P_\Phi(X,Y,U,V) = \epsilon X^2V^2 + 2UVY^2.\eqno(6.17.3)$$
Hence, the possible Ricci polynomials are, by O-type of the associated Witt basis:
$$P_\Phi = \cases{(\epsilon X^2V+2UY^2)V,&${\bf O^\bfplus_\bfplus}$;\cr (\epsilon Y^2U+2VX^2)U,&${\bf O^\bfminus_\bfminus}$;\cr (\epsilon U^2Y+2XV^2)Y,&${\bf O^\bfminus_\bfplus}$;\cr (\epsilon V^2X+2YU^2)X,&${\bf O^\bfplus_\bfminus}$.\cr}\eqno(6.17.4)$$
Observe that $\Phi_{ABA'B'}$ is of Ricci spinor type $(2,1)(0,1)$, for Witt bases of O-type {\bf SO}, and of type $(1,0)(1,2)$, for Witt bases of O-type {\bf ASO}.
For the form (6.17.2), $\mu^A\nu^{A'} = o^Ao^{A'}$ is a null eigenvector (with zero eigenvalue) and one computes, from (6.13) that $\alpha=0$ and $\beta=\epsilon$. Hence, the null eigenvector is a double (but not triple) point and a cusp. Note that, as expected, the zero eigenvalue has algebraic multiplicity at least three. Explicitly, introducing affine coordinates $y := Y/X$ and $v := V/U$, the singularity occurs at $(0,0)$. The affine description of the Ricci locus near $(0,0)$ is the zero set of
$$f(y,v) = (\epsilon v+2y^2)v.$$
In this case (6.3) is
$$0 = F(t) = f(\sigma t,\tau t) = \tau^2\epsilon t^2 + 2\sigma^2\tau\epsilon t^3,$$
and $(0,0)$ is clearly a double point. The tangents are given by $\tau=0$ and coincide, and moreover they make the coefficient of $t^3$ vanish, so the $y$-axis is the tangent at $(0,0)$ and has contact of order at least four, i.e., $(0,0)$ is a tacnode. In fact, the $y$-axis is a linear component of the Ricci locus in this affine picture, tangent at $(0,0)$ to a quadratic component; these components correspond to the $(0,1)$ and $(2,1)$ factors of the $(2,1)(0,1)$-factorization of the Ricci polynomial.

Analogous descriptions and results are obtained for the Ricci polynomials valid for Witt bases of the other O-types, which are listed in (6.17.4), because, as noted in (6.16), they differ from each other only by (simple linear) changes of homogeneous coordinates. I will not reiterate this point.
\vskip 24pt
\noindent {\bf 6.18 Type II}\hfil\break
For each traceless, self-adjoint endomorphism $\Phi^a{}_b$ of type II, there is a Witt basis $\{E_1,E_2,E_3,E_4\}$ with respect to which $\Phi^a{}_b$ has matrix representation (4.4.5) with $a=0$. One readily computes
$$\Phi^{ab} = E^a_1E^b_1 - E^a_2E^b_2 + 2b(E_1^{(a}E_4^{b)} - E_2^{(a}E_3^{b)}).\eqno(6.18.1)$$
From (4.10), the isotropy subgroup is $\{\pm1\}$, so there are distinct elements of type II for each component of {\bf O(2,2)}. For those with matrix representation (4.4.5) (with $\lambda=0$) with respect to Witt bases of O-type ${\bf O^\bfplus_\bfplus}$,
$$\Phi_{ABA'B'} = (o_Ao_B - \iota_A\iota_B)o_{A'}o_{B'} - b(o_Ao_B + \iota_A\iota_B)(o_{A'}\iota_{B'} + \iota_{A'}o_{B'}),\eqno(6.18.2)$$
with Ricci polynomial
$$P_\Phi(X,Y,V,W) = [X^2(2bU-V) + Y^2(2bU+V)]V.\eqno(6.18.3)$$
Hence, the possible Ricci polynomials are, by O-type of the associated Witt basis:
$$P_\Phi = \cases{[X^2(2bU-V) + Y^2(2bU+V)]V,&${\bf O^\bfplus_\bfplus}$;\cr [Y^2(2bV-U)+X^2(2bU+V)]U,&${\bf O^\bfminus_\bfminus}$;\cr [U^2(2bX-Y)+V^2(2bX+Y)]Y,&${\bf O^\bfminus_\bfplus}$;\cr [V^2(2bY-X)+U^2(2bY+X)]X,&${\bf O^\bfplus_\bfminus}$.\cr}\eqno(6.17.4)$$
Observe that $\Phi_{ABA'B'}$ is of Ricci spinor type $(2,1)(0,1)$, for Witt bases of O-type {\bf SO} or of type $(1,0)(1,2)$, for Witt bases of O-type {\bf ASO}.

Now $\Phi^a{}_b$ has complex null eigenvector $w^a = E_1^a + iE_2^a$, which for bases of O-type ${\bf O^\bfplus_\bfplus}$ takes the form
$$w^a = E_1^a + iE_2^a  = (o^A + i\iota^A)o^{A'} =: \mu^A\nu^{A'},\eqno(6.18.3)$$
and its conjugate. From (6.13), one readily computes for $w^a$:
$$\chi = ib \hskip .75in \alpha = 0 \hskip .75in \beta = -2.\eqno(6.18.4)$$
Thus, $4\chi^2 - \alpha\beta = -4b^2 < 0$, i.e., the complex null eigenvector $w$ is a node. Similarly for $\bar w$. Explicitly, with the definitions (6.16.1), the eigenvectors have homogeneous coordinates $([X,Y],[U,V]) = ([1,\pm i],[1,0])$. Taking affine coordinates $y := Y/X$ and $v := V/U$, the eigenvectors have coordinates $(\pm i,0)$. Consider then the lines given by $y = \pm i +\sigma t$ and $v = \tau t$. Then, in place of (6.3), one obtains
$$F(t) = \tau[2(-\tau\pm 2ib\sigma)t^2 + 2\sigma(b\sigma \pm i\tau)t^3 + \sigma^2\tau t^4],$$
and 
$$F(0) = 0 = F'(0), \hskip .5in F''(0) = 4\tau(-\tau \pm 2ib\sigma), \hskip .5in F'''(0) = 12\sigma\tau(b\sigma \pm i\tau)\hskip .5in F^{(4)}(0) = 24\sigma^2\tau^2.$$
Hence, the complex null eigenvectors $w$ and $\bar w$ are nodes, with tangents satisfying $\tau = 0$ and $\tau = \mp 2ib\sigma$ (signs respectively). The former is the $y$-axis and is tangent to both complex null eigenvectors and moreover makes $F'''(0) = F^{(4)}(0)$ vanish. Indeed, the $y$-axis corresponds in this affine picture to the linear factor in the Ricci polynomial, which defines a linear component of the Ricci locus that is tangent to both nodes. Being a component of the locus, it has infinite order of contact with the locus.
\vskip 24pt
\noindent {\bf 6.19 Types IIIa and IIIb}\hfil\break
For each generic traceless, self-adjoint endomorphism $\Phi^a{}_b$ of type IIIa, there is a $\Psi$-ON basis\hfil\break
$\{U^a,V^a,X^a,Y^a\}$ with respect to which $\Phi^a{}_b$ has matrix representation (4.5.1) with $\lambda = -3\mu\not= 0$. One readily computes
$$\Phi^{ab} = -\mu(3U^aU^b - V^aV^b + X^aX^b + Y^aY^b) - \sqrt2(V^{(a}X^{b)} + X^{(a}Y^{b)}),\eqno(6.19.1)$$
From (4.10), the isotropy subgroup is discrete but intersects ${\bf O^\bfminus_\bfplus(2,2)}$. For those $\Phi^a{}_b$ with matrix representation (4.5.1) ($\lambda=-3\mu\not=0$) with respect to $\Psi$-ON bases of O-type ${\bf O^\bfplus_\bfplus}$,
$$\eqalignno{\Phi_{ABA'B'} &= -\mu\bigl[2o_Ao_Bo_{A'}o_{B'} + 2\iota_A\iota_B\iota_{A'}\iota_{B'} + (o_A\iota_B + \iota_Ao_B)(o_{A'}\iota_{B'}+\iota_{A'}o_{B'})\bigr]\cr
&\qquad - {1 \over \sqrt{2}}\bigl[(o_A\iota_B + \iota_Ao_B)o_{A'}o_{B'} - \iota_A\iota_B(o_{A'}\iota_{B'} + \iota_{A'}o_{B'})\bigr],&(6.19.2)\cr}$$
with Ricci polynomial
$$P_\Phi(X,Y,U,V) = -2\mu(YV+XU)^2 - \sqrt{2}X(XU-YV)V.\eqno(6.19.3)$$
Hence, the possible Ricci polynomials are, by O-type of the associated $\Psi$-ON basis:
$$P_\Phi + \cases{-2\mu(YV+XU)^2-\sqrt{2}X(XU-YV)V,&${\bf O^\bfplus_\bfplus}$;\cr -2\mu(XU+YV)^2-\sqrt{2}Y(YV-XU)U,&${\bf O^\bfminus_\bfminus}$;\cr -2\mu(VY+XU)^2-\sqrt{2}Y(XU-YV)U,&${\bf O^\bfminus_\bfplus}$;\cr -2\mu(UX+YV)^2-\sqrt{2}X(YV-XU)V,&${\bf O^\bfplus_\bfminus}$.\cr}\eqno(6.19.4)$$
Recall that $Z_1$ (4.10.1) defines an element of ${\bf O^\bfminus_\bfplus(2,2)}$ in the isotropy subgroup of type IIIa. In particular, (6.19.1) is invariant under this transformation. I use this case to explicate the significance of the isotropy subgroup for the algebraic geometry. A generic $\Phi^a{}_b$ that has the matrix representation (4.5.1) (with $\lambda = -3\mu \not= 0$) with respect to some $\Psi$-ON basis $F$ of O-type ${\bf O^\bfminus_\bfplus}$ also has the same matrix representation with respect to the $\Psi$-ON basis of O-type ${\bf O^\bfplus_\bfplus}$ that is related to $F$ by $Z_1$, and thus the distinct generic $\Phi^a{}_b$ of type IIIa are given by taking the matrix form (4.5.1) (with $\lambda = -3\mu \not= 0$) with respect to bases of O-type ${\bf O^\bfplus_\bfplus} \amalg {\bf O^\bfplus_\bfminus} =: {\bf O^\bfplus}$. For a given $\Phi^a{}_b$, although its expression (6.19.1) is invariant with respect to the nontrivial element of the isotropy subgroup, the spinor forms are not, due to the different relationships between Witt bases of different O-types and null tetrads, summarized in (6.16); consequently the Ricci polynomials are not identical, though they are equivalent as regards their geometrical content. In any event, observe that the Ricci spinor type is $(2,2)$ ($\mu\not=0$).

For (6.19.1), the only null eigenvector is $V^a+Y^a$ (with eigenvalue $\mu$), which for (6.19.2) is $\tilde m^a$. By (6.13),
$$\chi = \mu \hskip .75in \alpha = -2\mu \hskip .75in \beta = -2\mu,$$
whence $4\chi^2 - \alpha\beta = 0$, i.e., the real null eigenvector defines a cusp. Explicitly, the homogeneous coordinates for $\tilde m^a$ are $([0,1],[1,0])$ so take affine coordinates $x := X/Y$ and $v := V/U$. The affine version of (6.19.3) is $f(x,v) = -2\mu(v+x)^2 - \sqrt{2}x(x-v)v$ and (6.3) is
$$F(t) = -2\mu(\tau + \sigma)^2t^2 - \sqrt{2}\sigma\tau(\sigma-\tau)t^3,$$
whence
$$F(0) = 0 = F'(0) \hskip .5in F''(0) = -4\mu(\tau+\sigma)^2 \hskip .5in F'''(0) = -6\sqrt{2}\sigma\tau(\sigma-\tau),$$
i.e., $(0,0)$ has two coincident tangents (given by $\sigma = -\tau$), making $(0,0)$ a cusp (but not a tacnode as $\sigma = -\tau$ does not make $F'''(0)$ vanish).

When $\mu=0$, $P_\Phi$ reduces to Ricci spinor type $(1,0)(1,1)(0,1)$; in the affine picture of the previous paragraph, $f(x,v) = -\sqrt2x(x-v)v$, i.e., consists of three distinct lines through $(0,0)$, i.e., the Ricci locus consists of three components, which intersect at the real null eigenvector, making it a triple point. The coincidence $\lambda=\mu$ creates no new null eigenvectors.

Type IIIb is completely analogous. I just record forms specific to the type. For each generic traceless, self-adjoint endomorphism $\Phi^a{}_b$ of type IIIb, there is a $\Psi$-ON basis $\{U^a,V^a,X^a,Y^a\}$ with respect to which $\Phi^a{}_b$ has matrix representation (4.5.2) with $\mu = -3\lambda$ and $\lambda \not= 0$. One readily computes that
$$\Phi^{ab} = \lambda(U^aU^b + V^aV^b - X^aX^b + 3Y^aY^b) + \sqrt2(U^{(a}V^{b)} + V^{(a}X^{b)}),\eqno(6.19.5)$$
and, for $\Psi$-ON bases of O-type ${\bf O^\bfplus_\bfplus}$,
$$\eqalignno{\Phi_{ABA'B'} &= \lambda\bigl[2\iota_A\iota_Bo_{A'}o_{B'} + 2o_Ao_B\iota_{A'}\iota_{B'} + (o_A\iota_B + \iota_Ao_B)(o_{A'}\iota_{B'} + \iota_{A'}o_{B'})\bigr]\cr
&\qquad + {1 \over \sqrt{2}}\bigl[(o_A\iota_B + \iota_Ao_B)o_{A'}o_{B'} - o_Ao_B(o_{A'}\iota_{B'} + \iota_{A'}o_{B'}).\bigr]&(6.19.6)\cr}$$
The Ricci polynomials are
$$P_\Phi = \cases{2\lambda(XV+UY)^2 + \sqrt{2}Y(YU-XV)V,&${\bf O^\bfplus_\bfplus}$;\cr 2\lambda(YU+VX)^2 + \sqrt{2}X(XV-YU)U,&${\bf O^\bfminus_\bfminus}$;\cr 2\lambda(YU+VX)^2 + \sqrt{2}Y(VX-YU)V,&${\bf O^\bfminus_\bfplus}$;\cr 2\lambda(XV+UY)^2 + \sqrt{2}X(YU-XV)U,&${\bf O^\bfplus_\bfminus}$.\cr}\eqno(6.19.7)$$
\vskip 24pt
\noindent {\bf 6.20 Type IV}\hfil\break
For each generic traceless, self-adjoint endomorphism $\Phi^a{}_b$ of type IV, there is a $\Psi$-ON basis $\{U^a,V^a,X^a,Y^a\}$ with respect to which $\Phi^a{}_b$ has matrix representation (4.5.4) with $\lambda + 2\mu + \nu = 0$ but with no coincidences amongst $\lambda$, $\mu$, and $\nu$. There is a real null eigenvector, with eigenvalue $\mu$, but the matrix (4.5.4) and its associated $\Psi$-ON basis is not well adapted to describing the singularity at this null eigenvector. Instead, consider the basis $\{v_1,v_2,v_3,v_4\}$ that gives the JCF $J_1(\lambda) \oplus J_2(\mu) \oplus J_1(\nu)$. Putting
$$E_1 := {v_1 + v_4 \over \sqrt{2}} \qquad E_2 := v_2 \qquad E_3 := {v_1 - v_4 \over \sqrt{2}} \qquad E_4 = \epsilon v_3,\eqno(6.20.1)$$
gives a Witt basis, with respect to which $S^a{}_b$ has matrix representation
$$\pmatrix{{\lambda + \nu \over 2}&0&{\lambda - \nu \over 2}&0\cr 0&\mu&0&\epsilon\cr {\lambda - \nu \over 2}&0&{\lambda+\nu \over 2}&0\cr 0&0&0&\mu\cr}.\eqno(6.20.2)$$
Hence, with respect to this Witt basis, a generic traceless, self-adjoint endomorphism $\Phi^a{}_b$ of type IV has matrix representation
$$\pmatrix{{\lambda + \nu \over 2}&0&{\lambda - \nu \over 2}&0\cr 0&-{\lambda+\nu \over 2}&0&\epsilon\cr {\lambda - \nu \over 2}&0&{\lambda+\nu \over 2}&0\cr 0&0&0&-{\lambda+\nu \over 2}\cr}.\eqno(6.20.3)$$
One readily computes
$$\Phi^{ab} = \left({\lambda-\nu \over 2}\right)(E_1^aE_1^b + E_3^aE_3^b) + \left({\lambda+\nu \over 2}\right)(2E_1^{(a}E_3^{b)} - 2E_2^{(a}E_4{^b)}) + \epsilon E_2^aE_2^b.\eqno(6.20.4)$$
From (4.10), the isotropy subgroup is discrete but intersects each component of {\bf O(2,2)}, so it suffices to consider Witt bases of O-type ${\bf O^\bfplus_\bfplus}$ to describe the distinct $\Phi^a{}_b$. One obtains
$$\Phi_{ABA'B'} = \left({\lambda-\nu \over 2}\right)[o_Ao_Bo_{A'}o_{B'} + \iota_A\iota_B\iota_{A'}\iota_{B'}] + \left({\lambda+\nu \over 2}\right)[(o_A\iota_B + \iota_Ao_B)(o_{A'}\iota_{B'} + \iota_{A'}o_{B'})] + \epsilon\iota_A\iota_Bo_{A'}o_{B'},\eqno(6.20.5)$$
and hence
$$P_\Phi(X,Y,U,V) = \left({\lambda-\nu \over 2}\right)(Y^2V^2+X^2U^2) + (\lambda+\nu)[2XYUV] + \epsilon X^2V^2,\eqno(6.20.6)$$
which has Ricci spinor type $(2,2)$ in the generic case.

The real null eigenvector is $v_2 = E_2^a = \tilde m^a = \iota^Ao^{A'}$. From (6.13), one computes
$$\chi = \mu \hskip .75in \alpha = {\lambda - \nu \over 2} = \beta,\eqno(6.20.7)$$
whence
$$4\chi^2 - \alpha\beta = {3\lambda^2 + 10\lambda\nu + 3\nu^2 \over 4} = {(3\lambda+\nu)(\lambda+3\nu) \over 4}.\eqno(6.20.8)$$
The numerator on the rhs vanishes for
$$\lambda = -3\nu (\ \Leftrightarrow\ \mu = \nu)\qquad\hbox{ and }\qquad \nu = -3\lambda (\ \Leftrightarrow\ \mu = \lambda).\eqno(6.20.9)$$
Thus, for the generic case of Type IV, the real null eigenvector is a node, with two real/complex tangents according as $(3\lambda+\nu(\lambda+3\nu)\ {> \atop <}\ 0$. For the coincidences of eigenvalues indicated in (6.20.9), the double point is at least a cusp.

Explicitly, the singularity has homogeneous coordinates $([0,1],[1,0])$, so choose affine coordinates $x = X/Y$ and $v = V/U$. Then,
$$f(x,v) = {\lambda \over 2}(v^2 + 4xv + x^2) - {\nu \over 2}(v^2 - 4xv + x^2) + \epsilon x^2v^2,$$
and (6.3) is
$$F(t) = {1 \over 2}\left[4(\lambda + \nu)\sigma\tau + (\lambda-\nu)(\sigma^2 + \tau^2)\right]t^2 + \epsilon\sigma^2\tau^2t^4.$$
Hence,
$$0 = F(0) = F'(0) = F'''(0), \hskip .5in F''(0) = 4(\lambda + \nu)\sigma\tau + (\lambda-\nu)(\sigma^2 + \tau^2), \hskip .5in F^{(4)}(0) = 24\epsilon\sigma^2\tau^2.$$
Thus, $(0,0)$ is a double point, its nature determined by the discrimninant of the expression for $F''(0)$ regarded as a quadratic in $\tau/\sigma$. This discriminat is
$$4(\lambda+\nu)^2 - (\lambda -\nu)^2 = 3\lambda^2 + 10\lambda\nu + 3\nu^2,$$
which confirms the result above.

Now consider coincidences of eigenvalues. Suppose $\lambda = \mu$ ($\ \Leftrightarrow\ \nu = -3\lambda$). Substituting $\nu = -3\lambda$ in the above computations, by (6.20.9), the singularity at the real null eigenvector is at least a cusp; since $F'''(0)$ is identically zero, technically, the cusp is a tacnode. The Ricci polynomial is
$$\eqalignno{P_\Phi &= 2\lambda(YV - XU)^2 + \epsilon X^2V^2&(6.20.10)\cr
& =\cases{\left[\sqrt{2\lambda}(YV-XU) - iXV\right]\left[\sqrt{2\lambda}(YV-XU) + iXV\right],&$\epsilon = 1$;\cr \left[\sqrt{2\lambda}(YV-XU) - XV\right]\left[\sqrt{2\lambda}(YV-XU) + XV\right],&$\epsilon = -1$;\cr}}$$
i.e., the Ricci spinor type is $(1,1)\overline{(1,1)}$ if $\epsilon = 1$ and $(1,1)(1,1)$ if $\epsilon = -1$.

The coincidnce $\nu = \mu$ ($\ \Leftrightarrow\ \lambda = -3\nu$) is analogous to the previous case. The double point is again a tacnode, and the Ricci polynomial is
$$\eqalignno{P_\Phi &= -2\nu(YV+XU)^2 + \epsilon X^2V^2&(6.20.11)\cr
&= \cases{\left[XV - \sqrt{2\nu}(YV+XU)\right]\left[XV + \sqrt{2\nu}(YV+XU)\right],&$\epsilon = 1$;\cr -\left[XV - i\sqrt{2\nu}(YV+XU)\right]\left[XV + i\sqrt{2\nu}(YV+XU)\right],&$\epsilon = -1$;\cr}}$$
i.e., the Ricci spinor type is $(1,1)\overline{(1,1)}$ if $\epsilon = -1$ and $(1,1)(1,1)$ if $\epsilon = 1$. One easily checks that $E_2^a$ is a point of intersection of the components of the Ricci locus determined by two quadratic factors of the Ricci polynomial in each case of (6.20.10--11). Neither this, nor the previous, coincidence creates new null eigenvectors.

For the coincidence $\lambda = \nu$, the Ricci polynomial is
$$P_\Phi = X(4\lambda YU + \epsilon XV)V,\eqno(6.20.12)$$
i.e., of Ricci spinor type $(1,0)(1,1)(0,1)$. As $\lambda = \nu$, $3\lambda^2 + 10\lambda\nu + 3\nu^2 = 16\lambda^2 > 0$, so $E_2^a$ is a node with real tangents. From the affine description above, the tangents are given by $\sigma = 0$ and $\tau=0$, i.e., the tangents are determined by the linear factors of the Ricci polynomial, i.e., $E_2^a$ is the point of intersection of these two linear components. Moreover, the coincidence $\lambda=\nu$ entails that the first and last summands in the decomposition of type IV span a single eigenspace; in particular, $E_1^a$ and $E_3^a$ are each real null eigenvectors (with eigenvalue $\lambda$). $E_1^a$ has homogeneus coordinates $([1,0],[1,0])$; choosing affine coordinates $y = Y/x$ and $v = V/U$, one obtains $f(y,v) = (4\lambda y + \epsilon v)v$ and $F(t) = (4\lambda\sigma\tau + \epsilon\tau^2)t^2$. One readily computes that $E_1^a$ is a node and that its tangents are given by $\tau = 0$ and $\epsilon\tau + 4\lambda\sigma = 0$. The first is the line $v=0$ and the second the line $\epsilon v + 4\lambda y = 0$, which in homogeneous coordinates are $V = 0$ and $4\lambda YU + \epsilon XY = 0$, respectively, i.e., the tangents are determined by these two components of the Ricci locus. Similarly, one finds that $E_3^a$ is a point of intersection of the linear component $X = 0$ and the quadratic component. Thus, the three components define, by their intersections, three nodes.

Finally, the triple coincidence $\lambda = \mu = \nu = 0$ yields $\Phi^{ab} = \epsilon E_2^aE_2^b$, whence $E_1^a$, $E_2^a$ and $E_3^a$ are each null eigenvectors of zero. The Ricci locus is
$$P_\Phi = \epsilon V^2X^2,\eqno(6.20.13)$$
i.e., consists of a pair of double lines. Note that $\langle E_1^a,E_2^a \rangle_{\bf K}$ and $\langle E_2^a,E_3^a \rangle_{\bf K}$ are each totally null eigensubspaces when $\lambda = \mu = \nu$. For Witt bases of O-type ${\bf O^\bfplus_\bfplus}$, their elements take the form $\alpha^Ao^{A'}$ and $\iota^A\beta^{A'}$, respectively, for arbitrary $\alpha^A$ and $\beta^{A'}$, i.e., constitute the components $V = 0$ and $X = 0$ respectively. Each of these components is double, whence singular, reflecting that the corresponding points are null eigenvectors.
\vskip 24pt
\noindent {\bf 6.21 Type V}\hfil\break
For each generic traceless, self-adjoint endomorphism $\Phi^a{}_b$ of type V, there is a $\Psi$-ON basis $\{U^a,V^a,X^a,Y^a\}$ with respect to which $\Phi^a{}_b$ has matrix representation (4.5.5) with $\lambda + 2a + \nu = 0$, but $\lambda \not= \nu$. One readily computes
$$\Phi^{ab} = \lambda U^aU^b - \nu Y^aY^b + \left({\lambda + \nu \over 2}\right)[X^aX^b - V^aV^b] - b[2V^{(a}X^{b)}].\eqno(6.21.1)$$
From (4.10), the isotropy subgroup intersects each component of {\bf O(2,2)} so it suffices to consider $\Psi$-ON bases of O-type ${\bf O^\bfplus_\bfplus}$ to describe the distinct $\Phi^a{}_b$. Hence, 
$$\eqalignno{\Phi_{ABA'B'} &= \left({3\lambda+\nu \over 4}\right)(o_Ao_Bo_{A'}o_{B'} + \iota_A\iota_B\iota_{A'}\iota_{B'}) - \left({\lambda+3\nu \over 4}\right)(\iota_A\iota_Bo_{A'}o_{B'} + o_Ao_B\iota_{A'}\iota_{B'})\cr
& +\left({\lambda - \nu \over 4}\right)(o_A\iota_B + \iota_Ao_B)(o_{A'}\iota_{B'}+\iota_{A'}o_{B'})&(6.21.2)\cr 
&- {b \over 2}\left[(o_A\iota_B + \iota_Ao_B)(o_{A'}o_{B'} + \iota_{A'}\iota_{B'}) - (o_Ao_B + \iota_A\iota_B)(o_{A'}\iota_{B'}+\iota_{A'}o_{B'})\right],\cr}$$
whence the Ricci polynomial is
$$P_\Phi(X,Y,U,V) = \left({3\lambda+\nu \over 4}\right)(YV+XU)^2 - \left({\lambda+3\nu \over 4}\right)(XV+YU)^2 + b\left[XY(V^2+U^2) - (X^2+Y^2)UV\right],\eqno(6.21.3)$$
which is Ricci spinor type $(2,2)$.

There no null eigenvectors, whence no singularities (the complex eigenvector $w = V + iX$, and its conjugate, are null in $\overline{\bf C}^{2,2}$ but not in ${\bf C}^{2,2})$.

With the coincidence $\lambda=\nu$, however, the first and last summands in the type V decomposition span a single eigenspace, so $U \pm Y$ are now real null eigenvectors. One computes
$$\Phi_{ABA'B'} = \lambda\epsilon_{AB}\epsilon_{A'B'} - {b \over 2}\left[(o_A\iota_B + \iota_Ao_B)(o_{A'}o_{B'} + \iota_{A'}\iota_{B'}) - (o_Ao_B + \iota_A\iota_B)(o_{A'}\iota_{B'}+\iota_{A'}o_{B'})\right],\eqno(6.21.4)$$
and
$$\eqalignno{P_\Phi &= \lambda\left[(YV+XU)^2 - (XV+YU)^2\right] + b\left[XY(V^2+U^2) - (X^2+Y^2)UV\right]\cr
&= \lambda\left[(YV-XU)^2 - (XV-YU)^2\right] + b(XV-YU)(VY-UX)&(6.21.5)\cr
&= \cases{\lambda[k^{-1}(YV-XU) - (XV-YU)][k(YV-XU) + (XV-YU)],& when $\lambda \not= 0$;\cr b(XV-YU)(VY-UX),& when $\lambda = 0$;\cr}\cr}$$
where
$$k := {-b \pm \sqrt{b^2+4\lambda^2} \over 2\lambda} \in {\bf R},\qquad \lambda \not = 0.\eqno(6.21.6)$$
Thus, the Ricci spinor type is $(1,1)(1,1)$ (further factorization when $\lambda \not= 0$ requires $k=\pm 1$, which is equivalent to $\lambda=0$). One can take the null eigenvectors to be 
$$N_+ := (o^A + \iota^A)(o^{A'} + \iota^{A'}) \hskip 1.25in N_- := (o^A - \iota^A)(o^{A'} - \iota^{A'}).\eqno(6.21.7)$$
Using (6.21.4), (6.13) yields,
$$\chi = \lambda\qquad \alpha = -b = -\beta\hbox{  for $N_+$} \hskip 1in \chi = \lambda\qquad \alpha = b = -\beta\hbox{  for $N_-$},\eqno(6.21.8)$$
whence, in each case $4\chi^2 - \alpha\beta = 4\lambda^2 + b^2 > 0$; each null eigenvector is a node with two real tangents.

For an affine picture, as $N_+$ has homogeneous coordinates $([1,1],[1,1])$, choose affine coordinates $x = X/Y$ and $u = U/V$ and consider lines given by $x = 1 + \sigma t$ and $u = 1 + \tau t$, i.e., that pass through $(1,1)\ \leftrightarrow N_+$. One has
$$f(x,u) = \lambda\left[(1-xu)^2 - (x-u)^2\right] + b(x-u)(1-xu)$$
and 
$$F(t) = \left[4\lambda\sigma\tau - b(\sigma^2-\tau^2)\right]t^2 + \sigma\tau\left[2\lambda(\sigma+\tau)-b(\sigma-\tau)\right]t^3 + \lambda\sigma^2\tau^2t^4,$$
from which one readily computes that the nature of the double point is determined by solving $4\lambda\sigma\tau - b(\sigma^2-\tau^2)$ as a quadratic in $\sigma/\tau$, say, which gives
$${\sigma \over \tau} = {2\lambda \pm 2\sqrt{4\lambda^2 + b^2} \over b},$$
(cf. (6.21.6)),confirming that the double point is a node with real tangents. Moreover, one can confirm that these tangents do not make $F'''(0)$ vanish. Taking affine coordinates $y = Y/X$ and $v = V/U$ gives an affine picture of $N_-$ with coordinates $(-1,-1)$ for which the affine picture of the locus near $(-1,-1)$ takes the same form as that just given for $N_+$, thus yielding the same result.
\vskip 24pt
\noindent {\bf 6.22 Type VI}\hfil\break
For this type, there is a real null eigenvector of $\lambda$; namely, $v_1$ in the basis $\{v_1,v_2,v_3,v_4\}$ giving the JCF in (4.5). Hence, rather than using the form (4.5.7), consider the Witt basis
$$E_1 := v_1 \qquad E_2 := {v_3 + v_4 \over \sqrt{2}} \qquad E_3 = \epsilon v_2 \qquad E_4 := {v_3 - v_4 \over \sqrt{2}},\eqno(6.22.1)$$
with respect to which $S^a{}_b$ has matrix representation
$$\pmatrix{\lambda&0&\epsilon&0\cr 0&a&0&-b\cr 0&0&\lambda&0\cr 0&b&0&a\cr}.\eqno(6.22.2)$$
Hence, a generic, traceless, self-adjoint endomorphism $\Phi^a{}_b$ of type VI has matrix representation, with respect to (6.22.1),
$$\pmatrix{\lambda&0&\epsilon&0\cr 0&-\lambda&0&-b\cr 0&0&\lambda&0\cr 0&b&0&-\lambda\cr},\eqno(6.22.3)$$
$\lambda\not=0$, and one computes that
$$\Phi^{ab} = \epsilon E_1^aE_1^b - b(E_2^aE_2^b - E_4^aE_4^b) + 2\lambda(E_1^{(a}E_3^{b)} - E_2^{(a}E_4{^b)}).\eqno(6.22.3)$$
For Witt bases of O-type ${\bf O^\bfplus_\bfplus}$, one computes
$$\Phi_{ABA'B'} = \epsilon o_Ao_Bo_{A'}o_{B'} - b(\iota_A\iota_Bo_{A'}o_{B'} - o_Ao_B\iota_{A'}\iota_{B'}) + \lambda(o_A\iota_B + \iota_Ao_B)(o_{A'}\iota_{B'}+ \iota_{A'}o_{B'}),\eqno(6.22.4)$$
and the Ricci polynomial is
$$P_\Phi(X,Y,U,V) = \epsilon Y^2V^2 - b(X^2V^2 - Y^2U^2) + 4\lambda XYUV.\eqno(6.22.5)$$
Hence, for Witt bases of the different O-types, the Ricci polynomial is
$$P_\Phi = \cases{\epsilon Y^2V^2 - b(X^2V^2 - Y^2U^2) + 4\lambda XYUV,&${\bf O^\bfplus_\bfplus}$;\cr \epsilon X^2U^2 - b(Y^2U^2 - X^2V^2) + 4\lambda XYUV,&${\bf O^\bfminus_\bfminus}$;\cr \epsilon Y^2V^2 - b(Y^2U^2 - X^2V^2) + 4\lambda XYUV,&${\bf O^\bfminus_\bfplus}$;\cr \epsilon X^2U^2 - b(X^2V^2 - Y^2U^2) + 4\lambda XYUV,&${\bf O^\bfplus_\bfminus}$,\cr}\eqno(6.22.6)$$
but only Witt bases of O-type ${\bf O^\bfplus_\bfplus}$ and either ${\bf O^\bfminus_\bfplus}$ or ${\bf O^\bfplus_\bfminus}$ are required to describe the distinct $\Phi^a{}_b$ of type VI due to the isotropy subgroup determined in (4.10). The Ricci spinor type is $(2,2)$.

The real null eigenvector is $E_1^a$. For bases of type ${\bf O^\bfplus_\bfplus}$, $E_1^a = \ell^a$, and one computes from (6.13) that
$$\chi =\lambda \hskip .75in \alpha = b \hskip .75in \beta = -b,\eqno(6.22.7)$$
whence
$$4\chi^2 - \alpha\beta = 4\lambda^2 + b^2 > 0.\eqno(6.22.8)$$
Hence, the null eigenvector is a node with real tangents. Explicitly, $\ell^a$ has homogeneous coordinates $([1,0],[1,0])$ so choose affine coordinates $y = Y/X$ and $v = V/U$. Then
$$f(y,v) = \epsilon y^2v^2 - b(v^2 - y^2) + 4\lambda yv,$$
and 
$$F(t) = \left[4\lambda\sigma\tau - b(\tau^2 - \sigma^2)\right]t^2 + \epsilon\sigma^2\tau^2t^4.$$
Thus, $(0,0)$ is a double point, the tangents at which are given by
$${\sigma \over \tau} = {-2\lambda \pm 2\sqrt{4\lambda^2 + b^2} \over b},$$
i.e., $(0,0)$ is indeed a node with real tangents.
\vskip 24pt
\noindent {\bf 6.23 Type VII}\hfil\break
For each generic traceless, self-adjoint endomorphism $\Phi^a{}_b$ of type VII, there is a Witt basis $\{E_1^a,E_2^a,E_3^a,E_4^a\}$ with respect to which $\Phi^a{}_b$ has matrix representation (4.5.8) with $\mu = -\lambda \not= 0$. One computes
$$\Phi^{ab} = \epsilon E_1^aE_1^b + \omega E_2^aE_2^b + 2\lambda(E_1^{(a}E_3^{b)} - E_2^{(a}E_4^{b)}).\eqno(6.23.1)$$
For Witt bases of type ${\bf O^\bfplus_\bfplus}$,
$$\Phi_{ABA'B'} = (\epsilon o_Ao_B + \omega\iota_A\iota_B)o_{A'}o_{B'} + \lambda(o_A\iota_B + \iota_Ao_B)(o_{A'}\iota_{B'} + \iota_{A'}o_{B'}),\eqno(6.23.2)$$
and the Ricci polynomial is
$$P_\Phi(X,Y,U,V) = \left[(\epsilon Y^2 + \omega X^2)V - 2\lambda UXY\right]V.\eqno(6.23.3)$$
Hence,
$$P_\Phi = \cases{\left[(\epsilon Y^2 + \omega X^2)V - 2\lambda UXY\right]V,&${\bf O^\bfplus_\bfplus}$;\cr \left[(\epsilon X^2 + \omega Y^2)U - 2\lambda VXY\right]U,&${\bf O^\bfminus_\bfminus}$;\cr \left[(\epsilon V^2 + \omega U^2)Y - 2\lambda UVX\right]Y,&${\bf O^\bfminus_\bfplus}$;\cr \left[(\epsilon U^2 + \omega V^2)X - 2\lambda UVY\right]X,&${\bf O^\bfplus_\bfminus}$.\cr}\eqno(6.23.4)$$
The isotropy subgroup is discrete and distinct forms are parametrised by bases of type ${\bf O^\bfplus_\bfplus}$ and ${\bf O^\bfplus_\bfminus}$, say. Hence, the Ricci spinor types are $(2,1)(0,1)$ and $(1,0)(1,2)$.

$\Phi^a{}_b$ has two null eigenvectors $E_1^a$ and $E_2^a$. Restricting attention to bases of type ${\bf O^\bfplus_\bfplus}$, the null eigenvectors take the form $\ell^a$, with eigenvalue $\lambda$, and $\tilde m^a$, with eigenvalue $-\lambda$. From (6.13), one computes
$$\chi = \lambda \qquad \alpha = 0 \qquad \beta = \omega/2$$
and
$$\chi = -\lambda \qquad \alpha = 0 \qquad \beta = \epsilon/2,$$
respectively. Hence, in both cases, $4\chi^2 - \alpha\beta = 4\lambda^2 > 0$, i.e., each null eigenvector determines a real node with real tangents. Explicitly, $\ell^a$ has homogeneous coordinates $([1,0],[1,0])$ so taking affine coordinates $y = Y/X$ and $v = V/U$ yields $f(y,v) = (\epsilon y^2 + \omega)v^2 - 2\lambda yv$. Hence, With $(y,v) = (\sigma t,\tau t)$, (6.3) gives
$$f(t) = (\omega\tau^2 - 2\lambda\sigma\tau)t^2 + \epsilon\sigma^2\tau^2t^4.$$
The tangents are therefore given by $\tau = 0$ and $\omega\tau = 2\lambda\sigma$. The tangent $\tau = 0$ is the component $V = 0$.

Similarly, $\tilde m^a$ has homogeneous coordinates $([0,1],[1,0])$ and taking affine coordinates $x = X/Y$ and $v = V/U$ yields $f(x,v) = (\epsilon + \omega x^2)v^2 - 2\lambda xv$. Hence, With $(x,v) = (\sigma t,\tau t)$, (6.3) gives
$$f(t) = (\epsilon\tau^2 - 2\lambda\sigma\tau)t^2 + \omega\sigma^2\tau^2t^4.$$
Thus, one sees that the linear component $V=0$ of the Ricci locus is a common tangent of the two nodes.

The coincidence  $\lambda = \mu$ forces $\lambda = \mu = 0$, whence the Ricci polynomial, for ${\bf O^\bfplus_\bfplus}$-type bases reduces to
$$P_\Phi = (\epsilon Y^2 + \omega X^2)V^2,$$
which has spinor type $(1,0)(1,0)(0,1)^2$ when $\epsilon = -\omega$ and type $(1,0)\overline{(1,0)}(0,1)^2$ when $\epsilon = \omega$. The two null eigenvectors $E_1^a$ and $E_2^a$ now have a common eigenvalue, whence $\langle E_1^a,E_2^a \rangle_{\bf K}$ is a totally null subspace of eigenvectors. For Witt bases of O-type ${\bf O^\bfplus_\bfplus}$, each null eigenvector is of the form $\alpha^Ao^{A'}$, for arbitrary $\alpha^A$, i.e., they have homogeneous coordinates $([a,b][1,0])$, i.e., constitute the component $V=0$, which being double is singular. When $\epsilon = -\omega$, the real Ricci locus is two `parallel' ${\bf S}^1$'s on the torus ${\bf S}^1 \times {\bf S}^1$, to which the double component $V=0$ is `orthogonal'. When $\epsilon=\omega$, there is an analogous description of the Ricci locus $\Omega$. Technically, the points of intersection of the double component with the other two components are tacnodes, with a common tangent ($V=0$), but the points of intersection are not $E_1^a$ and $E_2^a$.
\vskip 24pt
\noindent {\bf 6.24 Type VIII}\hfil\break
For each generic traceless, self-adjoint endomorphism $\Phi^a{}_b$ of type VIII, ther exists an $\Psi$-ON basis\hfil\break
 $\{U^a,V^a,X^a,Y^a\}$ with respect to which $\Phi^a{}_b$ has matrix representation (4.5.11) with $c=-a\not=0$. One computes that
$$\Phi^{ab} = a(U^aU^b - V^aV^b - X^aX^b + Y^aY^b) - 2bU^{(a}X^{b)} - 2dV^{(a}Y^{b)}.\eqno(6.24.1)$$
For $\Psi$-ON bases of O-type ${\bf O^\bfplus_\bfplus}$, one obtains
$$\Phi_{ABA'B'} = a(o_A\iota_B + \iota_Ao_B)(o_{A'}\iota_{B'} + \iota_{A'}o_{B'}) - b(o_Ao_Bo_{A'}o_{B'} - \iota_A\iota_B\iota_{A'}\iota_{B'}) - d(\iota_A\iota_Bo_{A'}o_{B'} - o_Ao_B\iota_{A'}\iota_{D'}),\eqno(6.24.2)$$
for which the Ricci polynomial is
$$P_\Phi(X,Y,U,V) = 4aUVXY + b(X^2U^2 - Y^2V^2) + d(Y^2U^2 - X^2V^2).$$
Hence,
$$P_\Phi = \cases{4aUVXY + b(X^2U^2 - Y^2V^2) + d(Y^2U^2 - X^2V^2),&${\bf O^\bfplus_\bfplus}$;\cr 4aUVXY + b(Y^2V^2 - X^2U^2) + d(X^2V^2 - Y^2U^2),&${\bf O^\bfminus_\bfminus}$;\cr 4aUVXY + b(X^2U^2 - Y^2V^2) + d(X^2V^2 - Y^2U^2),&${\bf O^\bfminus_\bfplus}$;\cr 4aUVXY + b(Y^2V^2 - X^2U^2) + d(Y^2U^2 - X^2V^2),&${\bf O^\bfplus_\bfminus}$.\cr}\eqno(4.24.3)$$
From (4.10), the isotropy subgroup is discrete and intersects ${\bf O^\bfminus_\bfminus(2,2)}$ so only bases of O-type ${\bf O^\bfplus_\bfplus}$ and either ${\bf O^\bfminus_\bfplus}$ or ${\bf O^\bfplus_\bfminus}$ are required to describe the distinct endomorphisms of type VIII. The Ricci spinor type is $(2,2)$.

$\Phi^a{}_b$ has two pairs of complex conjugate eigenvalues $(\lambda,\bar\lambda)$ and $(\mu,\bar\mu)$; $\lambda$ has eigenvalue $U^a + iX^a$ and $\mu$ has eigenavlue $V^a + iY^a$. These eigenvectors, and their conjugates, are null in $\overline{\bf C}^{2,2}$ but nonnull in ${\bf C}^{2,2}$. Hence, there are no singularities.

The coincidence $\lambda = \mu$ forces $c=d$ and $a=c=0$. For $\Psi$-ON bases of O-type ${\bf O^\bfplus_\bfplus}$, the Ricci polynomial reduces to
$$P_\Phi = b(X^2+Y^2)(U^2-V^2) = b(X+iY)(X-iY)(U-V)(U+V),\eqno(6.24.4)$$
i.e., type $(1,0)\overline{(1,0)}(0,1)(0,1)$, while for bases of O-type ${\bf O^\bfminus_\bfplus}$, the Ricci polynomial reduces to
$$P_\Phi = b(X^2-Y^2)(U^2+V^2) = b(X-Y)(X+Y)(U+iV)(U-iV),\eqno(6.24.5)$$
i.e., type $(1,0)(1,0)(0,1)\overline{(0,1)}$.

The eigenspace of $\lambda = \mu$ is $\langle U^a+iX^a,V^a+iY^a\rangle_{\bf C}$, which contains two linearly independent complex eigenvectors that are null in ${\bf C}^{2,2}$, viz., 
$$\displaylines{N_1 := (U^a + iX^a) + i(V^a+iY^a) = (U^a-Y^a) + i(X^a+V^a)\cr
\noalign{\vskip -6pt}
\hfill\llap(6.24.6)\cr
\noalign{\vskip -6pt}
N_2 := (U^a + iX^a) - i(V^a+iY^a) = (U^a+Y^a) + i(X^a-V^a).\cr}$$
Their complex conjugates are complex null eigenvectors of $\bar\lambda=\bar\mu$. From (6.13), one computes that for each complex null eigenvalue, $\alpha=\beta=0$, whence $4\chi^2 -\alpha\beta = -4b^2$ in each case and each is therefore a complex node (i.e., node on $\Omega\setminus \omega$). The four linear factors of $P_\Phi$ each define a line on $\Omega$; these four lines define four points of intersection between pairs of lines that are the four null eigenvectors, the pair of lines therefore providing the tangents at the node. For example, $N_1$ has homogeneous coordinates $([1+i,-1+i],[1,-1])$. With respect to $\Psi$-ON bases of O-type ${\bf O^\bfplus_\bfplus}$, (6.24.4) gives $P_\Phi$ and $N_1$ lies on the components given by $X+iY=0$ and $U+V=0$, i.e., is the intersection point of these two components. For an affine picture, and noting that $(-1+i)/(1+i) = i$, first rewrite the homogeneous coordinates as $([1,i],[1,-1])$, then take affine coordinates $y = Y/X$, $v = V/U$ and consider lines $(y = i + \sigma t,v = -1 + \tau t)$ through the point determined by $N_1$. One computes 
$$f(y,v) = b(4i\sigma\tau t^2 + 2\sigma\tau(\sigma-i\tau)t^3 - \sigma^2\tau^2t^4),$$
which shows that $N_1$ indeed determines a complex node and that the tangents, $\tau =0$ and $\sigma=0$, each determine a component of the curve; in this affine picture, $N_1$ is the intersection of $v=-1$ and $y = i$. One similarly computes that in the same affine coordinates $N_2$ is the intersection of $v=1$ and $y=i$. The other two complex nodes are the complex conjugates of these two.

Finally, note that the coincidence $\lambda = \bar\mu$ is not a distinct case but merely a relabelling of the complex conjugate pair of eigenvalues $(\mu,\bar\mu)$.
\vskip 24pt
\noindent {\bf 6.25 Type IX}\hfil\break
For a generic traceless, self-adjoint endomorphism $\Phi^a{}_b$ of type IX, there is a $\Psi$-ON basis $\{U^a,V^a,X^a,Y^a\}$ with respect to which $\Phi^a{}_b$ has matrix representation (4.10.16), with $\lambda+\mu+\nu+\sigma=0$, but no coincedences of eigenvalues. One therefore computes that
$$\Phi^{ab} = \lambda U^aU^b + \mu V^aV^b - \nu X^aX^b - \sigma Y^aY^b.\eqno(6.25.1)$$
From (4.10), the isotropy subgroup intersects all four components of {\bf O(2,2)} so it suffices to consider $\Psi$-ON bases of O-type ${\bf O^\bfplus_\bfplus}$. For such,
$$\eqalignno{\Phi_{ABA'B'} &= \left({\lambda-\nu \over 2}\right)(o_Ao_Bo_{A'}o_{B'} + \iota_A\iota_B\iota_{A'}\iota_{B'}) + \left({\mu-\sigma \over 2}\right)(\iota_A\iota_Bo_{A'}o_{B'} + o_Ao_B\iota_{A'}\iota_{B'})\cr
\noalign{\vskip 6pt}
&+ \left({\lambda+\nu \over 2}\right)(o_A\iota_B + \iota_Ao_B)(o_{A'}\iota_{B'} + \iota_{A'}o_{B'}),&(6.25.2)\cr}$$
and the Ricci polynomial is
$$P_\phi(X,Y,U,V) = \left({\lambda-\nu \over 2}\right)(Y^2V^2+X^2U^2) + \left({\mu-\sigma \over 2}\right)(X^2V^2+Y^2U^2) + 2(\lambda+\nu)UVXY,\eqno(6.25.3)$$
which is of type $(2,2)$, i.e., irreducible. There are no null eigenvectors, whence no singularities.

Each coincidence of two eigenvalues creates a two-dimensional eigenspace for that eigenvalue, isomorphic to ${\bf R}^{2,0}$, ${\bf R}^{0,2}$, or ${\bf R}^{1,1}$. In the last case there are two, linearly independent real null eigenvectors; in the cases of definite signature, there are two, linearly independent, complex eigenvectors that are null in ${\bf C}^{2,2}$. 

Consider first $\nu=\sigma$. The trace-free condition is
$$\nu=\sigma = -\left({\lambda+\mu \over 2}\right)\eqno(6.25.4)$$
while the conditions for genericity in this case are to exclude any other coincidences:
$$\nu=\sigma\not=\lambda (\Leftrightarrow 3\lambda+\mu\not=0) \qquad \nu=\sigma\not=\mu (\Leftrightarrow \lambda+3\mu\not=0) \qquad \lambda\not=\mu.\eqno(6.25.5)$$
Now, (6.25.2) becomes
$$\eqalignno{\Phi_{ABA'B'} &= \left({3\lambda+\mu \over 4}\right)(o_Ao_Bo_{A'}o_{B'} + \iota_A\iota_B\iota_{A'}\iota_{B'}) + \left({\lambda+3\mu \over 4}\right)(\iota_A\iota_Bo_{A'}o_{B'} + o_Ao_B\iota_{A'}\iota_{B'})\cr
\noalign{\vskip 6pt}
&+ \left({\lambda-\mu \over 4}\right)(o_A\iota_B + \iota_Ao_B)(o_{A'}\iota_{B'} + \iota_{A'}o_{B'}),&(6.25.6)\cr}$$
and Ricci polynomial is now
$$P_\Phi = \left({3\lambda+\mu \over 4}\right)(Y^2V^2+X^2U^2) + \left({\lambda+3\mu \over 4}\right)(X^2V^2+Y^2U^2) + (\lambda-\mu)UVXY.\eqno(6.25.7)$$
Defining
$$M_{\lambda,\mu} := {3\lambda+\mu \over 4} \hskip 1.25in N_{\lambda,\mu} := {\lambda + 3\mu \over 4},\eqno(6.25.8)$$
then (6.25.6) can be written
$$P_\Phi = M_{\lambda,\mu}(Y^2V^2+X^2U^2) + N_{\lambda,\mu}(X^2V^2+Y^2U^2) + 2(M_{\lambda,\mu} - N_{\lambda,\mu})UVXY,\eqno(6.25.9)$$
and the conditions in (6.25.5) are
$$M_{\lambda,\mu} \not= 0 \qquad N_{\lambda,\mu} \not= 0 \qquad M_{\lambda,\mu} \not= N_{\lambda,\mu}.\eqno(6.25.10)$$
Dropping the subscripts on $M_{\lambda,\mu}$ and $N_{\lambda,\mu}$, and with
$$s(M) := {M \over \vert M \vert} \hskip 1.25in s(N) := {N \over \vert N \vert},\eqno(6.25.11)$$
one finds that $P_\Phi$ factorizes as follows:
$$s(M)(\sqrt{\vert M \vert}UX \pm i\sqrt{\vert  N \vert}UY \mp i\sqrt{\vert N \vert}VX  + \sqrt{\vert M \vert}VY)(\sqrt{\vert M \vert}UX \mp i \sqrt{\vert N \vert}UY \pm i\sqrt{\vert N \vert}VX + \sqrt{\vert M \vert}VY),\eqno(6.25.12)$$
when $s(M) = s(N)$;
$$s(M)(\sqrt{\vert M \vert}XU \pm \sqrt{\vert  N \vert}UY \mp \sqrt{\vert N \vert}VX + \sqrt{\vert M \vert}VY)(\sqrt{\vert M \vert}UX \mp \sqrt{\vert N \vert}UY \pm \sqrt{\vert N \vert}VX + \sqrt{\vert M \vert}VY)\eqno(6.25.13)$$
when $s(M) = -s(N)$. I note that in (6.25.13) the two factors are not equal. Hence, the Ricci spinor type is $(1,1)\overline{(1,1)}$ and $(1,1)(1,1)$, respectively. Bearing in mind (6.25.10), in $(M,N)$-space, type $(1,1)\overline{(1,1)}$ occurs in the open first and third quadrants minus the line $M = N$; which is four connected components. Type $(1,1)(1,1)$ occurs in the open second and fourth quadrants. Hence, there are six connected components for this case, as noted in (4.10). These conditions are easily translated into $(\lambda,\mu)$-space if desired; the topology is unchanged, of course, only the shapes of the components change.

The coincidence $\lambda = \mu$ has a similar description. The tracefree condition is
$$\lambda = \mu = -\left({\nu+\sigma \over 2}\right),\eqno(6.25.14)$$
and the conditions of genericity are
$$3\nu+\sigma \not= 0 \qquad \nu+3\sigma \not= 0 \qquad \nu \not= \sigma.\eqno(6.25.15)$$
The Ricci spinor becomes
$$\eqalignno{\Phi_{ABA'B'} &= -\left({3\nu+\sigma \over 4}\right)(o_Ao_Bo_{A'}o_{B'} + \iota_A\iota_B\iota_{A'}\iota_{B'}) - \left({\nu+3\sigma \over 4}\right)(\iota_A\iota_Bo_{A'}o_{B'} + o_Ao_B\iota_{A'}\iota_{B'})\cr
\noalign{\vskip 6pt}
&+ \left({\nu-\sigma \over 4}\right)(o_A\iota_B + \iota_Ao_B)(o_{A'}\iota_{B'} + \iota_{A'}o_{B'}),&(6.25.16)\cr}$$
and the Ricci polynomial
$$P_\Phi = -M_{\nu,\sigma}(Y^2V^2+X^2U^2) - N_{\nu,\sigma}(X^2V^2+Y^2U^2) + 2(M_{\nu,\sigma} - N_{\nu,\sigma})UVXY.\eqno(6.25.17)$$
With the same notation as above, $P_\Phi$ factorizes as follows:
$$-s(M)(\sqrt{\vert M \vert}UX \pm i\sqrt{\vert  N \vert}UY \pm i\sqrt{\vert N \vert}VX  - \sqrt{\vert M \vert}VY)(\sqrt{\vert M \vert}UX \mp i \sqrt{\vert N \vert}UY \mp i\sqrt{\vert N \vert}VX - \sqrt{\vert M \vert}VY),\eqno(6.25.18)$$
when $s(M) = s(N)$;
$$-s(M)(\sqrt{\vert M \vert}XU \pm \sqrt{\vert  N \vert}UY \pm \sqrt{\vert N \vert}VX - \sqrt{\vert M \vert}VY)(\sqrt{\vert M \vert}UX \mp \sqrt{\vert N \vert}UY \mp \sqrt{\vert N \vert}VX - \sqrt{\vert M \vert}VY)\eqno(6.25.19)$$
when $s(M) = -s(N)$. The same comments apply as in the coincidence $\nu=\sigma$.

For the coincidence $\lambda=\sigma$, the null eigenvectors are real. The tracefree condition is
$$\lambda = \sigma = -\left({\mu+\nu \over 2}\right),\eqno(6.25.20)$$
and the conditions of genericity are
$$3\mu+\nu \not= 0 \qquad \mu+3\nu \not= 0 \qquad \mu \not= \nu.\eqno(6.25.21)$$
The Ricci spinor becomes
$$\eqalignno{\Phi_{ABA'B'} &= -\left({3\nu+\mu \over 4}\right)(o_Ao_Bo_{A'}o_{B'} + \iota_A\iota_B\iota_{A'}\iota_{B'}) + \left({\nu+3\mu \over 4}\right)(\iota_A\iota_Bo_{A'}o_{B'} + o_Ao_B\iota_{A'}\iota_{B'})\cr
\noalign{\vskip 6pt}
&+ \left({\nu-\mu \over 4}\right)(o_A\iota_B + \iota_Ao_B)(o_{A'}\iota_{B'} + \iota_{A'}o_{B'}),&(6.25.22)\cr}$$
and the Ricci polynomial
$$P_\Phi = -M_{\nu,\mu}(Y^2V^2+X^2U^2) + N_{\nu,\mu}(X^2V^2+Y^2U^2) + 2(M_{\nu,\mu} - N_{\nu,\mu})UVXY.\eqno(6.25.23)$$
The factorization of this $P_\Phi$ is as follows:
$$-s(M)(\sqrt{\vert M \vert}UX \pm i\sqrt{\vert  N \vert}UY \mp i\sqrt{\vert N \vert}VX  - \sqrt{\vert M \vert}VY)(\sqrt{\vert M \vert}UX \mp i \sqrt{\vert N \vert}UY \pm i\sqrt{\vert N \vert}VX - \sqrt{\vert M \vert}VY),\eqno(6.25.24)$$
when $s(M) = -s(N)$;
$$-s(M)(\sqrt{\vert M \vert}XU \pm \sqrt{\vert  N \vert}UY \mp \sqrt{\vert N \vert}VX - \sqrt{\vert M \vert}VY)(\sqrt{\vert M \vert}UX \mp \sqrt{\vert N \vert}UY \pm \sqrt{\vert N \vert}VX - \sqrt{\vert M \vert}VY)\eqno(6.25.25)$$
when $s(M) = s(N)$. The same comments apply as in the coincidence $\nu=\sigma$.

For the coincidence $\mu=\nu$ the null eigenvectors are real. The trace-free condition is
$$\nu=\mu = -\left({\lambda+\sigma \over 2}\right)\eqno(6.25.26)$$
while the conditions for genericity are:
$$3\lambda+\sigma\not=0 \qquad \lambda+3\sigma\not=0 \qquad \lambda\not=\sigma.\eqno(6.25.27)$$
The Ricci spinor is
$$\eqalignno{\Phi_{ABA'B'} &= \left({3\lambda+\sigma \over 4}\right)(o_Ao_Bo_{A'}o_{B'} + \iota_A\iota_B\iota_{A'}\iota_{B'}) - \left({\lambda+3\sigma \over 4}\right)(\iota_A\iota_Bo_{A'}o_{B'} + o_Ao_B\iota_{A'}\iota_{B'})\cr
\noalign{\vskip 6pt}
&+ \left({\lambda-\sigma \over 4}\right)(o_A\iota_B + \iota_Ao_B)(o_{A'}\iota_{B'} + \iota_{A'}o_{B'}),&(6.25.28)\cr}$$
and Ricci polynomial is 
$$P_\Phi = M_{\lambda,\sigma}(Y^2V^2+X^2U^2) - N_{\lambda,\sigma}(X^2V^2+Y^2U^2) + 2(M_{\lambda,\sigma}-N_{\lambda,\sigma})UVXY.\eqno(6.25.29)$$
$P_\Phi$ factorizes as follows:
$$s(M)(\sqrt{\vert M \vert}UX \pm i\sqrt{\vert  N \vert}UY \pm i\sqrt{\vert N \vert}VX  + \sqrt{\vert M \vert}VY)(\sqrt{\vert M \vert}UX \mp i \sqrt{\vert N \vert}UY \mp i\sqrt{\vert N \vert}VX + \sqrt{\vert M \vert}VY),\eqno(6.25.30)$$
when $s(M) = -s(N)$;
$$s(M)(\sqrt{\vert M \vert}XU \pm \sqrt{\vert  N \vert}UY \pm \sqrt{\vert N \vert}VX + \sqrt{\vert M \vert}VY)(\sqrt{\vert M \vert}UX \mp \sqrt{\vert N \vert}UY \mp \sqrt{\vert N \vert}VX + \sqrt{\vert M \vert}VY)\eqno(6.25.31)$$
when $s(M) = s(N)$. The same comments apply as in the coincidence $\nu=\sigma$.

For the coincidence $\lambda=\nu$, the null eigenvectors are real but have a simpler spinor description than in the previous cases, resulting in a simpler description of the Ricci spinor and polynomial. The tracefree condition is
$$\lambda = \nu = -\left({\mu+\sigma \over 2}\right),\eqno(6.25.32)$$
while the conditions for genericity are
$$3\mu+\sigma \not= 0 \qquad \mu+3\sigma \not= 0 \qquad \mu\not=\sigma.\eqno(6.25.33)$$
The Ricci spinor is
$$\Phi_{ABA'B'} = \left({\mu-\sigma \over 2}\right)(\iota_A\iota_Bo_{A'}o_{B'} + o_Ao_B\iota_{A'}\iota_{B'}) - \left({\mu+\sigma \over 2}\right)(o_A\iota_B + \iota_Ao_B)(o_{A'}\iota_{B'} + \iota_{A'}o_{B'}),\eqno(6.25.34)$$
and the Ricci polynomial is
$$P_\Phi = \left({\mu-\sigma \over 2}\right)(X^2V^2+Y^2U^2) - 2(\mu+\sigma)UVXY.\eqno(6.25.35)$$
Putting
$$k := {\mu+\sigma \over \mu-\sigma} \in {\bf R},\eqno(6.25.36)$$
let $r$ denote any root of the real quadratic
$$x^2 + 4kx + 1 = 0.\eqno(6.25.37)$$
The discriminant of this quadratic is
$$4k^2 - 1 = {3\mu^2 + 10\mu\sigma + 3\sigma^2 \over (\mu-\sigma)^2} = {(3\mu+\sigma)(\mu+3\sigma) \over (\mu-\sigma)^2}.\eqno(6.25.38)$$
Hence, $4k^2-1 = 0$ is excluded by (6.25.33). $P_\Phi$ factorizess as follows:
$$P_\Phi = \left({\mu-\sigma \over 2}\right)[rUY + VX]\left[{UY \over r} + VX\right],\eqno(6.25.39)$$
when $4k^2 - 1 > 0$ (whence $r$ is real), i.e., when $(3\mu+\sigma)(\mu+3\sigma) > 0$;
$$P_\Phi = \left({\mu-\sigma \over 2}\right)[rUY + VX][\bar r UY + VX],\eqno(6.25.40)$$
when $4k^2 - 1 < 0$ (whence $r$ is complex), i.e., when $(3\mu+\sigma)(\mu+3\sigma) < 0$.
These forms are type $(1,1)(1,1)$ and $(1,1)\overline{(1,1)}$, respectively, and the first cannot reduce to $(1,1)^2$ under the present assumptions.

Finally, The coincidence $\mu=\sigma$ is similar to the previous case. The tracefree condition is
$$\mu = \sigma = -\left({\lambda+\nu \over 2}\right),\eqno(6.25.41)$$
while the conditions for genericity are
$$3\lambda+\nu \not= 0 \qquad \lambda+3\nu \not= 0 \qquad \lambda\not=\nu.\eqno(6.25.42)$$
The Ricci spinor is
$$\Phi_{ABA'B'} = \left({\lambda-\nu \over 2}\right)(o_Ao_Bo_{A'}o_{B'} + \iota_A\iota_B\iota_{A'}\iota_{B'}) + \left({\lambda+\nu \over 2}\right)(o_A\iota_B + \iota_Ao_B)(o_{A'}\iota_{B'} + \iota_{A'}o_{B'}),\eqno(6.25.43)$$
and the Ricci polynomial is
$$P_\Phi = \left({\lambda-\nu \over 2}\right)(Y^2V^2+X^2U^2) + 2(\lambda+\nu)UVXY.\eqno(6.25.44)$$
With
$$k := {\lambda+\nu \over \lambda-\nu} \in {\bf R},\eqno(6.25.45)$$
let $r$ denote any root of the real quadratic
$$x^2 - 4kx + 1 = 0.\eqno(6.25.46)$$
The discriminant of this quadratic is
$$4k^2 - 1 = {3\lambda^2 + 10\lambda\nu + 3\nu^2 \over (\lambda-\nu)^2} = {(3\lambda+\nu)(\lambda+3\nu) \over (\lambda-\nu)^2}.\eqno(6.25.47)$$
Hence, $4k^2-1 = 0$ is excluded by (6.25.42). $P_\Phi$ factorizess as follows:
$$P_\Phi = \left({\lambda-\nu \over 2}\right)[rUX + VY]\left[{UX \over r} + VY\right],\eqno(6.25.48)$$
when $4k^2 - 1 > 0$ (whence $r$ is real), i.e., when $(3\lambda+\nu)(\lambda+3\nu) > 0$;
$$P_\Phi = \left({\lambda-\nu \over 2}\right)[rUX + VY][\bar r UX + VY],\eqno(6.25.49)$$
when $4k^2 - 1 < 0$ (whence $r$ is complex), i.e., when $(3\lambda+\nu)(\lambda+3\nu) < 0$.
These forms are type $(1,1)(1,1)$ and $(1,1)\overline{(1,1)}$, respectively, and the first cannot reduce to $(1,1)^2$ under the present assumptions.

The computations of (6.13) for each of these six coincidences may be described as follows. In each case, there are four real numbers $a$, $b$, $c$, and $d$ satisfying $a+b+c+d=0$, with $c=d$, say. Hence, 
$$c=d=-\left({a+b \over 2}\right).$$
That there are no further equalities amongst $a$, $b$, $c$ and $d$ is equivalent to
$$3a+b \not= 0 \qquad a+3b\not=0 \qquad a\not=b.$$
These facts are familiar from the preceding dicsussion. In each case, the two null eigenvectors have $c=d$ as eigenvalue, so
$$\chi = c = d = -\left({a+b \over 2}\right).$$
From (6.13), in each case and for each null eigenvector one finds
$$\alpha = \beta = \pm\left({a-b \over 2}\right),$$
whence
$$4\chi^2 - \alpha\beta = {3a^2 + 10ab + 3b^2 \over 4} = {(3a+b)(a+3b) \over 4}.\eqno(6.25.50)$$
Hence, in each case there is a pair of nodes, these nodes being the points of intersection of the two components of the Ricci locus. When the null eigenvectors are real (four cases), one can restrict attention to the real Ricci locus $\omega$, in which case both nodes have real tangents iff $(3a+b)(a+3b) > 0$ (and the Ricci spinor type is $(1,1)(1,1)$), while both are isolated points of $\omega$ iff $(3a+b)(a+3b) < 0$ (and the Ricci spinor type is $(1,1)\overline{(1,1)}$). In two cases, the nodes are both complex points, i.e., lie on $\Omega \setminus \omega$; in these two cases, $(3a+b)(a+3b) > 0$ is the condition for type $(1,1)\overline{(1,1)}$ while $(3a+b)(a+3b) < 0$ is the condition for type $(1,1)(1,1)$.

These six cases are not actually all distinct at a certain geometrical level, e.g., the cases $\lambda = \nu$ and $\lambda=\sigma$ are not distinct as regards the geometrical content presented above and may be converted into each other by relabelling. They do have distinct descriptions, however, reflecting the differences in the spinor descriptions of the two cases, resulting from the differing relationships between $\Psi$-ON bases of different O-types and null tetrads, which may be relevant in circmstances in which one does distinguish between $X^a$ and $Y^a$  and/or $U^a$ and $V^a$.

The affine picture in these cases adds little. To give but one example, when $\nu=\sigma$, $X^a+iY^a$ is a complex null eigenvector, with coordinates $([1,i],[1,i])$. Taking affine coordinates $y = Y/X$ and $v = V/U$, and considering lines $(y=i+\sigma t,v=i+ \tau t)$ through the point $(i,i)$ yields in place of (6.3)
$$F(t) = {1 \over 2}[(\mu-\lambda)(\sigma^2+\tau^2) - 4(\lambda+\mu)\sigma\tau]t^2 + \left({3\lambda + \mu \over 4}\right)[2i(\sigma\tau^2+\sigma^2\tau)t^3 + \sigma^2\tau^2t^4].$$
The discriminant of the quadratic coefficient of $t^2$ is indeed $3\lambda^2 + 10\lambda\mu + 3\mu^2$ and the roots do not make the coefficient of $t^3$ vanish.

Turning now to pairs of coincidences, when $\lambda=\mu\not=\nu=\sigma$, whence the tracefree condition imposes $\lambda+\nu=0$,
$$\Phi_{ABA'B'} = \lambda(o_Ao_B + \iota_A\iota_B)(o_{A'}o_{B'}+\iota_{A'}\iota_{B'}),\eqno(6.25.51)$$
and
$$P_\Phi = \lambda(X+iY)(X-iY)(U+iV)(U-iV),\eqno(6.25.52)$$
i.e., type $(1,0)\overline{(1,0)}(0,1)\overline{(0,1)}$. The Ricci locus has reduced to four distinct complex projective lines on ${\bf CP}^1 \times {\bf CP}^1$ that define four points of intersection, namely the four complex null eigenvectors of the separate cases of $\lambda=\mu$ and $\nu=\sigma$ from above; namely, $U^a \pm iV^a$ and $X^a \pm iY^a$, respectively. These four points are therefore complex nodes.

For $\lambda=\sigma\not=\mu=\nu$, the tracefree condition imposes $\lambda+\nu=0$, and
$$\Phi_{ABA'B'} = \lambda(o_Ao_B - \iota_A\iota_B)(o_{A'}o_{B'}-\iota_{A'}\iota_{B'}),\eqno(6.25.53)$$
whence
$$P_\Phi = \lambda(X-Y)(X+Y)(U-V)(U+V),\eqno(6.25.54)$$
i.e., type $(1,0)(1,0)(0,1)(0,1)$. The Ricci locus has reduced to four distinct projective lines. In this case, one can restrict to the real Ricci lcous $\omega$ in ${\bf RP}^1 \times {\bf RP}^1$. It consists of four distinct real projective lines, i.e., ${\bf S}^1$'s, on ${\bf S}^1 \times {\bf S}^1$ that define four points of intersection, namely the four real null eigenvectors from the separate cases $\lambda=\sigma$ and $\mu=\nu$ above. These null eigenvectors are thus real nodes.

The case $\lambda=\nu\not=\mu=\sigma$ is geometrically equivalent to the previous case at a certain level. One finds
$$\Phi_{ABA'B'} =\lambda(o_A\iota_A+\iota_Ao_B)(o_{A'}\iota_{B'}+\iota_{A'}o_{B'}),\eqno(6.25.55)$$
with $$P_\Phi = 4\lambda UVXY.\eqno(6.25.56)$$
Hence, the real Ricci locus is four projective lines that define four points of intersection, the null eigenvectors of the separate cases $\lambda=\nu$ and $\mu=\sigma$ above, which are real nodes.

I note that in each case, from (6.13) one computes that $\alpha=\beta=0$, whence $4\chi^2 -\alpha\beta = 4(\pm\lambda)^2 > 0$, in accord with the preceding geometric descriptions.

Now consider triple coincidences of eigenvalues. For $\lambda=\mu=\nu\not=\sigma$,
$$\Phi_{ABA'B'} = 2\lambda(\iota_A\iota_Ao_{A'}o_{B'} + o_Ao_B\iota_{A'}\iota_{B'}) + \lambda(o_A\iota_B+\iota_Ao_B)(o_{A'}\iota_{B'}+\iota_{A'}o_{B'}),\eqno(6.25.57)$$
and
$$P_\Phi = 2\lambda(XV+YU)^2.\eqno(6.25.58)$$
Hence, the type is $(1,1)^2$. The three-dimensional real eigenspace is $\langle U^a,V^a,X^a \rangle_{\bf R} \cong {\bf R}^{2,1}$. The null cone within this eigenspace provides an ${\bf S}^1$ of null eigenvectors that is precisley the projective line $XV+YU=0$ on ${\bf RP}^1 \times {\bf RP}^1$, i.e., the real Ricci locus $\omega$ is the double line of these null eigenvectors.

The case $\lambda=\mu=\sigma\not=\nu$ is geometrically similar to the previous case. Here
$$\Phi_{ABA'B'} = 2\lambda(o_Ao_Bo_{A'}o_{B'} + \iota_A\iota_B\iota_{A'}\iota_{B'}) - \lambda(o_A\iota_B+\iota_Ao_B)(o_{A'}\iota_{B'}+\iota_{A'}o_{B'}),\eqno(6.25.59)$$
with
$$P_\Phi = 2\lambda(YV-UX)^2.\eqno(6.25.60)$$
The real projective line $YV-UX=0$ describes the set of generators of the null cone of $\langle U^a,V^a,Y^a \rangle_{\bf R} \cong {\bf R}^{2,1}$.

The case $\lambda=\nu=\sigma\not=\mu$ yields
$$\Phi_{ABA'B'} = -2\lambda(\iota_A\iota_Ao_{A'}o_{B'} + o_Ao_B\iota_{A'}\iota_{B'}) + \lambda(o_A\iota_B+\iota_Ao_B)(o_{A'}\iota_{B'}+\iota_{A'}o_{B'}),\eqno(6.25.61)$$
with
$$P_\Phi = -2\lambda(XV-YU)^2.\eqno(6.25.62)$$
The real projective line $XV-YU=0$ describes the set of generators of the null cone in $\langle U^a,X^a,Y^a \rangle_{\bf R} \cong {\bf R}^{1,2}$.

The case $\mu=\nu=\sigma\not=\lambda$ yields
$$\Phi_{ABA'B'} = -2\nu(o_Ao_Bo_{A'}o_{B'} + \iota_A\iota_B\iota_{A'}\iota_{B'}) - \nu(o_A\iota_B+\iota_Ao_B)(o_{A'}\iota_{B'}+\iota_{A'}o_{B'}),\eqno(6.25.63)$$
with
$$P_\Phi = -2\nu(YV+XU)^2.\eqno(6.25.64)$$
The real projective line $YV+XU=0$ describes the set of generators of the null cone of $\langle V^a,X^a,Y^a \rangle_{\bf R} \cong {\bf R}^{1,2}$.

Finally, the coincidence $\lambda=\mu=\nu=\sigma$ together with the tracefree condition imposes $\lambda=\mu=\nu=\sigma=0$, whence $\Phi_{ABA'B'} = {\bf 0}$.

Table Two, beginning on the next page, summarizes the results of the paper. Recall that the topological structure of each subtype consists of orbits, with a fixed topological structure, trivially fibred over open domains in some ${\bf R}^N$, the latter providing the domain of values for the free parameters in the matrix form $\underline M$ characterizing the subtype. In the following table, df $= m+n$, where $n$ is the dimension of the orbit and $m$ the dimension of the space over which the orbits are fibred (i.e., the number of free parameters in $\underline M$); \# cc's is the number of connected components of the subtype, given as $p \times q$ where $q$ is the number of components of an orbit and $p$ that of the space over which the orbits are fibred.
\vfill\eject

$$\vbox{\offinterlineskip
\halign{&\strut\ #\ \cr
\multispan{8}\hfil\bf Table 2\hfil\cr
\noalign{\medskip}
\noalign{\hrule}
\noalign{\vskip 2pt}
\hfil${\scriptstyle{\bf Type}}$\hfil&\hfil${\scriptstyle{\bf JCF/coincidence\ of\ eigenvalues}}$\hfil&\hfil${\scriptstyle{\bf Traceless}}$\hfil&\hfil${\scriptstyle{\bf df}}$\hfil&\hfil${\scriptstyle{\bf spinor}}$\hfil&\hfil${\scriptstyle{\bf \#}}$\hfil&\hfil${\scriptstyle{\bf singularity}}$\hfil\cr
{}&\hfil${\scriptstyle{\bf Segre;\ genericity\ under\ tracelessness}}$\hfil&\hfil${\scriptstyle{\bf condition}}$\hfil&{}&\hfil${\scriptstyle{\bf type}}$\hfil&\hfil${\scriptstyle{\bf cc's}}$\hfil&\hfil${\scriptstyle {\bf structure}}$\hfil\cr
\noalign{\vskip 2pt}
\noalign{\hrule}
\noalign{\vskip 2pt}
\hfil ${\scriptstyle {\rm I}}$\hfil&\hfil${\scriptstyle \bigl({\bf R}^4_{\rm hb},J_4(\lambda)\bigr)}$\hfil&\hfil${\scriptstyle \lambda=0}$\hfil&\hfil${\scriptstyle 0+6}$\hfil&\hfil${\scriptstyle (2,1)(0,1)}$\hfil&\hfil${\scriptstyle 1 \times 2}$\hfil&\hfil${\scriptstyle {\rm tacnode\ at\ real}}$\hfil&\cr
{}&\hfil${\scriptstyle [4]}$\hfil&{}&{}&\hfil${\scriptstyle (1,0)(1,2)}$\hfil&\hfil${\scriptstyle 1 \times 2}$\hfil&\hfil${\scriptstyle {\rm null\ eigenvector}^1}$\hfil\cr
\noalign{\vskip 2pt}
\noalign{\hrule}
\noalign{\vskip 2pt}
\hfil ${\scriptstyle {\rm II}}$\hfil&\hfil${\scriptstyle \bigl({\bf R}^4_{\rm hb},K_2(a,b)\bigr)}$\hfil&\hfil${\scriptstyle a=0}$\hfil&\hfil${\scriptstyle 1+6}$\hfil&\hfil${\scriptstyle (2,1)(0,1)}$\hfil&\hfil${\scriptstyle 2 \times 2}$\hfil&\hfil${\scriptstyle {\rm pair\ of\ complex\ nodes}}$\hfil&\cr
{}&\hfil${\scriptstyle \ b \not= 0;\ [2\bar2]}$\hfil&{}&{}&\hfil${\scriptstyle (1,0)(1,2)}$\hfil&\hfil${\scriptstyle 2 \times 2}$\hfil&\hfil${\scriptstyle {\rm with\ common\ tangent}^2}$\hfil&\cr
\noalign{\vskip 2pt}
\noalign{\hrule}
\noalign{\vskip 2pt}
\hfil ${\scriptstyle {\rm IIIa}}$\hfil&\hfil${\scriptstyle \bigl({\bf R}^{1,0},J_1(\lambda)\bigr) \scoperp\, \bigl({\bf R}^{1,2},J_3(\mu)\bigr)}$\hfil&\hfil${\scriptstyle \lambda+3\mu=0}$\hfil&\hfil${\scriptstyle 1+6}$\hfil&\hfil${\scriptstyle (2,2)}$\hfil&\hfil${\scriptstyle 2 \times 2}$\hfil&\hfil${\scriptstyle {\rm cusp\ at\ real}}$\hfil&\cr
{}&\hfil${\scriptstyle [13];\ \lambda \not= 0}$\hfil&{}&{}&{}&{}&\hfil${\scriptstyle {\rm null\ eigenvector}^3}$\hfil&\cr
\noalign{\vskip 2pt}
\noalign{\hrule}
\noalign{\vskip 2pt}
\hfil ${\scriptstyle {\rm IIIa}}$\hfil&\hfil${\scriptstyle \lambda=\mu=0}$\hfil&\hfil${\scriptstyle \lambda=\mu=0}$\hfil&\hfil${\scriptstyle 0+5}$\hfil&\hfil${\scriptstyle (1,0)(1,1)(0,1)}$\hfil&\hfil${\scriptstyle 1 \times 2}$\hfil&\hfil${\scriptstyle {\rm triple\ point\ at\ real}}$\hfil&\cr
{}&\hfil${\scriptstyle [(13)]}$\hfil&{}&{}&{}&{}&\hfil${\scriptstyle {\rm null\ eigenvector}^4}$\hfil&\cr
\noalign{\vskip 2pt}
\noalign{\hrule}
\noalign{\vskip 2pt}
\hfil ${\scriptstyle {\rm IIIb}}$\hfil&\hfil${\scriptstyle \bigl({\bf R}^{2,1},J_3(\lambda)\bigr) \scoperp\, \bigl({\bf R}^{1,0},J_1(\mu)\bigr)}$\hfil&\hfil${\scriptstyle 3\lambda+\mu=0}$\hfil&\hfil${\scriptstyle 1+6}$\hfil&\hfil${\scriptstyle (2,2)}$\hfil&\hfil${\scriptstyle 2 \times 2}$\hfil&\hfil${\scriptstyle {\rm cusp\ at\ real}}$\hfil&\cr
{}&\hfil${\scriptstyle [31];\ \mu \not= 0}$\hfil&{}&{}&{}&{}&\hfil${\scriptstyle {\rm null\ eigenvector}^3}$\hfil&\cr
\noalign{\vskip 2pt}
\noalign{\hrule}
\noalign{\vskip 2pt}
\hfil ${\scriptstyle {\rm IIIb}}$\hfil&\hfil${\scriptstyle \lambda=\mu=0}$\hfil&\hfil${\scriptstyle \lambda=\mu=0}$\hfil&\hfil${\scriptstyle 0+5}$\hfil&\hfil${\scriptstyle (1,0)(1,1)(0,1)}$\hfil&\hfil${\scriptstyle 1 \times 2}$\hfil&\hfil${\scriptstyle {\rm triple\ point\ at\ real}}$\hfil&\cr
{}&\hfil${\scriptstyle [(31)]}$\hfil&{}&{}&{}&{}&\hfil${\scriptstyle {\rm null\ eigenvector}^4}$\hfil&\cr
\noalign{\vskip 2pt}
\noalign{\hrule}
\noalign{\vskip 2pt}
\hfil ${\scriptstyle {\rm IV}}$\hfil&\hfil${\scriptstyle \bigl({\bf R}^{1,0},J_1(\lambda)\bigr) \scoperp\, \bigl({\bf R}^2_{\rm hb},J_{\pm2}(\mu)\bigr) \scoperp\, \bigl({\bf R}^{0,1},J_1(\nu)\bigr)}$\hfil&\hfil${\scriptstyle \lambda+\nu=-2\mu}$\hfil&\hfil${\scriptstyle 2+6}$\hfil&\hfil${\scriptstyle (2,2)}$\hfil&\hfil${\scriptstyle 6 \times 1}$\hfil&\hfil${\scriptstyle {\rm node\ at}}$\hfil&\cr
{}&\hfil${\scriptstyle [121];\ 3\lambda + \nu \not= 0;\ \lambda + 3\nu \not= 0;\ \lambda \not = \nu}$\hfil&{}&{}&{}&{}&\hfil${\scriptstyle {\rm real\ null\ eigenvector}^5}$\hfil&\cr
\noalign{\vskip 2pt}
\noalign{\hrule}
\noalign{\vskip 2pt}
\hfil ${\scriptstyle {\rm IV}}$\hfil&\hfil${\scriptstyle \mu=\nu}$\hfil&\hfil${\scriptstyle \lambda+3\nu=0}$\hfil&\hfil${\scriptstyle 1+5}$\hfil&\hfil${\scriptstyle (1,1)(1,1)}$\hfil&\hfil${\scriptstyle 2 \times 1}$\hfil&\hfil${\scriptstyle {\rm tacnode\ at}}$\hfil&\cr
{}&\hfil${\scriptstyle [1(21)];\ 3\lambda + \nu \not= 0;\ \lambda \not = \nu}$\hfil&{}&{}&\hfil${\scriptstyle (1,1)\overline{(1,1)}}$\hfil&\hfil${\scriptstyle 2 \times 1}$\hfil&\hfil${\scriptstyle {\rm real\ null\ eigenvector}^6}$\hfil&\cr
\noalign{\vskip 2pt}
\noalign{\hrule}
\noalign{\vskip 2pt}
\hfil ${\scriptstyle {\rm IV}}$\hfil&\hfil${\scriptstyle \lambda=\mu}$\hfil&\hfil${\scriptstyle 3\lambda+\nu=0}$\hfil&\hfil${\scriptstyle 1+5}$\hfil&\hfil${\scriptstyle (1,1)\overline{(1,1)}}$\hfil&\hfil${\scriptstyle 2 \times 1}$\hfil&\hfil${\scriptstyle {\rm tacnode\ at}}$\hfil&\cr
{}&\hfil${\scriptstyle [(12)1];\ \lambda + 3\nu \not= 0;\ \lambda \not = \nu}$\hfil&{}&{}&\hfil${\scriptstyle (1,1)(1,1)}$\hfil&\hfil${\scriptstyle 2 \times 1}$\hfil&\hfil${\scriptstyle {\rm real\ null\ eigenvector}^7}$\hfil&\cr
\noalign{\vskip 2pt}
\noalign{\hrule}
\noalign{\vskip 2pt}
\hfil ${\scriptstyle {\rm IV}}$\hfil&\hfil${\scriptstyle \lambda=\nu}$\hfil&\hfil${\scriptstyle \lambda+\mu=0}$\hfil&\hfil${\scriptstyle 1+5}$\hfil&\hfil${\scriptstyle (1,0)(1,1)(0,1)}$\hfil&\hfil${\scriptstyle 2 \times 1}$\hfil&\hfil${\scriptstyle {\rm three\ real\ nodes}}$\hfil&\cr
{}&\hfil${\scriptstyle [(1\vert2\vert1)];\ \lambda + 3\nu \not= 0;\ 3\lambda + \nu\not= 0;\ {\rm i.e.,}\ \lambda \not= 0}$\hfil&{}&{}&{}&{}&\hfil${\scriptstyle {\rm with\ real\ tangents}^8}$\hfil&\cr
\noalign{\vskip 2pt}
\noalign{\hrule}
\noalign{\vskip 2pt}
\hfil ${\scriptstyle {\rm IV}}$\hfil&\hfil${\scriptstyle \lambda=\mu=\nu}$\hfil&\hfil${\scriptstyle \lambda=\mu=\nu=0}$\hfil&\hfil${\scriptstyle 0+3}$\hfil&\hfil${\scriptstyle (1,0)^2(0,1)^2}$\hfil&\hfil${\scriptstyle 1 \times 1}$\hfil&\hfil${\scriptstyle {\rm Ricci\ locus\ is\ a}}$\hfil&\cr
{}&\hfil${\scriptstyle [(121)]}$\hfil&{}&{}&{}&{}&\hfil${\scriptstyle {\rm pair\ of\ double\ lines}^9}$\hfil&\cr
\noalign{\vskip 2pt}
\noalign{\hrule}
\noalign{\vskip 2pt}
\hfil ${\scriptstyle {\rm V}}$\hfil&\hfil${\scriptstyle \bigl({\bf R}^{1,0},J_1(\lambda)\bigr) \scoperp\, \bigl({\bf R}^{1,1},K_1(a,b)\bigr) \scoperp\, \bigl({\bf R}^{0,1},J_1(\nu)\bigr)}$\hfil&\hfil${\scriptstyle \lambda+\nu=-2a}$\hfil&\hfil${\scriptstyle 3+6}$\hfil&\hfil${\scriptstyle (2,2)}$\hfil&\hfil${\scriptstyle 4 \times 1}$\hfil&\hfil${\scriptstyle{\rm no\ singularities}^{10}}$\hfil&\cr
{}&\hfil${\scriptstyle b \not= 0;\ [11\bar 11];\ \lambda \not= \nu}$\hfil&{}&{}&{}&{}&{}&\cr
\noalign{\vskip 2pt}
\noalign{\hrule}
\noalign{\vskip 2pt}
\hfil ${\scriptstyle {\rm V}}$\hfil&\hfil${\scriptstyle \lambda=\nu}$\hfil&\hfil${\scriptstyle \lambda=\nu=-a}$\hfil&\hfil${\scriptstyle 2+5}$\hfil&\hfil${\scriptstyle (1,1)(1,1)}$\hfil&\hfil${\scriptstyle 2 \times 1}$\hfil&\hfil${\scriptstyle {\rm two\ real\ nodes}}$\hfil&\cr
{}&\hfil${\scriptstyle b \not= 0;\ [(1\vert1\bar 1\vert1)]}$\hfil&{}&{}&{}&{}&\hfil${\scriptstyle {\rm with\ real\ tangents}^{11}}$\hfil&\cr
\noalign{\vskip 2pt}
\noalign{\hrule}
\noalign{\vskip 2pt}
\hfil ${\scriptstyle {\rm VI}}$\hfil&\hfil${\scriptstyle \bigl({\bf R}^2_{\rm hb},J_{\pm2}(\lambda)\bigr) \scoperp\, \bigl({\bf R}^{1,1},K_1(a,b)\bigr)}$\hfil&\hfil${\scriptstyle \lambda+a=0}$\hfil&\hfil${\scriptstyle 2+6}$\hfil&\hfil${\scriptstyle (2,2)}$\hfil&\hfil${\scriptstyle 2 \times 2}$\hfil&\hfil${\scriptstyle {\rm a\ real\ node}}$\hfil&\cr
{}&\hfil${\scriptstyle b \not= 0;\ [21\bar 1]}$\hfil&{}&{}&{}&{}&\hfil${\scriptstyle{\rm with\ real\ tangents}^{12}}$\hfil&\cr
\noalign{\vskip 2pt}
\noalign{\hrule}
\noalign{\vskip 2pt}
\hfil ${\scriptstyle {\rm VII}}$\hfil&\hfil${\scriptstyle \bigl({\bf R}^2_{\rm hb},J_{\pm2}(\lambda)\bigr) \scoperp\, \bigl({\bf R}^2_{\rm hb},J_{\pm2}(\mu)\bigr)}$\hfil&\hfil${\scriptstyle \lambda+\mu=0}$\hfil&\hfil${\scriptstyle 1+6}$\hfil&\hfil${\scriptstyle (2,1)(0,1)}$\hfil&\hfil${\scriptstyle 2 \times 1}$\hfil&\hfil${\scriptstyle {\rm two\ real\ nodes\  with}}$\hfil&\cr
{}&\hfil${\scriptstyle [22];\ \lambda \not= 0;\ \mu \not= 0}$\hfil&{}&{}&\hfil${\scriptstyle (1,0)(1,2)}$\hfil&\hfil${\scriptstyle 2 \times 1}$\hfil&\hfil${\scriptstyle {\rm common\ tangent}^{13}}$\hfil&\cr
\noalign{\vskip 2pt}
\noalign{\hrule}
\noalign{\vskip 2pt}
\hfil ${\scriptstyle {\rm VII}}$\hfil&\hfil${\scriptstyle \lambda=\mu}$\hfil&\hfil${\scriptstyle \lambda=\mu=0}$\hfil&\hfil${\scriptstyle 0+4}$\hfil&\hfil${\scriptstyle (1,0)(1,0)(0,1)^2}$\hfil&\hfil${\scriptstyle 1 \times 1}$\hfil&\hfil${\scriptstyle {\rm double\ line}}$\hfil&\cr
{}&\hfil${\scriptstyle [(22)]}$\hfil&{}&{}&\hfil${\scriptstyle (1,0)\overline{(1,0)}(0,1)^2}$\hfil&\hfil${\scriptstyle 1 \times 1}$\hfil&\hfil${\scriptstyle {\rm intersecting}}$\hfil&\cr
{}&{}&{}&{}&\hfil${\scriptstyle (1,0)^2(0,1)(0,1)}$\hfil&\hfil${\scriptstyle 1 \times 1}$\hfil&\hfil${\scriptstyle {\rm two\ other}}$\hfil&\cr
{}&{}&{}&{}&\hfil${\scriptstyle (1,0)^2(0,1)\overline{(0,1)}}$\hfil&\hfil${\scriptstyle 1 \times 1}$\hfil&\hfil${\scriptstyle {\rm components}^{14}}$\hfil&\cr
\noalign{\vskip 2pt}
\noalign{\hrule}
\noalign{\vskip 2pt}
\hfil ${\scriptstyle {\rm VIII}}$\hfil&\hfil${\scriptstyle \bigl({\bf R}^{1,1},K_1(a,b)\bigr) \scoperp\, \bigl({\bf R}^{1,1},K_1(c,d)\bigr)}$\hfil&\hfil${\scriptstyle a+c=0}$\hfil&\hfil${\scriptstyle 3+6}$\hfil&\hfil${\scriptstyle (2,2)}$\hfil&\hfil${\scriptstyle 8 \times 2}$\hfil&\hfil${\scriptstyle {\rm no\ singularities}^{15}}$\hfil&\cr
{}&\hfil${\scriptstyle b \not= 0;\ d \not= 0;\ [1\bar 11\bar 1];\ -c = a \not= 0}$\hfil&{}&{}&{}&{}&{}&\cr
\noalign{\vskip 2pt}
\noalign{\hrule}
\noalign{\vskip 2pt}
\hfil ${\scriptstyle {\rm VIII}}$\hfil&\hfil${\scriptstyle a+ib = c+id}$\hfil&\hfil${\scriptstyle a=c=0}$\hfil&\hfil${\scriptstyle 1+4}$\hfil&\hfil${\scriptstyle (1,0)\overline{(1,0)}(0,1)(0,1)}$\hfil&\hfil${\scriptstyle 2 \times 1}$\hfil&\hfil${\scriptstyle {\rm 4\ complex}}$\hfil&\cr
{}&\hfil${\scriptstyle b=d \not= 0;\ [(1\bar 11\bar 1)]}$\hfil&{}&{}&\hfil${\scriptstyle (1,0)(1,0)(0,1)\overline{(0,1)}}$\hfil&\hfil${\scriptstyle 2 \times 1}$\hfil&\hfil${\scriptstyle {\rm nodes}^{16}}$\hfil&\cr
\noalign{\vskip 2pt}
\noalign{\hrule}
\cr}}$$

$$\vbox{\offinterlineskip
\halign{&\strut\ #\ \cr
\multispan{8}\hfil\bf Table 2 Continued\hfil\cr
\noalign{\medskip}
\noalign{\hrule}
\noalign{\vskip 2pt}
\hfil${\scriptstyle {\bf Type}}$\hfil&\hfil${\scriptstyle {\bf JCF/coincidence\ of\ eigenvalues}}$\hfil&\hfil${\scriptstyle {\bf Traceless}}$\hfil&\hfil${\scriptstyle {\bf df}}$\hfil&\hfil${\scriptstyle {\bf spinor}}$\hfil&\hfil${\scriptstyle {\bf \#}}$\hfil&\hfil${\scriptstyle {\bf singularity}}$\hfil\cr
{}&\hfil${\scriptstyle {\bf Segre;\ genericity\ under\ tracelessness}}$\hfil&\hfil${\scriptstyle {\bf condition}}$\hfil&{}&\hfil${\scriptstyle {\bf type}}$\hfil&\hfil${\scriptstyle {\bf cc's}}$\hfil&\hfil${\scriptstyle {\bf structure}}$\hfil\cr
\noalign{\vskip 2pt}
\noalign{\hrule}
\noalign{\vskip 2pt}
\hfil ${\scriptstyle {\rm IX}}$\hfil&\hfil${\scriptstyle \bigl({\bf R}^{1,0},J_1(\lambda)\bigr) \scoperp\, \bigl({\bf R}^{1,0},J_1(\mu)\bigr) \scoperp\, \bigl({\bf R}^{0,1},J_1(\nu)\bigr) \scoperp\, \bigl({\bf R}^{0,1},J_1(\sigma)\bigr)}$\hfil&\hfil${\scriptstyle \lambda+\mu+\nu+\sigma=0}$\hfil&\hfil${\scriptstyle 3+6}$\hfil&\hfil${\scriptstyle (2,2)}$\hfil&\hfil${\scriptstyle 24 \times 1}$\hfil&\hfil${\scriptstyle {\rm no\ singularities}^{17}}$\hfil&\cr
{}&\hfil${\scriptstyle [1111]\ {\rm no\ two\ eigenvalues\ equal}}$\hfil&{}&{}&{}&{}&{}\cr
\noalign{\vskip 2pt}
\noalign{\hrule}
\noalign{\vskip 2pt}
\hfil ${\scriptstyle {\rm IX}}$\hfil&\hfil${\scriptstyle \nu=\sigma}$\hfil&\hfil${\scriptstyle \lambda+\mu+\nu+\sigma=0}$\hfil&\hfil${\scriptstyle 2+5}$\hfil&\hfil${\scriptstyle (1,1)\overline{(1,1)}}$\hfil&\hfil${\scriptstyle 4 \times 1}$\hfil&\hfil${\scriptstyle {\rm two\ complex}}$\hfil&\cr
{}&\hfil${\scriptstyle [11(11)];\ 3\lambda+\mu\not=0;\ \lambda+3\mu\not=0;\ \lambda\not=\mu}$\hfil&{}&{}&\hfil${\scriptstyle (1,1)(1,1)}$\hfil&\hfil${\scriptstyle 2 \times 1}$\hfil&\hfil${\scriptstyle {\rm nodes}^{18}}$\hfil\cr
\noalign{\vskip 2pt}
\noalign{\hrule}
\noalign{\vskip 2pt}
\hfil ${\scriptstyle {\rm IX}}$\hfil&\hfil${\scriptstyle \lambda=\mu}$\hfil&\hfil${\scriptstyle \lambda+\mu+\nu+\sigma=0}$\hfil&\hfil${\scriptstyle 2+5}$\hfil&\hfil${\scriptstyle (1,1)\overline{(1,1)}}$\hfil&\hfil${\scriptstyle 4 \times 1}$\hfil&\hfil${\scriptstyle {\rm two\ complex}}$\hfil&\cr
{}&\hfil${\scriptstyle [(11)11];\ 3\nu+\sigma\not=0;\ \nu+3\sigma\not=0;\ \nu\not=\sigma}$\hfil&{}&{}&\hfil${\scriptstyle (1,1)(1,1)}$\hfil&\hfil${\scriptstyle 2 \times 1}$\hfil&\hfil${\scriptstyle {\rm nodes}^{19}}$\hfil\cr
\noalign{\vskip 2pt}
\noalign{\hrule}
\noalign{\vskip 2pt}
\hfil ${\scriptstyle {\rm IX}}$\hfil&\hfil${\scriptstyle \lambda=\sigma}$\hfil&\hfil${\scriptstyle \lambda+\mu+\nu+\sigma=0}$\hfil&\hfil${\scriptstyle 2+5}$\hfil&\hfil${\scriptstyle (1,1)\overline{(1,1)}}$\hfil&\hfil${\scriptstyle 2 \times 1}$\hfil&\hfil${\scriptstyle {\rm two\ real}}$\hfil&\cr
{}&\hfil${\scriptstyle [(1\vert11\vert1)];\ 3\mu+\nu\not=0;\ \mu+3\nu\not=0;\ \mu\not=\nu}$\hfil&{}&{}&\hfil${\scriptstyle (1,1)(1,1)}$\hfil&\hfil${\scriptstyle 4 \times 1}$\hfil&\hfil${\scriptstyle {\rm nodes}^{19}}$\hfil\cr
\noalign{\vskip 2pt}
\noalign{\hrule}
\noalign{\vskip 2pt}
\hfil ${\scriptstyle {\rm IX}}$\hfil&\hfil${\scriptstyle \mu=\nu}$\hfil&\hfil${\scriptstyle \lambda+\mu+\nu+\sigma=0}$\hfil&\hfil${\scriptstyle 2+5}$\hfil&\hfil${\scriptstyle (1,1)\overline{(1,1)}}$\hfil&\hfil${\scriptstyle 2 \times 1}$\hfil&\hfil${\scriptstyle {\rm two\ real}}$\hfil&\cr
{}&\hfil${\scriptstyle [1(11)1];\ 3\lambda+\sigma\not=0;\ \lambda+3\sigma\not=0;\ \lambda\not=\sigma}$\hfil&{}&{}&\hfil${\scriptstyle (1,1)(1,1)}$\hfil&\hfil${\scriptstyle 4 \times 1}$\hfil&\hfil${\scriptstyle {\rm nodes}^{19}}$\hfil\cr
\noalign{\vskip 2pt}
\noalign{\hrule}
\noalign{\vskip 2pt}
\hfil ${\scriptstyle {\rm IX}}$\hfil&\hfil${\scriptstyle \lambda=\nu}$\hfil&\hfil${\scriptstyle \lambda+\mu+\nu+\sigma=0}$\hfil&\hfil${\scriptstyle 2+5}$\hfil&\hfil${\scriptstyle (1,1)\overline{(1,1)}}$\hfil&\hfil${\scriptstyle 2 \times 1}$\hfil&\hfil${\scriptstyle {\rm two\ real}}$\hfil&\cr
{}&\hfil${\scriptstyle [(1\vert1\vert1)1];\ 3\mu+\sigma\not=0;\ \mu+3\sigma\not=0;\ \mu\not=\sigma}$\hfil&{}&{}&\hfil${\scriptstyle (1,1)(1,1)}$\hfil&\hfil${\scriptstyle 4 \times 1}$\hfil&\hfil${\scriptstyle {\rm nodes}^{19}}$\hfil\cr
\noalign{\vskip 2pt}
\noalign{\hrule}
\noalign{\vskip 2pt}
\hfil ${\scriptstyle {\rm IX}}$\hfil&\hfil${\scriptstyle \mu=\sigma}$\hfil&\hfil${\scriptstyle \lambda+\mu+\nu+\sigma=0}$\hfil&\hfil${\scriptstyle 2+5}$\hfil&\hfil${\scriptstyle (1,1)\overline{(1,1)}}$\hfil&\hfil${\scriptstyle 2 \times 1}$\hfil&\hfil${\scriptstyle {\rm two\ real}}$\hfil&\cr
{}&\hfil${\scriptstyle [1(1\vert1\vert1)];\ 3\lambda+\nu\not=0;\ \lambda+3\nu\not=0;\ \lambda\not=\nu}$\hfil&{}&{}&\hfil${\scriptstyle (1,1)(1,1)}$\hfil&\hfil${\scriptstyle 4 \times 1}$\hfil&\hfil${\scriptstyle {\rm nodes}^{19}}$\hfil\cr
\noalign{\vskip 2pt}
\noalign{\hrule}
\noalign{\vskip 2pt}
\hfil ${\scriptstyle {\rm IX}}$\hfil&\hfil${\scriptstyle \lambda=\mu\not=\nu=\sigma}$\hfil&\hfil${\scriptstyle \lambda+\mu+\nu+\sigma=0}$\hfil&\hfil${\scriptstyle 1+4}$\hfil&\hfil${\scriptstyle (1,0)\overline{(1,0)}}$\hfil&\hfil${\scriptstyle 2 \times 1}$\hfil&\hfil${\scriptstyle {\rm four\ complex}}$\hfil&\cr
{}&\hfil${\scriptstyle [(11)(11)];}$\hfil&{}&{}&\hfil${\scriptstyle (0,1)\overline{(0,1)}}$\hfil&{}&\hfil${\scriptstyle {\rm nodes}^{20}}$\hfil\cr
\noalign{\vskip 2pt}
\noalign{\hrule}
\noalign{\vskip 2pt}
\hfil ${\scriptstyle {\rm IX}}$\hfil&\hfil${\scriptstyle \lambda=\nu\not=\mu=\sigma}$\hfil&\hfil${\scriptstyle \lambda+\mu+\nu+\sigma=0}$\hfil&\hfil${\scriptstyle 1+4}$\hfil&\hfil${\scriptstyle (1,0)(1,0)}$\hfil&\hfil${\scriptstyle 2 \times 1}$\hfil&\hfil${\scriptstyle {\rm four\ real}}$\hfil&\cr
{}&\hfil${\scriptstyle [(1\vert(1\vert1)\vert1)];}$\hfil&{}&{}&\hfil${\scriptstyle (0,1)(0,1)}$\hfil&{}&\hfil${\scriptstyle {\rm nodes}^{21}}$\hfil\cr
\noalign{\vskip 2pt}
\noalign{\hrule}
\noalign{\vskip 2pt}
\hfil ${\scriptstyle {\rm IX}}$\hfil&\hfil${\scriptstyle \lambda=\sigma\not=\mu=\nu}$\hfil&\hfil${\scriptstyle \lambda+\mu+\nu+\sigma=0}$\hfil&\hfil${\scriptstyle 1+4}$\hfil&\hfil${\scriptstyle (1,0)(1,0)}$\hfil&\hfil${\scriptstyle 2 \times 1}$\hfil&\hfil${\scriptstyle {\rm four\ real}}$\hfil&\cr
{}&\hfil${\scriptstyle [(1\vert(11)\vert1)];}$\hfil&{}&{}&\hfil${\scriptstyle (0,1)(0,1)}$\hfil&{}&\hfil${\scriptstyle {\rm nodes}^{21}}$\hfil\cr
\noalign{\vskip 2pt}
\noalign{\hrule}
\noalign{\vskip 2pt}
\hfil ${\scriptstyle {\rm IX}}$\hfil&\hfil${\scriptstyle {\rm exactly\ three\ coincident\ eigenvalues}}$\hfil&\hfil${\scriptstyle \lambda+\mu+\nu+\sigma=0}$\hfil&\hfil${\scriptstyle 1+3}$\hfil&\hfil${\scriptstyle (1,1)^2}$\hfil&\hfil${\scriptstyle 2 \times 1}$\hfil&\hfil${\scriptstyle {\rm double\ projective}}$\hfil&\cr
{}&\hfil${\scriptstyle [(111)1)];\ [(1\vert 1 \vert11)];\ [(11\vert 1\vert 1)];\ [1(111)]}$\hfil&{}&{}&{}&{}&\hfil${\scriptstyle {\rm line}^{22}}$\hfil\cr
\noalign{\vskip 2pt}
\noalign{\hrule}
\noalign{\vskip 2pt}
\hfil ${\scriptstyle {\rm IX}}$\hfil&\hfil${\scriptstyle {\rm all\ 4\ eigenvalues\ coincide}}$\hfil&\hfil${\scriptstyle \lambda=\mu=\nu=\sigma=0}$\hfil&\hfil${\scriptstyle 0+0}$\hfil&\hfil${\scriptstyle \{-\}}$\hfil&\hfil${\scriptstyle 1 \times 1}$\hfil&\hfil${\scriptstyle {{\bf KP}^1 \times {\bf KP}^1}^{23}}$\hfil&\cr
{}&\hfil${\scriptstyle [(1111)]}$\hfil&{}&{}&{}&{}&{}\cr
\noalign{\vskip 2pt}
\noalign{\hrule}

\cr}}$$
$^1$ The linear factor in the Ricci polynomial defines a linear component in the Ricci locus, which is the tangent to the other component of the Ricci locus at the real null eigenvector and thus has infinite-order of contact with the locus, thus making the real null eigenvector a tacnode.\hfil\break
$^2$ There is a complex conjugate pair of null eigenvectors, each of which is a node; the linear factor in the Ricci polynomial defines a linear component of the Ricci locus, which is a common tangent to the two nodes. Its intersection with the other component yields the two nodes.\hfil\break
$^3$ The real Ricci locus is a simple closed curve on the torus with a cusp.\hfil\break
$^4$ The three components of the Ricci locus intersect at the null eigenvector to create the triple point.\hfil\break
$^5$ The real null eigenvector (of eigenvalue $\mu$) is a node with real tangents iff $(3\lambda+\nu)(\lambda+3\nu) > 0$ and a node with complex tangents (i.e., an isolated point of the real Ricci locus) iff $(3\lambda+\nu)(\lambda+3\nu)$. $(3\lambda+\nu)(\lambda+3\nu) = 0$ iff $\mu = \lambda$ or $\mu = \nu$, and thus does not occur in the generic scenario of type IV. The real Ricci locus is a closed curve on the torus with this single self-intersection.\hfil\break
$^6$ Ricci spinor type is $(1,1)(1,1)$/$(1,1)\overline{1,1)}$ according as $\epsilon = \pm 1$, where $s_{2,2}(v_2,v_3) = \epsilon$ and $\{v_2,v_3\}$ gives the JCF for $\Phi^a{}_b$ restricted to the middle summand of the decomposition for type IV. The null eigenvector is a point of intersection of the two components of the Ricci locus.\hfil\break
$^7$ Ricci spinor type is $(1,1)\overline{1,1)}$/$(1,1)(1,1)$/ according as $\epsilon = \pm 1$, where $s_{2,2}(v_2,v_3) = \epsilon$ and $\{v_2,v_3\}$ gives the JCF for $\Phi^a{}_b$ restricted to the middle summand of the decomposition for type IV. The null eigenvector is a point of intersection of the two components of the Ricci locus.\hfil\break
$^8$ The real null eigenvector of $\mu$ is the intersection of the two linear components of the Ricci locus; the two real null eigenvectors of $\lambda = \nu$ are each an intersection of a linear component and the quadratic component of the Ricci locus.\hfil\break
$^9$ The two distinct linear components of the Ricci locus intersect in a point corresponding to the real null eigenvector $E$ of the generic case of no coincident eigenvalues, and each component consists of further real null eigenvectors (each component is double, thus singular) determining a totally null eigensubspace, these two spaces intersecting in $E$.\hfil\break
$^{10}$ No null eigenvectors\hfil\break
$^{11}$ The two nodes are precisely the intersections of the two quadratic factors of the Ricci locus; they are the two null eigendirections of $\lambda=\nu$.\hfil\break
$^{12}$ The node is at the null eigendirection of $\lambda$. The real Ricci locus is a closed curve on the torus with this single self-intersection.\hfil\break
$^{13}$ The two nodes are the intersections of the linear component of the Ricci locus with the other component, and are real with real tangents; the linear component of the Ricci locus is in fact a common tangent to the two nodes.\hfil\break
$^{14}$ In the decomposition for this type, let $\{v_1,v_2\}$ be the basis for the first summand giving the JCF on that summand and such that $\{v_1,\epsilon v_2\}$ is a Witt basis for that summand, and let $\{v_3,v_4\}$ be the basis for the second summand giving the JCF on that summand and such that $\{v_3,\omega v_4\}$ is a Witt basis for that summand. Then, for those $\Phi^a{}_b$ that take the type's matrix form with respect to Witt bases of O-type ${\bf SO}$, the Ricci spinor type is $(1,0)(1,0)(0,1)^2$ if $\epsilon = -\omega$ and is of type $(1,0)\overline{(1,0)}(0,1)^2$ if $\epsilon = \omega$; for those $\Phi^a{}_b$ that take the type's matrix form with respect to Witt bases of O-type ${\bf ASO}$, the Ricci spinor type is $(1,0)^2(0,1)(0,1)$ if $\epsilon = -\omega$ and is of type $(1,0)^2(0,1)\overline{(0,1)}$ if $\epsilon = \omega$. The double component corresponds to the totally null eigensubspace spanned by the two real null eigenvectors $v_1$ and $v_3$ (with common eigenvalue zero) and intersects the other two components of the Ricci locus; these other two components may be taken to be real when $\epsilon = -\omega$ but are complex when $\epsilon=\omega$. Technically, these points of intersection are tacnodes with a common tangent; they are not the points $v_1$ and $v_3$.\hfil\break
$^{15}$ No null eigenvectors\hfil\break
$^{16}$ The Ricci polynomial consists of four linear factors, each of which defines a nonsingular component of the Ricci locus. These four components define four points of intersection, which are complex nodes (i.e., lie on $\Omega\setminus\omega$), with the two components that intersect at a node providing the tangents at that node. Note that the coincidence $\lambda = \bar\mu$ is not a distinct case but just a relabelling of the pair of complex conjugate eigenvalues $(\mu,\bar\mu)$.\hfil\break
$^{17}$ No null eigenvectors.\hfil\break
$^{18}$ The two nodes are the points of intersection of the two components of the Ricci locus. Let $a$ and $b$ denote the two eigenvalues not involved in the coincidence. For both null eigenvectors, the discriminant defining the character of the null eigenvector as a double point is $3a^2 + 10ab + 3b^2 = (3a+b)(a+3b)$; moreover, when $3a^2 + 10ab + 3b^2 = (3a+b)(a+3b) > 0$, the Ricci spinor type is $(1,1)\overline{(1,1)}$ and when $3a^2 + 10ab + 3b^2 = (3a+b)(a+3b) < 0$ the Ricci spinor type is $(1,1)(1,1)$. \hfil\break
$^{19}$ The two nodes are the points of intersection of the two components of the Ricci locus. Let $a$ and $b$ be as in footnote 17. The nodes are real, they both have real tangents iff $3a^2 + 10ab + 3b^2 = (3a+b)(a+3b) > 0$ and this is the case of Ricci spinor type $(1,1)(1,1)$. When $3a^2 + 10ab + 3b^2 = (3a+b)(a+3b) < 0$, the tangents are complex, i.e., the nodes are isolated points of the real Ricci locus, and the Ricci spinor type is $(1,1)\overline{(1,1)}$.\hfil\break
$^{20}$ The Ricci locus $\Omega$ is four distinct complex projective lines that have four points of intersection. These four complex nodes are the four complex null eigenvectors; the Ricci spinor type is $(1,0)\overline{(1,0)}(0,1)\overline{(0,1)}$.\hfil\break
$^{21}$ In these two cases, geometrically equivalent at a certain level, it suffices to consider the real Ricci locus, which consists of four distinct real projective lines that have four points of intersection. These four real nodes are the four real null eigenvectors; the spinor Ricci type is $(1,0)(1,)(0,1)(0,1)$.\hfil\break
$^{22}$ In each case, the coincident eigenvalues have (real) eigenspace isomorphic to ${\bf R}^{21,}$ or ${\bf R}^{1,2}$; the generators of the null cone in the relevant eigenspace form a real projective line; the real Ricci locus is the square of this line.\hfil\break
$^{23}$ The zero Ricci spinor is the only case.

\vskip 24pt
\noindent {\section Bibliography}
\vskip 24pt
\frenchspacing
\vskip 1pt
\hangindent=20pt \hangafter=1
\noindent Beem, J. K., Ehrlich, P. E. \& Easley, K. L. 1966 {\sl Global Lorentzian Geometry}. Second Edition. Marcel Dekker, New York.
\vskip 1pt
\hangindent=20pt \hangafter=1
\noindent Carter, R., Segal, G. \& MacDonald, I. 1995 {\sl Lectures on Lie Groups and Lie Algebras}. London Mathematical Society Student Text {\bf 32}, Cambridge University Press, Cambridge.
\vskip 1pt
\hangindent=20pt \hangafter=1
\noindent Churchill, R. V. 1932 Canonical Forms for Symmetric linear Vector Functions in Pseudo-Euclidean Space. {\sl Trans. Amer. Math. Soc.} {\bf 34}, 784--794.
\vskip 1pt
\hangindent=20pt \hangafter=1
\noindent Derdzinski, A. 2000 Einstein Metrics in dimension four; in {\sl Handbook of Differential Geometry}, Vol. I, F. J. E. Dillen \& L. C. A. Verstraelen (eds), N-H Elsevier, Amsterdam, pp. 419--707.
\vskip 1pt
\hangindent=20pt \hangafter=1
\noindent Kramer, D., Stephani, H., MacCallum, M. \& Herlt, E. 1980 {\sl Exact solutions of Einstein's field equations}. Cambridge University Press, Cambridge.
\vskip 1pt
\hangindent=20pt \hangafter=1
\noindent Law, P. R. 1991 Neutral Einstein metrics in four dimensions. {\sl J. Math. Phys.\/} {\bf 32}, 3039--3042.
\vskip 1pt
\hangindent=20pt \hangafter=1
\noindent Law, P. R. 1992 Neutral Geometry and the Gauss-Bonnet Theorem for Two-Dimensional Pseudo-Riemannian Manifolds. {\sl Rocky Mountain J. Math.\/} {\bf 22}, 1365--1383.
\vskip 1pt
\hangindent=20pt \hangafter=1
\noindent Law, P. R. 2006 Classification of the Weyl curvature spinors of neutral metrics in four dimensions. {\sl J. Geo. Phys.} {\bf 56}, 2093--2108.
\vskip 1pt
\hangindent=20pt \hangafter=1
\noindent Law, P. R. 2009 Spin Coefficients for four-dimensional neutral metrics, and null Geometry. {\sl J. Geo. Phys.} {\bf 59(8)}, 1087--1126. arXiv:0802.1761v2 [math.DG] (26 Aug 2009).
\vskip 1pt
\hangindent=20pt \hangafter=1
\noindent Law, P. R. \& Matsushita, Y. 2008 A Spinor Approach to Walker Geometry. {\sl Communications in Mathematical Physics} {\bf 282}, 577-623. (arXiv:math/0612804v4 [math.DG] 7 Apr 2009.)
\vskip 1pt
\hangindent=20pt \hangafter=1
\noindent O'Neill, B. 1983 {\sl Semi-Riemannian Geometry With Applications to Relativity}. Academic Press, Orlando, FL.
\vskip 1pt
\hangindent=20pt \hangafter=1
\noindent Penrose, R. \& Rindler, W. 1984 {\sl Spinors and Space-Time, Vol. 1: Two-spinor calculus and relativistic fields\/}, Cambridge University Press, Cambridge.
\vskip 1pt
\hangindent=20pt \hangafter=1
\noindent Penrose, R. \& Rindler, W. 1986 {\sl Spinors and Space-Time, Vol. 2: Spinor and twistor methods in space-time geometry\/}, Cambridge University Press, Cambridge.
\vskip 1pt
\hangindent=20pt \hangafter=1
\noindent \noindent Petrov, A. Z. 1969 {\sl Einstein Spaces}. Pergamon Press, Oxford.
\vskip 1pt
\hangindent=20pt \hangafter=1
\noindent Porteous, I. R. 1981 {\sl Topological Geometry\/}, 2nd Edition, Cambridge University Press, Cambridge.
\vskip 1pt
\hangindent=20pt \hangafter=1
\noindent Porteous, I. R. 1995 {\sl Cifford Algebras and the Classical Groups}. Cambridge University Press, Cambridge.
\vskip 1pt
\hangindent=20pt \hangafter=1
\noindent Walker, R. J. 1962 {\sl Algebraic Curves}. Dover, New York.
\vskip 1pt
\hangindent=20pt \hangafter=1
\noindent Warner, F. W. 1983 {\sl Foundations of Differentiable Manifolds and Lie Groups}. GTM 94, Springer-Verlag, New York.
\vskip 1pt

\bye